\documentclass[12pt,reqno]{amsart}
\usepackage[headings]{fullpage}
\usepackage{amssymb,amsmath,amscd,bbm,mathrsfs,extarrows}
\usepackage{graphicx}
\usepackage{texdraw}
\usepackage{pb-diagram}
\usepackage[all,cmtip]{xy}
\usepackage{url}
\usepackage[bookmarks=true,%
    colorlinks=true,%
    linkcolor=blue,%
    citecolor=blue,%
    filecolor=blue,%
    menucolor=blue,%
    urlcolor=blue,%
    breaklinks=true]{hyperref}
\usepackage{slashed}    
\usepackage{tikz-cd}
\usepackage{stmaryrd}
\usepackage{cancel,mathtools}
\usepackage{standalone,tikz,tikz-cd}
\usetikzlibrary{positioning,knots,patterns,hobby,decorations.pathreplacing}

\usepackage{soul,xcolor,fixmath}

\usepackage{verbatim}

\newcommand{\no}[1]{}

\no{
\usepackage {autobreak}
\usepackage{amsthm}
\usepackage{latexsym}
\usepackage{here}
\usepackage{setspace}
\usepackage{caption}
\usepackage{amscd}
\usepackage{extarrows}
\usepackage{mathrsfs}
\usepackage{multicol}
\graphicspath{{tex2/}}
\def\qed{\hfill $\Box$} 
\usepackage[dvipsnames]{xcolor}
\usepackage{mathrsfs}

\usepackage{comment}
\allowdisplaybreaks
\usepackage{accents}
\usepackage{extarrows}

}

\usepackage{soul,xcolor}

\setstcolor{red}

\DeclareMathOperator{\tr}{tr}

\DeclareMathOperator{\Fr}{Fr}
\DeclareMathOperator{\Gr}{Gr}
\DeclareMathOperator{\MaxSpec}{MaxSpec}

\DeclareMathOperator{\Hom}{Hom}
\DeclareMathOperator{\ord}{ord}

\DeclareMathOperator{\id}{id}
\DeclareMathOperator{\pr}{pr}

\newcommand{\qeq}{{\,\, \stackrel{(q)}{=}\, \,}}

\usepackage[headings]{fullpage}
\usepackage{amssymb,epic,eepic,epsfig,amsbsy,amsmath,amscd,color}
\usepackage[all]{xy}
\usepackage{graphicx,subcaption}
\usepackage{texdraw}
\usepackage{url}
\usepackage{bbm}
\usepackage{mathrsfs}
\usepackage{standalone}
\usepackage{tikz}
\usetikzlibrary{cd,positioning,knots,patterns,hobby}
\usepackage{accents}
\usepackage[normalem]{ulem}
\usepackage{comment}

\definecolor{ForestGreen}{RGB}{34,139,34}

\newcommand{\red}[1]{{\color{red}#1}}
\newcommand{\orange}[1]{{\color{orange}#1}}

\newtheorem{theorem}{Theorem}[section]
\newtheorem{THM}{Theorem}
\newtheorem{lemma}[theorem]{Lemma}
\newtheorem{cor}[theorem]{Corollary}
\newtheorem{corrr}[THM]{Corollary}

\newtheorem{proposition}[theorem]{Proposition}
\newtheorem{conjecture}{Conjecture}
\newtheorem{question}{Question}

\newtheorem{definition}[theorem]{Definition}
\theoremstyle{definition}
\newtheorem{remark}[theorem]{Remark}
\newtheorem{example}[theorem]{Example}

\newcommand{\bcon}{\begin{conjecture}}
\newcommand{\econ}{\end{conjecture}}
\newcommand{\bcor}{\begin{cor}}
\newcommand{\ecor}{\end{cor}}
\newcommand{\bdf}{\begin{definition}}
\newcommand{\edf}{\end{definition}}
\newcommand{\benu}{\begin{enumerate}}
\newcommand{\eenu}{\end{enumerate}}
\newcommand{\beq}{\begin{equation}}
\newcommand{\eeq}{\end{equation}}
\newcommand{\bexa}{\begin{example}}
\newcommand{\eexa}{\end{example}}
\newcommand{\bexe}{\begin{exercise}}
\newcommand{\eexe}{\end{exercise}}
\newcommand{\bfac}{\begin{fact}}
\newcommand{\efac}{\end{fact}}
\newcommand{\bite}{\begin{itemize}}
\newcommand{\eite}{\end{itemize}}
\newcommand{\blem}{\begin{lemma}}
\newcommand{\elem}{\end{lemma}}
\newcommand{\bmat}{\begin{pmatrix}}
\newcommand{\emat}{\end{pmatrix}}
\newcommand{\bprb}{\begin{problem}}
\newcommand{\eprb}{\end{problem}}
\newcommand{\bpro}{\begin{proposition}}
\newcommand{\epro}{\end{proposition}}

\newcommand{\bque}{\begin{question}}
\newcommand{\eque}{\end{question}}
\newcommand{\brem}{\begin{remark}}
\newcommand{\erem}{\end{remark}}
\newcommand{\bthm}{\begin{theorem}}
\newcommand{\ethm}{\end{theorem}}

\newcommand{\bpr}{\begin{proof}}
\newcommand{\epr}{\end{proof}}
\newcommand{\ignore}[1]{}
\newcommand{\xra}{\xrightarrow}

\def\onto{\twoheadrightarrow}

\newcommand{\cC}{\mathcal{C}}
\newcommand{\bm}{\mathbf{m}}

\def\BC{\mathbb C}
\def\BN{\mathbb N}
\def\BZ{\mathbb Z}
\def\BR{\mathbb R}

\def\cA{\mathscr A}

\def\ft{\mathfrak t}

\def\la{\langle}
\def\ra{\rangle}

\def\cM{\mathcal M}

\def\tx{\tilde x}
\def\tQ{{\tilde Q}}

\def\cS{{\mathcal S}}
\def\ot{\otimes}

\def\cE{\mathcal E}
\def\cF{\mathcal F}

\def\Mat{\mathrm{Mat}}
\def\bk{{\mathbf k}}
\def\bn{{\mathbf n}}

\def\Id{\mathrm{Id}}
\def\fS{\mathfrak S}
\def\bfS{\overline{\fS}}
\def\pfS{\partial \fS}

\def\cY{\mathcal Y}
\def\ev{{\mathrm{ev}}}

\def\embed{\hookrightarrow}

\def\sX{\mathfrak X}

\newcommand{\al}{\alpha}

\newcommand{\fm}{\mathfrak m}

 \def\bm{\mathbf{m}}

 \def\cC{\mathcal C}

\def\bl{{\mathbf{l}} }

\def\tr{\mathrm{tr}}

\def\utr{{\underline{\tr}}}
\def\bbl{{\mathbf l}}
\def\mP{{\mathsf P}}
\def\uD{{\underline{\Delta}}}
\def\uA{{\underline{\cA}}}
\def\uP{{\underline{\mP}}}
\def\uv{\overleftarrow }
\def\pbfS{{\partial\bfS}}

\def\al{\alpha}
\def\ve{\varepsilon}
\def\be { \begin{equation} }
\def\ee { \end{equation} }

\def\bD{{\bar \Delta }}

\def\bt{ {\mathbf t}}

\def\CC{\mathcal C}

\def\bT{\mathbb T}
\def\CC{\mathcal C}

\def\lt{\mathrm{lt}}

\def\nc{\newcommand}
\def\Fr{\mathrm{Fr}}
\def\Weyl{{\mathrm{Weyl}}}

                  \nc\FI[2]{\begin{figure}
    \begin{center}\input{#1.pstex_t}\end{center}
    \caption{#2}
    \lbl{#1}
  \end{figure}}
\nc\FIG[3]{\begin{figure}
    \includegraphics[#3]{#1.eps}
    \caption{#2}
    \lbl{fig:#1}
    \end{figure}}
\nc\FF[3]{\begin{figure}
    \includegraphics[#3]{#1.eps}
    \caption{#2}
    \lbl{#1}
    \end{figure}}
    \nc\FIGc[3]{\begin{figure}[htpb]
    \includegraphics[height=#3]{draws/#1.eps}
    \caption{#2}
    \label{fig:#1}
    \end{figure}}
    
    \nc\FIGjpg[3]{\begin{figure}[htpb]
    \includegraphics[height=#3]{draws/#1.jpg}
    \caption{#2}
    \label{fig:#1}
    \end{figure}}

    \nc\FIGh[3]{\begin{figure}[htpb]
    \includegraphics[height=#3]{#1.eps}
    \caption{#2}
    \lbl{fig:#1}
    \end{figure}}

\def\ord{\mathrm{ord}}

     \def\tr{\mathrm{tr}}

\def\lt{\mathrm{lt}}

\def\cR{{\mathcal R}}
\def\BC{{\mathbb C}}

\def\G{{SL_2(\BC)}}
\def\X{{\mathcal X}}
\def\vz{{\vec z}}

\def\w{{\mathsf w}}
\def\vw{{\vec{\mathsf w}}}

\newcommand\incl[2]{{\includegraphics[height=#1]{draws/#2.eps}}}

\title
{Degenerations on sliced skein algebras}

\title[Geometric Kauffman bracket]
{ Sliced skein algebras and geometric Kauffman bracket}

\author{Charles D. Frohman}

\address{Department of Mathematics, The University of Iowa, Iowa City, IA
52242, USA}

\email{\tt charles-frohman@uiowa.edu}

\author{Joanna Kania-Bartoszynska}

\address{National Science Foundation, Arlington, VA, 22230, USA} 

\email{\tt jkaniaba@nsf.gov}

\author{ Thang T.Q. L\^{e}}

\address{Georgia Tech}

\email{\tt letu@math.gatech.edu}

    \def\RS{{ \cR_\fS }}
    \def\SS{{  \mathscr{S}(\fS)}}

\def\Zq{{\BZ[q^{\pm 1/2}]}}

\def\qq{{q^{1/2}}}

\def\term#1{{\bf #1}}
\def\eqdef{:=}
\def\tfS{ { \widetilde{\fS}  }}

\def\ptfS{ {\partial { \tfS }}}
\def\pal{{\partial \al}}

\def\ori{{\mathfrak o}}

\def\SSR{{\cS_\xi(\fS;\cR)}}
\def\SslR{{\cS^{\sli, \vw}_\xi(\fS;\cR)}}
\def\sli{{\mathsf{sl}}}

\def\Mo{{\mathring{\cM}}}
\def\ddd{\mathbbm{d}}
\def\bnu{{\boldsymbol{\nu}}}
\def\pv{{\fd _v}}
\def\SslS{\cS^\sli(\fS)}

\def\cI{{\mathcal I}}

\def\RR{{\mathsf R}}

\def\fd{{\mathfrak{d}}}
\def\BslS{\B^\sli(\fS)}

\def\X{\chi}
\def\B{{\mathsf{B}}}

\def\Res{\mathrm{Res}}
\def\Cr{{\mathsf{Cr}}}
\def\XS{\X(\fS)}
\def\XSe{{ \X^\vw(\fS)}}
\def\Sx{\cS_\xi}
\def\Sslx{{\cS_\xi^{\sli, \vw  }(\fS;\cR)}}
\def\SxS{\Sx(\fS)}
\def\bi{{\mathbf i}}
\def\csh{{\mathsf{csh}}}
\def\SeS{{\cS_\ve(\fS)}}
\def\Ssl{\cS^{\sli}}
\def\PIdeg{{ \text{\rm{PI-deg }}}}

\newcommand{\towrite}[1]{{ \red{to write} }}

\begin{document}

\begin{abstract} The sliced skein algebra of a closed surface of genus $g$ with $m$ punctures, $\fS=\Sigma_{g,m}$,  is the quotient of the Kauffman bracket skein algebra $\SxS$ corresponding to fixing the  scalar values of its  peripheral curves. We show that the sliced skein algebra of a finite type surface is  a domain if the ground ring is a domain. When the quantum parameter $\xi$ is a root of unity we calculate the center of the sliced skein algebra and its PI-degree. Among applications we show that any smooth point of a sliced character variety is  an Azumaya point of the skein algebra $\SxS$.

For any $SL_2(\BC)$--representation $\rho$ of the fundamental group of an oriented connected 3-manifold $M$ and a root  of unity $\xi$ with the order of $\xi^2$ odd, we introduce the $\rho$-reduced skein module $\cS_{\xi, \rho}(M)$. We show that $\cS_{\xi, \rho}(M)$ has dimension 1  
 when $M$ is  closed and $\rho$ is irreducible. 
We also show that if $\rho$ is irreducible the $\rho$-reduced skein module of a handlebody, as a module over  the skein algebra of its boundary, is simple and has the dimension equal to the PI-degree of the skein algebra of its boundary.

\end{abstract}

\maketitle

\section{Introduction}

Let $\fS=\Sigma_{g,m}=\Sigma_g \setminus \{v_1,\dots, v_m\}$   where $\Sigma_g$ is an oriented closed surface of genus $g$. We assume $\fS$ has negative Euler characteristic. 
Let the ground ring $\mathcal{R}$ be a commutative domain with a distinguished invertible element $\xi$.

\subsection{Integrality of the  sliced  skein algebra}

The Kauffman bracket skein algebra $\SSR$, introduced by Przytycki and Turaev \cite{Prz,Turaev}, is the  $\cR$-module freely spanned by isotopy classes of link diagrams on $\fS$ subject to the Kauffman bracket relations \cite{Kauffman}

\begin{align}
\label{eqSkein0}
\begin{tikzpicture}[scale=0.8,baseline=0.3cm]
\fill[gray!20!white] (-0.1,0)rectangle(1.1,1);
\begin{knot}[clip width=8,background color=gray!20!white]
\strand[very thick] (1,1)--(0,0);
\strand[very thick] (0,1)--(1,0);
\end{knot}
\end{tikzpicture}
&=\xi
\begin{tikzpicture}[scale=0.8,baseline=0.3cm]
\fill[gray!20!white] (-0.1,0)rectangle(1.1,1);
\draw[very thick] (0,0)..controls (0.5,0.5)..(0,1);
\draw[very thick] (1,0)..controls (0.5,0.5)..(1,1);
\end{tikzpicture}
+\xi^{-1}
\begin{tikzpicture}[scale=0.8,baseline=0.3cm]
\fill[gray!20!white] (-0.1,0)rectangle(1.1,1);
\draw[very thick] (0,0)..controls (0.5,0.5)..(1,0);
\draw[very thick] (0,1)..controls (0.5,0.5)..(1,1);
\end{tikzpicture}\, , \quad 
\begin{tikzpicture}[scale=0.8,baseline=0.3cm]
\fill[gray!20!white] (0,0)rectangle(1,1);
\draw[very thick] (0.5,0.5)circle(0.3);
\end{tikzpicture}
=(-\xi^2 -\xi^{-2})
\begin{tikzpicture}[scale=0.8,baseline=0.3cm]
\fill[gray!20!white] (0,0)rectangle(1,1);
\end{tikzpicture}\, . 
\end{align}

\def\SxC{\cS_\xi(\fS; \BC)}

The product is given by stacking. See Section \ref{skein} for details.

The skein algebra  and its analogs play an important role in low dimensional topology and quantum algebra.  The skein algebra $\cS_{-1}(\fS;\BC)$ is canonically isomorphic to the coordinate ring of the $SL_2(\BC)$-character variety $\XS$ of the surface \cite{Bullock,PS1,CM}. 

Understanding the representation theory of $\SxC$ for a root of unity $\xi\in \BC$ is a step in the construction of geometric field theories associated to the Jones polynomial. The first step in analyzing the representation theory of an algebra is to identify its center. In our case, the center of $\SxC$ is almost the same as $\BC[\XS]$, the ring of regular functions on $\XS$. 
The work of \cite{GJS} introduced the use of techniques from symplectic geometry into the study of the representation theory of skein algebras.   The character variety $\XS$ has the Atiyah-Bott-Goldman  Poisson structure  and is stratified by sliced character varieties $\XSe$, each of which is a symplectic leaf  with singularities.
Here $\vw=(\w_1, \dots, \w_m) \in \BC^m$ and $\XSe\subset \XS$ is the subset of characters taking value $\w_i$ on the peripheral loop $\fd_i$, a small loop surrounding the puncture $v_i$.  For details and references see  Section \ref{secMRY}. 

In order to study $\SSR$ we introduce the sliced skein algebra $\SslR$, 
$$ \SslR=  \SSR/( \fd _i - \w_i).$$
When $R=\BC$ and $\xi=-1$ the sliced skein algebra is isomorphic to $\BC[\XSe]$.
 
As most results in representation theory concern prime algebras, we want to know if the sliced skein algebra is prime. We prove a stronger result:
\begin{THM}[See Theorem \ref{thmMain2}] \label{thmDomain1} Let  $\cR$ be a commutative domain,  $\xi\in \cR$ be invertible,  $\vw\in \cR^m$, and $\fS=\Sigma_{g,m}$. 
The  sliced skein algebra $\SslR$ is a domain, meaning if $ab =0$ where $a,b\in \SslR$ then either $a=0$ or $b=0$.
\end{THM}
In particular, for $\cR=\BC$ and $\xi=-1$, we get
\begin{corrr} For $\vw\in \BC^m$, the sliced character variety $\XSe$ is an irreducible affine  variety.
\end{corrr}
The corollary was recently proved by P. Whang \cite{Whang} by a different, quite non-trivial proof. 
Our proof is based on ideas coming from the non-commutative nature of the skein algebras.

Actually we  prove the following stronger result, which implies  Theorem \ref{thmDomain1}. It is not difficult to see that the set $\BslS$ of all isotopy classes of link diagrams not containing (i) crossings, (ii) contractible loops, and (iii) peripheral loops $\fd_i$, is a free $\cR$-basis of $\SslR$.

\begin{THM}[See Theorem \ref{thmMain2}]

\label{thm12} Assume  $\fS=\Sigma_{g,m}$ has
negative Euler characteristic. Let $r= 3-3g+m$. 
There exist a submonoid $\Lambda\subset \BZ^{2r}$, an antisymmetric $2r\times 2r$ integral matrix $\tQ$,  and a bijection $\mu: \Lambda \xra{\cong} \BslS$ such that in $\SslR$,
\be 
\mu(\bn) \mu(\bm) = \xi^{ \frac12\la \bn, \bm\ra_\tQ} \mu(\bn + \bm) + F_{  < \bn+\bm}, \label{Topdeg0}
\ee
where $\la \bn, \bm\ra_\tQ := \sum_{i,j} \tQ_{ij} n_i m_j $ is always even, and
\begin{align*}
F_{  < \bn+\bm} &:= \cR\text{-span of } \ \{ \mu(\bk) \mid \bk < \bn + \bm \ \text{in lexicographic order}\}.
\end{align*}
\end{THM}

The theorem shows that there is a filtration on $\Sslx$ whose associated graded algebra is a monomial subalgebra of a quantum torus.

The form $\tQ$ is defined combinatorially  from a  pants decomposition of $\fS$, together with its dual graph. See Section \ref{secDegen}.

\subsection{The role of the center and the PI-degree}

Let us recall some representation theory. 
Let $A$ be an affine $\BC$-algebra with center $Z$. Assume that $A$ is a domain, and as a $Z$-module $A$ is finitely generated. 
Let $\Fr(Z)$ be the field of fractions of $Z$. By Wedderburn's theorem $A \ot_Z \Fr(Z)$ is a  
division algebra, and consequently has dimension $D^2$ over $\Fr(Z)$, where $D$ is a positive integer called the PI-degree of $A$. By Artin-Tate's lemma $Z$ is an affine $\BC$-algebra. Let $V= \MaxSpec(Z)$, considered as an affine variety. A point $\fm \in V$ is called Azumaya if $A/ \fm A\cong M_D(\BC)$, the algebra of $D \times D$ matrices with entries in $\BC$. The {\bf Azumaya locus} $\cA(A)\subset V$ of all Azumaya points is known  to be Zariski open dense in $V$. Any irreducible representation $\rho$ has its central character $\csh(\rho)\in V$, and if $\csh(\rho)$ is not Azumaya then the dimension of $\rho$ is less than $  D$. On the other hand if $\csh(\rho)\in \cA(A)$ then it has dimension $D$, and if $\csh(\rho)= \csh(\rho')$ then $\rho$ and $\rho'$ are conjugate. We call $\csh(\rho)$ the classical shadow of $\rho$, and $V= \MaxSpec(Z)$ the  variety of shadows of $A$.

\subsection{The center and  the PI-degree of the sliced skein algebra} Until the end of the introduction since $\cR=\BC$ we will drop $\BC$ in the notation of $\SxC$.

Let $\xi\in \BC$ be a root of 1.  In \cite{FKL,FKL2} we calculated the center of $\SxS$ and its PI-degree, based on earlier work of Bonahon and Wong \cite{BW2}.

Let $N= \ord(\xi^4)$. Then $\ve:= \xi^{N^2}\in \{ \pm1, \pm \bi\}$, where $\bi$ is the complex unit. Bonahon and Wong constructed an algebra embedding $\Phi_\xi: \SeS \embed \SxS$, called the Chebyshev-Frobenius map, which is recalled in Section \ref{skein}. It is not difficult to show that $\Phi_\xi$ descends to an algebra embedding 
$$\Phi_\xi^\vw: \cS_\ve^{\sli, T_N(\vw)}(\fS)\embed \cS_\xi^{\sli, \vw}(\fS).$$
 Here $T_N$ is the $N$-th Chebyshev polynomial of type 1 and $T_N(\vw)=( T_N(\w_1), \dots, T_N(\w_m))$. 

Using the top degree formula \eqref{Topdeg0} and the nature of the form $\tQ$ we prove
\begin{THM}[See Theorem \ref{thmCenter}] \label{thmCenter00}
Let $\xi\in \BC$ be a root of 1, $\fS= \Sigma_{g,m}$, and $\vw \in \BC^m$. Let $N=\ord(\xi^4)$ and 
 $r=3g-3+m$. 

(a)  The center of $Z(\cS_\xi^{\sli, \vw}(\fS))$ is the image of 
$\Phi_\xi^\vw$ or the even image of $\Phi_\xi^\vw$ according as $\ord(\xi^2)$ is odd or even. In other words,
\be 
Z(\cS_\xi^{\sli, \vw}(\fS)) = \begin{cases}  \Phi_\xi^\vw \left( \cS_\ve^{\sli, T_N(\vw)}(\fS) \right) & \ord(\xi^2) \not \in 2 \BZ \\
 \Phi_\xi^\vw \left( \cS_\ve^{\sli, T_N(\vw)}(\fS)^\ev \right)  & \ord(\xi^2) \in 2 \BZ 
\end{cases}
\ee
where $\cS_\ve^{\sli, T_N(\vw)}(\fS)^\ev$ is the $\BC$-subalgebra generated by link diagrams having even intersection number with any loop on $\fS$.

(b) The sliced skein algebra  $\cS_\xi^{\sli, \vw}(\fS)$  is a $\BC$-domain, is finitely generated over its center, and has PI-degree equal to that of $\SxS$, given by

\be
\PIdeg (\cS_\xi^{\sli, \vw}(\fS)) = \PIdeg (\SxS) = \begin{cases} N^{r} & \ord(\xi^2) \not \in 2 \BZ \\
2^{g} N^{r} & \ord(\xi^2) \in 2 \BZ.
\end{cases}
\label{eqPIdeg0}
\ee

\end{THM}

\def\bu{{\bullet}}
\def\smooth{{\mathrm {sm}}}
\subsection{Azumaya loci of skein algebras} Theorem \ref{thmCenter00} allows to use the theory of Poisson orders \cite{BG} to analyze the Azumaya locus of $\SxS$, by focusing on the quotients $\cS_\xi^{\sli, \vw}(\fS)$.

Let $\xi$ be root of 1, with $N= \ord(\xi^4)$. Let $\fS=\Sigma_{g,m}$ and $V_\xi: =\MaxSpec(Z(\SxS))$ be the classical shadow variety of $\SxS$.
By \cite{FKL}  there is finite morphism  of degree $N^{m}$
\be 
p_\xi: V_\xi \onto \X(\fS)^\bullet,
\ee
where $\X(\fS)^\bu=\X(\fS)$ if $\ord(\xi^2)$ is odd; otherwise, 
\be 
\X(\fS)^\bu= \X(\fS)^\ev:= \X(\fS)// H^1(\bar\fS, \BZ/2).
\ee
Here $\bar \fS= \Sigma_g$, and the action of $u \in H^1(\bar \fS,\BZ/2)\equiv \Hom(\pi_1(\bar\fS), \{\pm 1\})$ on $\X(\fS)$ is given by $[\rho] \to [u \rho]$, where $(u\rho)(\al)= u(\al) \rho(\al)$ for $\al\in \pi_1(\fS)$. We call $\X(\fS)^\ev$ the even character variety of $\fS$, which is closely related to the $PSL_2(\BC)$ character variety; see Section \ref{secChar}.

For $\vw=(\w_1, \dots, \w_m)\in \BC^m $ we define the {\bf sliced character variety} $\X^{\vw}(\fS)^\bu\subset \X(\fS)^\bu$ to be the set of characters $[\rho]$ for which $\tr(\rho(\fd_i)) = \w_i$. We have the stratification 
\be 
\X(\fS)^\bu = \bigsqcup_{\vw \in \BC^m } \X^{\vw}(\fS)^\bu.
\ee

\begin{THM}[ Theorem \ref{thmSliAzu} ] \label{thmAzu0}
Let $\xi\in \BC$ be a root of 1 and $\fS=\Sigma_{g,m}$. If $x$ is a smooth point of a sliced character variety $ \X^{\vw}(\fS)^\bu$ for some $\vw\in \BC^m$, then  the fiber $(p_\xi)^{-1}(x)$ is in the Azumaya locus of $\SxS$. In other words, if $ (\X^{\vw}(\fS)^\bu) ^\smooth$ is the smooth locus of $ \X^{\vw}(\fS)^\bu$, then
\be 
(p_\xi)^{-1} \left (  \bigsqcup_{\vw\in \BC^m} ( \X^{\vw}(\fS)^\bu) ^\smooth \right) \subset \cA(\SxS).
\ee
\end{THM}

When $m=0$, i.e.,  $\fS$ is a closed surface without punctures, and at the same time $\ord(\xi)$ is odd, the theorem was proved by Ganev-Jordan-Safronov \cite{GJS}, using the theory of Poisson orders \cite{BG}, which is also used here. The case when $\ord(\xi^2)$ is even, i.e. when we  deal   with $\X(\fS)^\ev$, is technically more complicated.

In light of Theorem \ref{thmAzu0}, we want to know when a point of  $ \XSe$ is smooth. By Proposition \ref{rSliceSmooth}, a smooth point of $ \XSe$ is smooth as a point of $\XS$.
There is a generic condition of $\vw$ which guarantees that the whole sliced character variety $\XSe$ is smooth, see Proposition \ref{rKostov}. When $\rho: \pi_1(\Sigma_{0,g}) \to SL_2(\BC)$ is irreducible and $\tr(\rho(\fd_i)) \neq \pm 2$ for all $i=1, \dots, g$, we show in  Proposition \ref{rSmooth} that $[\rho]$ is a smooth point of the sliced character variety containing it. This fact will be used in the proof of Theorem \ref{thmGeoKau}.

\brem For results concerning the abelian characters see \cite{KK}. 

\erem

\def\G{SL_2(\BC)}
\def\SeM{\cS_\ve(M)}
\subsection{Geometric skein modules} Let $\xi\in \BC$ be a root of unity. Assume $\ord(\xi^2)$ is odd. 

When $m=0$, the surface $\fS= \Sigma_g$ is closed and  $V_\xi= \MaxSpec(Z(\SxS))= \X(\fS)$. 
  It is known that every irreducible character $x$ is a smooth point of $\X$, and hence  is Azumaya. Thus up to equivalence there is a unique irreducible representation of $\SxS$ having $x$ as its classical shadow, called the Azumaya representation above $x$. Its dimension is
  the PI-degree of $\SxS$, which is $N^{3g-3}$, where $N= \ord(\xi^4)$. 
  We show that in some cases there is a geometric realization of the Azumaya representation. 

First, for an oriented 3-manifold $M$, define the skein module $\Sx(M)$ using the same relations as in~\eqref{eqSkein0}. There is an natural action of $\Sx(\partial M)$ on $\Sx(M)$ obtained by stacking.

The Chebyshev-Frobenius map can be defined for 3-manifolds, and it gives rise to an action of $\SeM$ on $\Sx(M)$; see Section \ref{skein}. If $\rho:\pi_1(M)\to SL_2(\BC)$ is a representation, then its character is a point of $\X(M)$, and consequently determines a maximal ideal $\fm_\rho$ of $\SeM$. The {\bf skein module of $M$ reduced by $\rho$} is
 \be 
 \cS_{\xi, \rho}(M) = \Sx(M) / \fm_\rho \Sx(M).
 \ee

We prove that given an irreducible representation of the fundamental group of the handlebody, the reduced  skein module of a handlebody is an irredicuble representation of the skein algebra of its boundary surface. More precisely, 
\begin{THM}[See Theorem \ref{thmAzuRep}  ]\label{thmAzuRep0} 
Let $\xi\in \BC$ be a root of unity  with  $\ord(\xi^2)$ odd, and $H$ a handlebody of genus $g>1$. Let $\fS= \partial M= \Sigma_g$.
Assume $\rho:\pi_1(H) \rightarrow \G$ is an irreducible representation. 

The action of $\SxS$ on $\Sx(H)$ descends to an irreducible action of $\SxS$ on $\cS_{\xi, \rho}(H)$ which has dimension equal to $N^{3g-3}$, the PI-degree of $\SxS$. The classical shadow of this representation is the character of $\bar \rho$, which  is the restriction of $\rho$ onto $\pi_1(\fS)$. 
\end{THM}
Based on this theorem, we prove

\begin{THM}[See Theorem \ref{thmGeoKau}]  \label{thmGeoKau0} Let $\xi\in \BC$ be a root of unity  with  $\ord(\xi^2)$ odd, and $M$ a connected oriented closed 3-manifold. Let $\rho: \pi_1(M) \to SL_2(\BC)$ be an irreducible representation. Then the $\rho$-reduced skein module $\cS_{\xi, \rho}(M)$ is isomorphic to $\BC$.
\end{THM}

\begin{remark}There are full analogs of Theorems \ref{thmAzuRep0} and \ref{thmGeoKau0} for even $\ord(\xi^2)$. Their proofs require a strengthened version of \cite{FKL3},  and we defer them to a future work. \end{remark}

\subsection{On proofs} A large part of the paper is devoted to the proof of Theorem \ref{thm12} about a degeneration of the sliced skein algebra into a quantum torus. We  cut the surface into pairs of pants. For this we need to extend the definition of skein algebra to a version involving the boundary, and define  an analog of  well-known quantum trace homomorphisms \cite{BW,LY2} that helps us to glue back the pairs of pants and  capture the top degree term of the product of skeins. We also need to modify the definition of the Dehn-Thurston coordinates that help us to pick up the top degree term, as the ordinary Dehn-Thurston coordinates do not behave well under skein product.

\subsection{Organization of the paper}
Section \ref{sec.alg} contains notations and algebraic facts.  We recall known facts about Azumaya loci and character varieties in Sections~\ref{secAzu} and ~\ref{secChar}. In  Section \ref{skein} we discuss skein modules, and in Section \ref{secMRY} we introduce the sliced skein algebra. The modified Dehn-Thurston coordinates are defined in Section \ref{secDT}. In Section \ref{secQT} we introduce a quantum trace map for the basic pairs of pants. In section \ref{secDegen} we prove a finer version of Theorem \ref{thm12}. We prove finer versions of Theorems \ref{thmCenter00} and \ref{thmAzu0} in Sections \ref{secCenter} and \ref{secthmAzu}. Theorems \ref{thmAzuRep0} and \ref{thmGeoKau0} are proved in Section \ref{secGeoKau}.

\subsection{Acknowledgments} The authors  thank  F. Bleher, F. Bonahon,  F. Costantino, R. Kinser, 
 E. Letellier and   A. Sikora for helpful discussions.  {We also thank T. Yu for suggesting Proposition \ref{rSliceSmooth}.  This material is based upon work supported by the NSF grant DMS-2203255 (T.L.) and done while serving (J.K.B.) at the National Science Foundation. Any opinion, findings, and conclusions or recommendations expressed in this material are those of the authors and do not necessarily reflect the views of the National Science Foundation.

The authors presented the results of the current paper at various conferences, in particular at Joint Math Meetings  in Denver in January 2020, at virtual JMM in January 2021, and at the "Quantum Topology and
Geometry" conference in Paris in June 2022, and would like to thank the organizers for allowing us to present our work.

\section{Notation and algebraic preliminaries}\label{sec.alg}
\def\BCx{{\BC^\times}}
\def\ZQ{{\BZ[q^{\pm1}]}}
In this section we review the theory of quantum tori and filtered algebras.
We start the section by fixing notation. 

\subsection{Notation and  conventions} Let $\BN, \BZ, \BC$ be respectively the set of non-negative integers, the set of integers, and the set of complex numbers. Note that $\BN$ contains~0. Let $\BCx$ be the set of non-zero complex numbers.  For a root or unity $\xi\in \BC$ let $\ord(\xi)$ be the smallest positive integer $n$ such that $\xi^n=1$. 

All  rings are associative and  unitary, and ring  morphisms preserve the unit.
A (potentially  non-commutative) ring is a {\bf domain} if $ab=0$ implies $a=0$ or $b=0$. For a commutative algebra $\cR$, an {\bf $\cR$-domain} is a potentially non-commutative   $\cR$-algebra which is a domain.

In the paper $q$ or $q^{1/2}$ is always a formal variable.
Denote by $\Zq$  the ring of Laurent polynomials in  $q^{\pm 1/2}$ with integer coefficients.  If $\cR$ is a commutative $\Zq$-domain and $A$ is a $\Zq$ module or $\Zq$-algebra denote $A_\cR= A \ot _\Zq \cR$. 
Many properties of $A$ transfer to $A_\cR$, for example freeness. But, notably, ``being a domain" does {\em not} transfer.

In this section we fix a commutative $\Zq$-domain $\cR$.

\subsection{Weyl ordering, algebra with reflection}  \label{ss.reflection}

Two elements $x,y$ in a $\Zq$-algebra $A$ are {\bf{$q$-proportional}}, denoted by $x\qeq y$, if there is $k\in \BZ$ such that $x = q^{k} y$. Two elements $x,y\in A$ are {\bf $q$-commuting} if $xy$ and $yx$ are $q$-proportional.

 Suppose $x_1,x_2,\dots,x_n$ are pairwise $q$-commuting elements, $x_i x_j = q^{c_{ij}} x_j x_i$.
 The well-known {\bf Weyl normalization} of the product $x_1x_2\dots x_n$ is
\[[x_1x_2\dots x_n]_\Weyl=q^{-\frac{1}{2}\sum_{i<j}c(x_i,x_j)}x_1x_2\dots x_n.\]
It is easy to check that if $\sigma$ is a permutation of $\{1,2,\dots,n\}$, then 
$$[x_1x_2\dots x_n]_\Weyl=[x_{\sigma(1)}x_{\sigma(2)}\dots x_{\sigma(n)}]_\Weyl.$$

 A \term{$\Zq$-algebra with reflection} is a $\Zq$-algebra $A$ equipped with a $\BZ$-linear anti-involution $\omega$, called the \term{reflection}, such that $\omega(\qq)=q^{-1/2}$. In other words, $\omega : A \to A$ is a $\BZ$-linear map such that for all $x,y \in A$,
\[\omega(xy)=\omega(y)\omega(x),\qquad \omega(\qq x)=q^{-1/2} \omega(x).\]

An element $z\in A$ is called \term{reflection invariant} if $\omega(z)=z$. If $A'$ is another $\Zq$-algebra with reflection $\omega'$, then a map $f:A\to A'$ is \term{reflection invariant} if $f\circ \omega=\omega'\circ f$.

The following trivial statement is very helpful  as it allows to do calculations up to powers of $q$ and fix the exact powers at the end.
\blem\label{rReflection} Let $a, b$ be reflection invariant elements of a $\Zq$-algebra with reflection.

(a) If $a \qeq b$  then $a=a'$.

(b) If $a$ and $b$ are $q$-commuting, then $[ab]_\Weyl$ is reflection invariant.

\elem

\subsection{Quantum tori} \label{ssQtorus}
Let $Q$ be an antisymmetric $r\times r$ integral matrix. 
The {\bf quantum torus} associated to $Q$ is the algebra
\be 
\mathbb{T}(Q)\eqdef \Zq\langle x_1^{\pm1},\dots,x_r^{\pm1}\rangle/\langle x_ix_j=q^{Q_{ij}}x_jx_i\rangle.
\label{eqT}
\ee
A quantum torus $\bT(Q)_\cR$ is an Ore domain \cite{GW}. 

For $\bk=(k_1,\dots, k_r)\in \BZ^r$, define the Weyl normalized monomial
\[ x^\bk \eqdef [ x_1 ^{k_1} x_2^{k_2} \dots x_r ^{k_r} ]_\Weyl = q^{-\frac{1}{2} \sum_{i<j} Q_{ij} k_i k_j} x_1 ^{k_1} x_2^{k_2} \dots x_r ^{k_r}\]
 Note that $\{ x^\bk \mid \bk \in \BZ^r\}$ is a free $\Zq$-basis of $\bT(Q)$. For $\bk, \bk'\in \BZ^r$ let
 \be \la \bk, \bk'\ra_Q := \sum_{1\le i, j \le r} Q_{ij} k_i k'_j. 
 \label{eqQform}
 \ee It is easy to check that
\begin{align}
\label{eq.prod}
x^\bk x ^{\bk'}& = q ^{\frac 12 \la \bk, \bk'\ra_Q} x^{\bk + \bk'}.
\end{align}
\no{Thus the following decomposition gives the algebra $\bT(Q)$ a $\BZ^r$-grading:
\begin{equation}\label{eq.grad}
\bT(Q) = \bigoplus_{\bk\in \BZ^r} \cR x^\bk
\end{equation}

Suppose $Q'$ is another antisymmetric $r'\times r'$ integral matrix such that $HQ' H^T= Q$, where $H$ is an $r\times r'$ integral matrix and $H^T$ is its transpose. Then the $\Zq$-linear map $\bT(Q)\to \bT(Q')$ given on the basis by $x^\bk \mapsto x^{\bk H}$ is an algebra homomorphism, called a \emph{multiplicatively linear homomorphism}. Here $\bk H$ is the product of the row vector $\bk$ and the matrix $H$.
}

The quantum torus $\bT(Q)$ has a reflection anti-involution
\[\omega:\mathbb{T}(Q)\to\mathbb{T}(Q), \quad \text{given by} \quad
\omega(q^{1/2})=q^{-1/2},\quad \omega(x_i)=x_i.\]
All normalized monomials $x^\bk$ are reflection invariant.

If $\Lambda\subset \BZ^r$ is a submonoid, then the $\Zq$-submodule
$\bT(Q,\Lambda)\subset \bT(Q)$ spanned by $\{ x^\bk\mid \bk \in \Lambda \}$ is a $\Zq$-subalgebra of $\bT(Q)$, called a {\bf monomial subalgebra}.

\subsection{Filtrations and associated graded algebras} 
\label{ssFiltr} In this paper, an {\bf ordered monoid} $\Gamma$ is a submonoid of $\BZ^r$ for some $r\in \BN$, equipped with a linear order $\le$ such that if $a \le a'$ then $a+b \le a'+b$ for all $b\in \Gamma$. Fix an ordered monoid $\Gamma$.

A {\bf $\Gamma$-filtered $\cR$-module} $A$
 is an $\cR$-module equipped with a  {\bf $\Gamma$-filtration}, which is a family $F= (F_k(A))_{k\in \Gamma}$ of $R$-submodules of $A$ such that $F_k(A)\subset F_l(A)$ if $k \le l$ and $\cup_{k\in \Gamma} F_k(A) = A$. The {\bf associated graded module} of $F$ is 
$$
  \Gr^F(A):=\bigoplus_{k\in\Gamma}\Gr^F_k, \ \text{where $\Gr_k^F:=F_k/F_{<k}$ and $F_{<k}:=\displaystyle{\sum_{k'<k}F_{k'}}$. }
$$
An $\cR$-linear map $f: A \to A'$ between two $\Gamma$-filtered module {\bf respects the $\Gamma$-filtrations} if $f( F_k(A)) \subset F_k(A')$ for all $k \in \Gamma$. Such a map induces the {\bf associated graded map}
$$ \Gr^F(f): \Gr^F(A) \to \Gr^F(A'),\   a  + F_{<k}(A) \to f(a) + F_{<k} (A')\ \text{for } \ a\in F_k(A).$$

 A $\Gamma$-filtration $F$ is {\bf good} if for every non-zero $a\in A$, there is $k\in \Gamma$, denoted by $\deg^F(a)$, such that $a \in F_k(A) \setminus F_{<k}(A)$. Let the {\bf lead term} of $a$ be defined by $$\lt^F(a) := p_k(a)\in \Gr^F(A),$$ where $p_k: F_k \onto \Gr^F_k$ is the natural projection. 
 By convention $\lt^F(0)=0$. Note that $\lt^F(a) \neq 0$ unless $a=0$.
 
 \blem \label{rLift1}
  Let $A, A'$ be $\Gamma$-filtered $\cR$-modules and $f: A \to A'$ be an $\cR$-linear map respecting the $\Gamma$-filtrations. Assume the $\Gamma$-filtration of $A$ is good.
If $\Gr^F(f): \Gr^F(A) \to \Gr^F(A')$ is an isomorphism then $f^{-1}(F_k(A')) = F_k(A)$ and $f$ is injective.
\elem
\bpr
Since $f(F_k(A)) \subset F_k(A')$  it is clear that $f^{-1}(F_k(A')) \supset F_k(A)$.
Assume the contrary that there exists $a \in f^{-1}(F_k(A')) \setminus F_k(A)$.
Then $l:=\deg^F(a) >k$. Because $f(a) \in F_k(A') \subset F_{< l}(A')$ we have $\Gr^F(f)(a)=0$, contradicting the fact that $\Gr^F(f)$ is an isomorphism.
\epr

\no{ Assume $A$ is a free $\cR$-module. An $\cR$-basis $S$ of $A$ respects to a $\Gamma$-filtration $F$ if $S \cap F_k(A)$ is a free $\cR$-basis of $F_k(A)$ for all $k\in \Gamma$.}

When $A$ is an $\cR$-algebra, we say a  $\Gamma$-filtration $F$ {\bf respects the product}, or we call $F$ an {\bf algebra $\Gamma$-filtration}, if $1\in F_0(A)$ and $F_k(A) F_l(A) \subset F_{k+l}(A)$. In this case $\Gr^F(A)$ has an $\cR$-algebra structure defined by
$
p_k(x)p_{k'}(y)=p_{k+k'}(xy).
$

Using the lead term, one can easily prove
\begin{proposition}
\label{liftfacts} 
  If the filtration is good and 
 $\mathrm{Gr}^F(A)$ is a domain, then $A$ is a domain.\qed

\end{proposition}

\def\reordno{  \raisebox{-8pt}{\incl{.8 cm}{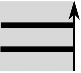}} }

\def\reordnod{  \raisebox{-8pt}{\incl{.8 cm}{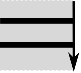}} }

\def\reordoneall{  \raisebox{-8pt}{\incl{.8 cm}{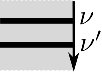}} }

\def\reordoneallp{  \raisebox{-8pt}{\incl{.8 cm}{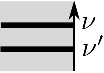}} }

\def\reordthree{  \raisebox{-8pt}{\incl{.8 cm}{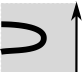}} }
\def\Coeff{{\mathrm{Coef}}}

\def\bqq{{\, \overset{(q)}= \, }}

 \def\sS{{\mathscr S}}
\def\bSS{{\overline \sS(\fS)}}

\def\Do{{\mathring \Delta}}
\def\bQ{{\overline{ \mathsf Q }}}
\def\hDd{{\hat \Delta_\partial}}
\def\DD{{ \tilde \Delta}  }
\def\mQ{{\mathsf Q}}
\def\bsX{{\overline{\sX}}}

\newcommand{\facedef}[2]{
\begin{tikzpicture}[scale=0.8,baseline=0.28cm]
\draw[fill=gray!20!white] (0,1)--(0.6,0)--(1,1);
\draw[inner sep=0pt] (0.1,0.5)node{\vphantom{$b$}#1} (1,0.5)node{\vphantom{$b$}#2};
\draw[fill=white] (0.6,0)circle(2pt);
\end{tikzpicture}
}

\def\bD{\DD}
\def\bx{{\bar x}}
\def\sXo{{\overset{\diamond}{\sX}}}

\def\fT{{\mathfrak T}}
\def\bST{{\overline{\mathscr{S}}(\ft)}}
\def\BAT{{\bT^A(\ft)}}
\def\BXT{{\overline {\sX}(\ft)}}
\def\btr{{ \overline{\tr}  }}

\def\uB{{\underline{\mathrm{B}}}}
\def\PP{{\mathbb P}}
\def\PPb{{\check{\mathbb P}}}

\newcommand{\Move}[1]{\orange{Move: #1}}
\def\Irrep{{\mathrm{Irrep}}}
\def\End{{\mathrm{End}}}
\def\csh{{\mathsf{csh}}}
\def\A{{\mathscr A}}

\section{Azumaya loci} \label{secAzu}
We collect necessary facts about the Azumaya loci and Poisson orders in this section. As we will deal only with $\BC$-domains, we restrict the treatment to this class of algebras.

\subsection{Almost Azumaya $\BC$-domain} We recall  the definition of Azumaya loci.

\bdf A $\BC$-domain $A$ is {\bf almost Azumaya} if \begin{itemize}
\item $A$ is an affine $\BC$-algebra, meaning it is a finitely generated $\BC$-algebra, and
\item $A$ is finitely generated as a  $Z$-module, where $Z$ is its center.
\end{itemize}
\edf

Fix an almost Azumaya $\BC$-domain $A$ with center $Z$. Let $\Fr(Z)$ be the field of fractions of $Z$. Define $$\dim_ZA:=\dim_{\Fr(Z)} \Fr(Z) \ot_Z A.$$
By Posner's Theorem \cite[Theorem 13.6.5]{MR}, $\Fr(Z) \ot_Z A$ is a division algebra with center $\Fr(Z)$.
Hence there is  $D\in \BN$, called the {\bf PI-degree of $A$}, such that  $\dim_ZA= D^2$.

By the Artin-Tate lemma \cite[Theorem 13.9.10]{MR}, the center $Z$ is an affine $\BC$-algebra. The maximal spectrum $\MaxSpec(Z)$ is a $\BC$-affine irreducible variety.

A non-zero irreducible representation $\rho: A \to \End(V)$ defines a point $\csh(\rho)\in \MaxSpec(Z)$, called the the {\bf classical shadow} of $\rho$, as follows. By Schur's lemma, for $z\in Z$ there is a scalar $\chi_\rho(z)\in \BC$ such that $\rho(z) = \chi_\rho(z)\, \Id$. Since $\chi_\rho: Z \to \BC$ is a surjective $\BC$-algebra map, its kernel $\ker(\chi_\rho)$ is  a maximal ideal, or an element of $\MaxSpec(Z)$. Define $\csh(\rho)= \ker(\chi_\rho)$.
 We  call $\MaxSpec(Z)$ the {\bf classical shadow variety of $A$}.

\bdf A point $\fm\in \MaxSpec(Z)$ is {\bf Azumaya} if $A/\fm A \cong \Mat_D(\BC)$, the algebra of $D \times D$ matrices with entries in $\BC$. The set $\A(A)\subset \MaxSpec(Z)$ of all Azumaya points is called the {\bf Azumaya locus}.
\edf
Since $\Mat_D(\BC)$ is simple and its only irreducible representation is $\BC^D$ having dimension $D$, if $\csh(\rho)$ is Azumaya, then $\dim \rho= D$.
The following facts are well-known.
\bthm[ See Section III.1 of \cite{BG}] \label{thmAzu} Let $A$ be an almost Azumaya $\BC$-domain with PI-degree $D$.

(a) 
The Azumaya locus $\A(A)$  is  Zariski open dense in $\MaxSpec(Z)$. 

(b) If  $\fm \in \cA(A)$, then there is a unique, up to conjugations, irreducible representation with classical shadow $\fm$, called the Azumaya representation of $\fm$. The Azumaya representation has dimension $D$.

 (c) If $\csh(\rho)\not \in \A(A)$ then $\dim \rho < D$.
 
 (d) If $\fm \in \cA(A)$ and $V$ is a non-zero $A/\fm A$-module with $\dim_\BC(V) \le D$,  then $V$ is the Azumaya representation of $\fm$.

\ethm

\def\hZ{{\hat Z}}
\def\BCx{{\BC^\times}}
\def\Zz{ Z_\zeta}
\def\Cq{{\BC[q^{\pm 1}]}}
\def\cL{{\mathcal L}}
\def\smooth{{\mathrm{sm}}}

\bcor \label{rQuotient}
Let $f: A \to B$ be a surjective $\BC$-algebra homomorphism between two almost
Azumaya $\BC$-domains having the same $PI$-degree $D$. The natural map $\MaxSpec(Z(B)) \to \MaxSpec(Z(A))$ sends $\cA(B)$ into $\cA(A)$.
\ecor
\bpr  From surjectivity we have $f(Z(A)) \subset Z(B)$. Denote the restriction by
  $\bar f: Z(A) \to Z(B)$. Assume $\fm \in \cA(B)$, with $\bar \fm = \bar f^{-1}(\fm)\in \MaxSpec(Z(A))$. The composition 
$$ A / \bar \fm A \onto B/\fm B \cong M_D(\BC)$$
shows that $A/\bar \fm A$ has an irreducible representation of dimension $D$, which is the PI-degree of $A$. By Theorem \ref{thmAzu}(c), we have $\bar \fm \in \cA(A)$. This proves the Corollary.
\epr

\subsection{Poisson algebras} A {\bf Poisson algebra} is a   $\BC$-algebra $A$ equipped with a Lie  bracket $\{\cdot, \cdot\}: A \ot_\BC A \to A$ such that $\{xy, z\}=x  \{y, z\} + y\{x, z\} $. An element $x$  of a Poisson algebra $A$ is {\bf Casimir} if $\{x, y\}=0$ for all $y\in A$.

A $\BC$-algebra map $f: A \to A'$ from one Poisson algebra to another is a {\bf Poisson morphism} if there is a non-zero $k\in \BC$  such that $f(\{a, b\}) = k \{ f(a), f(b)\}$ for all $a,b,\in A$.

\bdf {\bf A  symplectic variety with singularities} is 
an irreducible $\BC$-affine  variety $X$ equipped with a Poisson structure on  $\BC[X]$ such that the Poisson bracket is non-degenerate on the smooth locus $X^\smooth$ of $X$. In other words, there is a symplectic $2$-form $\omega$ on the smooth locus $X ^\smooth$ such that for $f, g\in \BC[X]$ we have $\{ f, g\} = \omega( df, dg)$ on  $X ^\smooth$.
\edf
In this case $X^\smooth$ is the symplectic leaf of $X$ of largest dimension in the sense of \cite{Weinstein}.

\subsection{Poisson structure, Azumaya locus and Poisson leaves}

We present  a connection between symplectic leaves in $\MaxSpec(Z)$ and the Azumaya locus, following \cite{BG2}.

Let $A$ be a torsion free $\Cq$-algebra. For $\xi \in \BCx$ define the $\BC$-algebra
$$   A_\xi= A/(q-\xi)Z \quad \text{with projection $\pi: A \to A_\xi$.}$$
Since $A$ has  no $\Cq$-torsion, if $\pi(x)=0$, then $x = (q-\zeta)y$ for a {\em unique} $y \in A$; or  $x /(q-\zeta) \in A$ is well-defined.

Let $Z(A_\xi)$ be the center of $A_\xi$. The following is well-known.
\blem \label{rPoisson}
Let $x,y\in Z(A_\xi)$. Choose lifts $\hat x, \hat y\in A$, that is, $\pi(\hat x)=x, \pi(\hat y) =y$. 

(a)  The following element does not depend on the choices of $\hat x, \hat y$:
\be  \{ x, y\} := \pi \left (\frac{[\hat x, \hat y ]}{q-\xi}  \right)
\label{eqBra}
\ee

(b) One has $\{x, y\} \in Z(A_\xi)$.

(c) The bracket $\{\cdot, \cdot\}$ is a Poisson bracket on $Z(A_\xi)$.

\elem
\bpr (a) Another lift of $x$ must have the form $\hat x + (q-\xi) \hat u$, where  $\hat u \in A$. Then
$$ \pi \left (\frac{[\hat x+ (q-\xi)\hat u, \hat y ]}{q-\xi}  \right) =  \pi \left (\frac{[\hat x, \hat y ]}{q-\xi}  \right) + \pi ([\hat u, \hat y])= \pi \left (\frac{[\hat x, \hat y ]}{q-\xi}  \right) +[\pi \hat u, y] = \pi \left (\frac{[\hat x, \hat y ]}{q-\xi}  \right). $$
Similarly one shows that the right hand side of \eqref{eqBra} does not depend on a choice of $\hat y$.

(b) 
Let $a\in A_\xi$ with lift $\hat a\in A$. The Jacobi identity for commutators gives
\be\left  [ \frac{[\hat x, \hat y]}{q-\xi} , \hat a \right] = \left [\hat x, \frac{[\hat y, \hat a]}{q-\xi} \right] + \left  [ \frac{[\hat x, \hat a]}{q-\xi} , \hat y \right].
\label{eqJac}
\ee
Applying $\pi$ to both sides, we get $[\{x,y\}, a]=0$ for all $a\in A_\xi$. Hence $\{x, y \} \in Z(A_\xi)$.

(c) From \eqref{eqBra} we have $\{x, y \}= - \{ y, x \}$. In  \eqref{eqJac} let $a\in Z(A_\xi)$. Divide  \eqref{eqJac} by $(q-\xi)$ and apply $\pi$, we get the Jacobi identity for $\{ \cdot, \cdot\}$. From identity
$$ [ \hat x, \hat y \hat z] =  \hat y [ \hat x,  \hat z] +  [ \hat x, \hat y] \hat z$$
we get that $\{x, \cdot\}$ is a derivation. Hence $\{\cdot, \cdot \}$ is a Poisson bracket.
\epr
We call $\{\cdot, \cdot\}$ the {\bf quantization Poisson bracket} of the center of $A_\xi$  at $q=\xi$.

The following is clear from the definition, and will be very handy later.
\bpro \label{rPoiMor}
Let $A$ and $B$ be  torsion free $\Cq$-algebras, and $f: A\to B$  a $\Cq$-algebra  homomorphism. Assume $a,a'\in Z(A_\xi)$ such that $f_\xi(a), f_\xi(a')\in Z(B_\xi)$. Then 
\be  f_\xi(\{ a, a'\} ) = \{f_\xi( a), f_\xi(a')\}.
\ee
In particular, if 
$f_\xi(Z(A_\xi)) \subset Z(B_\xi)$ (for example when $f: A \to B$ is surjective) 
then the restriction $f_\xi: Z(A_\xi) \to Z(B_\xi)$ is Poisson. 
\epro

The following is a main fact about Azumaya locus that we will use.  

\bpro[ See \cite{BG2}]
 \label{rPoiAzu}  Let $A$ be a torsion free $\Cq$-algebra and $\xi \in \BCx$. 
 Assume
\begin{enumerate}
\item[(a)] $A_\xi$ is an almost Azumaya $\BC$-domain, and
\item [(b)]   $\MaxSpec(Z(A_\xi))$ is a symplectic variety with singularities.

\end{enumerate}
 Then the smooth locus $\MaxSpec(Z(A_\xi)) ^\smooth$ is in the Azumaya locus $\A(A_\xi)$.
\epro
\bpr Let us prove that $A_\xi$ is  a Poisson $Z(A_\xi)$-order in the sense of \cite{BG2}.

Let $\{x_i, i \in I\}$ be a $\BC$-basis of $Z(A_\xi)$. Fix once and for all lifts $\hat x_i$ of $x_i$. For $a\in A_\xi$ choose a lift $\hat a\in A$ and define
$$ D_{x_i}(a)= \pi \left (\frac{[\hat x_i, \hat a] }{q-\xi}  \right) .$$
It is easy to check that $D_{x_i}(a)$ does not depend on the choice of $\hat a$, and that $D_{x_i}$ is a $\BC$-derivation of $A_\xi$. By linearity define $D_x$ for $x\in Z(A_\xi)$. From the definition it follows that for $x,y\in Z(A_\xi)$ we have $D_x(y) = \{ x,y\}$.
Thus we have the following 
\begin{itemize}
\item $Z(A_\xi)$ is stable under the action of $D_x, x\in Z(A_\xi)$,
\item The map $(x,y) \to D_x(y)$ gives a Poisson bracket on $Z(A_\xi)$, and
\item $A_\xi$ is an affine $\BC$-algebra finitely generated over its center.
\end{itemize}
This is the definition of a
 Poisson $Z(A_\xi)$-order \cite{BG2}.

By Condition (b) we have that  $(\MaxSpec(Z(A_\xi)) ^\smooth$ is a symplectic leaf.  By \cite[Theorem 4.2]{BG2}
 if $\fm, \fm'\in \MaxSpec(Z(A_\xi))) ^\smooth$ then $A/\fm A \cong A / \fm' A$. Hence if one point of $\MaxSpec(Z(A_\xi))^\smooth$  is in $\A(A_\xi)$, then the whole $\MaxSpec(Z(A_\xi)) ^\smooth$  is in $\A(A_\xi)$. Irreducibility  implies that $\MaxSpec(Z(A_\xi)) ^\smooth$ is open and dense in $\MaxSpec(Z(A_\xi))$. Since 
  $\A(A_\xi)$ is also open and dense in $\MaxSpec(Z(A_\xi))$, it has a common point with $\MaxSpec(Z(A_\xi)) ^\smooth$. The statement follows.
\epr

\section{Character varieties}

\label{secChar}

We recall the necessary facts about the $SL_2(\BC)$ character variety, its sliced subvariety, and the Poisson structure. We also discuss the $PSL_2(\BC)$ character variety, and give a criterion when a point is smooth on a sliced character variety.

\subsection{Character variety} We recall  the $SL_2(\BC)$ character variety, its cousin $PSL_2(\BC)$ character variety, and define the even part of the character variety of a surface.

Let $M$ be a connected manifold. The set $\Hom(\pi_1(M), \G)$ of all group homomorphisms from the fundamental group $\pi_1(M)$ to $SL_2(\BC)$ is a $\BC$-algebraic set. The group $\G$ acts on $\Hom(\pi_1(M), \G)$ by conjugation. The {\bf $SL_2(\BC)$ character variety $\X(M)$} is the set of closed points of the GIT quotient  $\Hom(\pi_1(M), \G)//\G$. If we replace $SL_2(\BC)$ by $PSL_2(\BC)$ we get the $PSL_2(\BC)$ character variety.

Alternatively, one can define $\X(M)$ by
$$ \X(M) = \MaxSpec(\BC[ \Hom(\pi_1(M), \G)]^{\G}).$$
Two representations $\rho, \rho'\in \Hom(\pi_1(M), \G)$ descend to the same element of $\X(M)$ if and only if they have the same character. Thus, we consider $\X(M)$ as the set of all $SL_2(\BC)$-characters of $\pi_1(M)$. For a representation $\rho:\pi_1(M)\to SL_2(\BC)$ denote by $[\rho]$ its character, considered as an element of $\X(M)$.

A closed curve $\al\subset M$ defines a regular function $ T_\al:\X(M)\to \BC$, depending only on the homotopy class of $\al$, as follows. Let $\tilde \al \in \pi_1(M)$ be any element  represented by $\al$. Define
\be  T_\al([\rho]):= \tr(\rho(\tilde \al)).
\label{eqCXM}
\ee
\def\BfS{\bar \fS}

\def\XSw{{\X^\vw(\fS)}}
\def\XeSw{{\X^{\vw}(\fS)^\ev}}
\def\XSev{\X(\fS)^\ev}
\def\li{{\mathsf {li}}}

Consider $\fS= \Sigma_{g,m}\subset \Sigma_g$. Let $\BfS= \Sigma_g$. It is known (see for example \cite{Letellier}) that $\XS$ is irreducible as an affine algebraic set.

 The group $G:= H^1(\BfS;\BZ/2)\equiv \Hom(\pi_1(\BfS), \{\pm 1\})$ acts on
  $\Hom(\pi_1(\fS); SL_2(\BC))$ by
\be  u * \rho(\al) = u(\al) \rho(\al), u \in G.
\label{eqAction}
\ee

The action descends to an action of $G$ on  $\X(\fS)$. The  GIT quotient $\XSev:=\X(\fS)//G$ will be called the {\bf even} character variety. In other words, 
 $$ \XSev =\MaxSpec( \BC[\X(\fS)]^G).$$
 As a quotient of an irreducible variety, $\XSev $ is also irreducible.
 
The following explains the relation between $\XSev$ and the $PSL_2(\BC)$-character variety, though we don't need it in the current paper.
\bpro (a) If $m=0$ then $\XSev$ is the connected component of the $PSL_2(\BC)$ character variety $\X^{PSL_2(\BC)}(\fS)$ containing the trivial character.

(b) If $m>0$ then $\XSev/(\BZ/2)^m$ is the $PSL_2(\BC)$-character variety $\X^{PSL_2(\BC)}(\fS)$. Here $(\BZ/2)^m $ is identified with $H^1(\BfS, \fS;\BZ/2))$ and the  action is given by the same formula~\eqref{eqAction}.
\epro
\bpr Let $H^1(\fS;\BZ/2 ) $ act on $\XS$ by the  \eqref{eqAction}.
Then $\XSev/(\BZ/2)^m= \X(\fS)//H^1(\fS;\BZ/2 )$, which by \cite[Proposition 4.2]{HP} is equal to
 $\X^{PSL_2(\BC), \li }(\fS)$, 
 the subvariety of $\X^{PSL_2(\BC)}(\fS)$
  consisting of characters of 
  $PLS_2(\BC)$-representations which lift to $SL_2(\BC)$-representations. 
  When $m >0$, one has  $\X^{PSL_2(\BC), \li }(\fS)= \X^{PSL_2(\BC)}(\fS)$ by  \cite{HP}, which proves (b).
 
Assume $m=0$. Then $(\BZ/2)^m$. Hence $\XSev= \X^{PSL_2(\BC), \li }(\fS)$, which, by \cite{Culler}, consists of several connected components of $\X^{PSL_2(\BC), \li }(\fS)$. As $\XSev$ is connected and contains the trivial character, we get (a).
\epr

\subsection{Sliced character variety} We now define the sliced character varieties.

Let $\fS=  \Sigma_{g,m}$ with  peripheral   loops $\fd_1, \dots , \fd_m$ and let $\BfS= \Sigma_g$. Consider the regular function
\be 
\fd: \X(\fS)\to \BC^m, \quad \fd([\rho])= (\tr( \rho(\fd_1)), \dots, \tr( \rho(\fd_m)).
\ee

For $\vw \in \BC^m$ the affine set
$ \XSw:= \fd^{-1}(\vw),$ is called 
 a {\bf sliced character variety}. 

\def\bG{{\bar G }}
\def\Gp{G_\partial}

The action of the group  $G=H^1(\bfS;\BZ/2)$ on $\XS$ given by \eqref{eqAction} 
 preserves each slice $\X^\vw(\fS)$.
 Let  $\XeSw=\XSw//G$, which is the image of $\XSw$ under the projection $\X(\fS) \onto \XSev$.
 We have the stratifications
\be 
\X(\fS) = \bigsqcup_{\vw \in \BC^m  }\XSw, \quad \XSev = \bigsqcup_{\vw \in \BC^m  }\XeSw.
\ee
By \cite{Whang}, see also Corollary \ref{rInt}, 
each sliced character variety $\XSw$ is irreducible.
This was known \cite{HLR2} before for
 generic $\vw$.  Consequently $\XeSw$ is also irreducible.

\subsection{Atiyah-Bott-Goldman Poisson structure} 

There exists a Poisson structure on $\BC[\X(\fS)]$, following work of Atiyah-Bott and Goldman, see \cite{Goldman,Lawton,GHJW}, which will be called the ABG-Poisson structure.  It is defined up to a non-zero scalar constant. 
In \cite{Lawton,GHJW} it is shown that the ABG Poisson bracket can be restricted to each sliced character variety $\XSw$ and  that the smooth part $\XSw^\smooth$ is a symplectic complex manifold. As $\XSw$ and $\XeSw$ are irreducible,  we have the following.
\bpro \label{rSliceSym}
For each $\vw\in \BC^m$ the sliced character variety $\XSw$ and its quotient $\X^\vw(\fS)^\ev$, with the ABG-Poisson structure, are symplectic varieties with singularities.  
\epro

{
\subsection{Smooth points of the sliced character variety}  A character $x\in \XS$ is {\bf sliced-smooth} if $x$  is smooth in the sliced character variety containing $x$. 
\bpro \label{rSliceSmooth}
(a) A smooth point of $\XSw$ is smooth as a point in $\XS$.

(b) Similarly, a smooth point of $\X^\vw(\fS)^\ev$  is smooth as a point in $\X(\fS)^\ev$. 
\epro

\blem \label{rSmoothX}
Let $X$ be an irreducible $\BC$-affine variety and $Y$ be a subvariety of $X$ cut out by $m$ regular functions, i.e. there are $m$ regular functions $g_1, \dots , g_m : X \to \BC$ such that $Y = \bigcap_{i=1}^m g^{-1}(0)$.  Assume $Y$ is  irreducible and $\dim(X) = \dim(Y) +m$, then any smooth point of $Y$ is smooth in $X$.
\elem

\bpr Assume
 $X$ is  defined by a regular function $F=(f_1, \dots, f_k): \BC^n \to \BC^k, \ X = F^{-1}(0)$.  
For $x\in X$ let $dF_x: T_x X \to T_{f(x)} \BC^k$ be the derivative map. 
Then
\be \dim \ker (dF_x) \ge \dim(X),
\ee
and  irreducibility of $X$ implies that equality happens exactly when $x$ is smooth in $X$.
Hence we have a criterion: $x$ is smooth  in $X$ if and only if 
\be \dim \{ d(f_1)_x, \dots, d(f_k)_x\} \ge n -\dim(X)
\label{eq1}
\ee

 By this criterion, at a point $x$ smooth in $Y$ we have
$$ \dim \{ d(f_1)_x, \dots, d(f_k)_x, d(g_1)_x, \dots, d(g_m)_x\} \ge n - \dim (Y)= n -\dim(X)+m.$$

By dropping the last $m$ vectors, we have \eqref{eq1}, proving that $x$ is smooth in $X$.
\epr
\bpr[Proof of Proposition \ref{rSliceSmooth}] (a) Both  $\XS$ and $\XSw$  are irreducible affine varieties of dimensions respectively  $d=6g-6+3m $ and $d-m=6g-6+2m $ (see \cite{Whang2}). Besides, $\XSw$ is cut out of $\XS$ by $m$ regular functions. Hence Lemma \ref{rSmoothX} implies Proposition \ref{rSliceSmooth}(a). The proof for (b) is identical, noting that $\dim(\XS)= \dim(\X(\fS)^\ev)$ and $\dim(\XSw)= \dim(\X^\vw(\fS)^\ev)$.
\epr

\subsection{Sufficient conditions for sliced-smoothness}
We present here some sufficient conditions for a character to be smooth in a sliced character variety.

 Note that  ``sliced-smooth" is different from ``smooth in $\XS$".


\def\irr{{\mathsf{irr}}}
\def\Xirr{\X^\irr}
\def\eee{\mathbf e}
\bpro \label{rSmooth}
Suppose  $\fS= \Sigma_{0,m}$ has genus 0. Assume $\vw \in (\BC \setminus \{2, -2\})^m$, and  $\tau\in \XSw$ is the character of an irreducible presentation.
 Then $\tau$ is  smooth in $\XSw$. 
\epro
\bpr 
We can assume that $\pi_1(\fS)$ is free with generators $x_1, \dots, x_{m-1}$ such that the free homotopy classes of $x_1, \dots, x_{m-1}$, and $x_m:= x_1 \dots x_{m-1}$ represent respectively $\fd_1, \dots, \fd_m$.
The subset $\Xirr\subset \XS$ of all irreducible characters is a Zariski open subset of the smooth locus $\XS^\smooth$. It is enough to show that the function $F: \Xirr \to \BC^m, F([\rho]) = (\tr(\rho(x_1), \dots, \tr(\rho(x_m))$ has $\vw$ as a regular value. Assume $[\rho] \in F^{-1}(\vw) \cap \Xirr$. We will show  $[\rho]$ is a regular for $F$ by showing that the derivative map $dF_{[\rho]}$ is surjective. Let 
$\eee_i$ be the $i$-th unit vector of $\BC^m$, which has entries 0 everywhere except for the $i$-th position which is 1. 

Let $X_i = \rho(x_i)$. For each $i\in  \{1,\dots, m-1\}$, by varying $X_i$ in $SL_2(\BC)$ and fixing $X_j$ for all other   $j\in \{1,\dots, m-1\}$, we see that there is a scalar $f_i$ such that $\eee_i+ f_i \eee_m$ is in the image of $dF_{[\rho]}$. It remains to show that $\eee_m$ is in the image of $dF_{[\rho]}$.

Since $\tr(X_m)\neq 2$ there is basis in which $X_m$ is diagonal with eigenvalues $\neq \pm 1$. The commutator group of $X_m$ is the group of diagonal matrices, which is an abelian group.
 Since $\rho$ is irreducible, there must be an $X_i$ not in this abelian group, or not commuting with $X_m$. Choose smallest such $i$. 
Then $A:= X_1 \dots X_i$ does not commute with $X_m$. 
Let $B= X_{i+1} \dots X_{m-1}$. Since  $AB= X_m$ and $A$ and $X_m$ do not commute, $A$ and $B$ do not commute. 
Consider two cases.

Case 1. One of $A, B$ has trace not equal to $\pm 2$. Say, $\tr(A)\neq \pm 2$; the case $\tr(B)\neq \pm 2$ is similar. Choose a basis where $A$ is diagonal
$$ A = \begin{pmatrix}
u & 0 \\ 0 & u^{-1}
\end{pmatrix},   u\neq \pm 1, \quad B = \begin{pmatrix}
a & b \\ c & d
\end{pmatrix}.
$$
Non-commutativity of $A$ and $B$ implies one of $b,c$ is non-zero. By exchanging the two basis vectors we can assume $b \neq 0$. Define the curve of representations $\rho_t$, with $t\in \BR$,  by
\be  \rho_t(x_j) =  \begin{cases} X_j & j \le i \\
 U_t \,  X_j \, U_t ^{-1} & i< j \le m-1 
\end{cases}, \quad \text{where} \ U_t = \begin{pmatrix} 
1 & t \\ 0 & 1
\end{pmatrix} .
\label{eqUt}
\ee
Note that $\rho_0=\rho$. A simple calculation shows that 
\be \frac{d F( [\rho_t] )  }{dt}\big|_{t=0} = b(u^{-1} -u) \eee_m  
\label{eqUt2}
\ee
which shows that $\eee_m$ is in the image of $dF_{[\rho]}$.

Case 2. Both $\tr(A), \tr(B) \in \{ \pm 2\}$. In the basis consisting of one eigenvector of $A$ and one eigenvector of $B$, we have
$$  A = \pm \begin{pmatrix}
1 & a \\ 0 & 1
\end{pmatrix},   \quad B = \pm \begin{pmatrix}
1 & 0 \\ b & 1
\end{pmatrix} .  $$
Non-commutativity of $A$ and $B$ means $ab\neq 0$. Define the curve $\rho_t$ by equation \eqref{eqUt}, with $U_t$ replaced by $\begin{pmatrix} 
1 & t \\ 1& 1+t
\end{pmatrix}.$
The right hand side of \eqref{eqUt2} now becomes $ \pm 2ab\, \eee_m  $, which shows that $\eee_m$ is in the image of $dF_{[\rho]}$ and completes the proof.
\epr 

The following is a sufficient condition for a slice to be smooth. For $\w\in \BC$ let $E(\w) = \{ z \in \BC \mid  z+ z^{-1} = \w\}$. For $S, S' \subset \BC$ let $SS'=\{ s s \mid s\in S, s'\in S'\}$. 
The $GL_2(\BC)$ version of the following was in \cite[ Theorem 3.2]{Letellier}.

\def\bode{{\boldsymbol \delta}}
\bpro  \label{rKostov} Let $\fS= \Sigma_{g,m}$. 
Assume  $\vw=(\w_1, \dots, \w_m)\in \BC^m$ is {\em Kostov-generic} in the sense that  $1 \not \in \prod_{i=1}^m E(\w_j)$, and besides, each $\w_i \neq \pm 2$. The sliced character variety $\XSw$ is a smooth algebraic irreducible variety.
\epro
\bpr 
Since $\w_i\neq \pm 2 $, the class of $GL_2(\BC)$-matrices with trace $\w_i$ is semi-simple. By \cite[ Theorem 3.2]{Letellier} the Kostov genericity implies the sliced  $GL_2(\BC)$-character variety $\X^{\vw,GL_2(\BC)}(\fS)$ is a smooth irreducible affine variety. We have a presentation
\be  \pi_1 (\fS) = \la a_1, b_1, a_2, b_2, \dots, a_g, b_g, \fd_1, \dots, \fd_m \mid \prod [a_i b_i] \fd_1 \dots \fd_m =1 \ra.
\label{eqPresent}
\ee

For $\rho: \pi_1(\fS)\to GL_2(\BC)$ let $[\rho]$ be its character.
Consider the function
$$ F: \X^{\vw,GL_2(\BC)}(\fS) \to \BC^{m-1},\  [\rho] \to (\det(\rho(\fd_1), \dots, \det \rho(\fd_{m-1}).$$
 If $\det(\rho(\fd_i))=1$ for $i=1,\dots,  m-1$ then $\det(\rho(\fd_m))=1$. Hence 
$ \XSw = F^{-1}(1,\dots, 1).$
It is enough to show that $(1,\dots,1 )\in \BC^{m-1}$ is a regular value of $F$.

Assume $[\rho]\in F^{-1}(1,\dots, 1)$. For $i=1, \dots, m-1$ let 
$$ \rho(\fd_i) = \begin{pmatrix}
a_i & b_i \\ c_i & d_i
\end{pmatrix} \in SL_2(\BC).
$$
Since $\tr(\rho(\fd_i)) \neq 2$, we have $c_i\neq 0$. Fix  $i\in \{1, \dots, m-1\}$.
For $t\in \BR$ 
let $\rho_t: \pi_1(\fS) \to GL_2(\BC)$ be defined so that $\rho_t(x)= \rho(x)$ for all generators of the presentation \eqref{eqPresent}, except for $\fd_i$ and $\fd_m$, where
$$ \rho_t( \fd_i) =  \begin{pmatrix}
a_i & b_i+t  \\c_i  & d_i
\end{pmatrix}, \ \rho_t(\fd_m) = \rho_t^{-1} \left( \prod [a_i b_i] \fd_1 \dots \fd_{m-1}  \right   ).
$$
Then $\frac{d}{dt}F(\rho_t)|_{t=0}= -c_i \eee_i$, where $\eee_i\in \BC^{m-1}$ is the $i$-th unit vector. Consequently $(1,\dots, 1)$ is a regular value.
\epr

Using the theory of Azumaya loci we can present examples of irreducible, but non-smooth points of the sliced character variety, see Corollary \ref{rNonSmooth}.

\section{Kauffman bracket Skein Modules}\label{skein}

In this section 
we define skein modules/algebras and  present basic facts: relation to character variety, Barrett's spin isomorphism and a similar map for  even parts, Chebyshev-Frobenius map and transparency, and center of skein algebras. We show that the quantization Poisson structure and the Atiyah-Bott-Goldman Poisson structure are equal.

In this section $M$ is an oriented 3-manifold and $\fS$ is an oriented  surface.


\def\KzM{K_\zeta(M)}
\def\KzuM{K_\zeta^u(M)}
\def\SM{{\cS(M)}}
\def\Sx{{\cS_\xi}}
\def\UM{{\mathsf{U}(M)}}
\def\KzuM{{\cS^u(M)}}
\def\UM{{\mathsf U}M}
\def\SzM{\cS_\zeta(M)}
\def\SzS{\cS_\zeta(M)}
\def\G{{SL_2(\BC)}}
\def\ChiM{{\X(M)}}
\def\rdd{{\mathrm{rd}}}
\def\Sx{{\cS_\xi}}
 \def\SxM{{\cS_\xi(M)}}
 \def\SMR{{\cS_\xi(M;\cR)}}
\subsection{Kauffman bracket skein module}
\label{ssKBSK}
We recall  the definitions and properties of skein modules/algebras.

Let $\cR$ be a commutative domain and let $\xi\in \cR$ be invertible.
  The Kauffman bracket skein module  $\SMR$
is the $\cR$-module freely spanned by isotopy classes of unoriented framed links in $M$, including the empty link,  subject to the Kauffman relations

\begin{align}
\label{eq.skein}
\begin{tikzpicture}[scale=0.8,baseline=0.3cm]
\fill[gray!20!white] (-0.1,0)rectangle(1.1,1);
\begin{knot}[clip width=8,background color=gray!20!white]
\strand[very thick] (1,1)--(0,0);
\strand[very thick] (0,1)--(1,0);
\end{knot}
\end{tikzpicture}
&=\xi 
\begin{tikzpicture}[scale=0.8,baseline=0.3cm]
\fill[gray!20!white] (-0.1,0)rectangle(1.1,1);
\draw[very thick] (0,0)..controls (0.5,0.5)..(0,1);
\draw[very thick] (1,0)..controls (0.5,0.5)..(1,1);
\end{tikzpicture}
+\xi^{-1}
\begin{tikzpicture}[scale=0.8,baseline=0.3cm]
\fill[gray!20!white] (-0.1,0)rectangle(1.1,1);
\draw[very thick] (0,0)..controls (0.5,0.5)..(1,0);
\draw[very thick] (0,1)..controls (0.5,0.5)..(1,1);
\end{tikzpicture}\,    \\
\label{eq.loop}
\begin{tikzpicture}[scale=0.8,baseline=0.3cm]
\fill[gray!20!white] (0,0)rectangle(1,1);
\draw[very thick] (0.5,0.5)circle(0.3);
\end{tikzpicture}
&=(-\xi^2 -\xi^{-2})
\begin{tikzpicture}[scale=0.8,baseline=0.3cm]
\fill[gray!20!white] (0,0)rectangle(1,1);
\end{tikzpicture}\, 
\end{align}

\def\Sq{\cS_q}
\def\SqM{\Sq(M)}
Very often  we identify a framed link $L$ with the element it represents in $\SM$. 

\no{An embedding  $f: M \embed M'$  of oriented 3-manifolds induces a $\cR$-linear homomorphism $f_*:\SM \to \SM'$, by $f_*(\al)= f(\al)$.}

Let $\fS$ be an oriented surface. Define
$$ \SS = \cS( \fS \times [-1,1]).$$
In this case the module $\SS$ has an algebra structure where the product $\al \beta$ of two framed links is given by placing $\al$ above $\beta$ in the direction given by the interval. 

Framed links in $\fS \times [-1,1]$ can be described by link diagrams on $\fS$ with blackboard framing. 
A link diagram on $\fS$ is {\bf simple} if it does not have crossings nor trivial loops. The set  $\B(\fS)$ of isotopy classes of simple diagrams is a free $\cR$-basis of $\SS$, see \cite{Prz}.

We use the notation $\Sq(M)$ for $\cS_q(M; \BC[q^{\pm 1}])$.  Following our notation, the formal variable $q$ is the variable used in the Kauffman bracket skein relations.
Let $D$ be a link diagram on $\fS$ with $\Cr$ the set of crossings.
As $\B(\fS)$ is a free $\Cq$-basis of $\SS$, in $\SS$ we have a presentation 
\be 
D = \sum_{\al \in S} C_\al \al, \ 0 \neq C_\al \in \Cq, S \subset \B(\fS), |S | < \infty.
\ee
\blem \label{rCross}
 We have $C_\al \in q^{e} \BZ[q^{\pm 2}]$, where $e\in \{0,1\}$ is determined by $e \equiv |\Cr| \mod 2$.
\elem
\bpr  Call the first and second figures on the right hand side of the skein relation \eqref{eq.skein}  respectively the positive and negative smoothings of the left hand side crossing. 
Using  skein and loop relations we have, in $\SS$, 
\begin{align}
D &= \sum_{\rho:\Cr \to \{ \pm 1\}} q^{\sigma(\rho)}  D _\rho, \quad \sigma(\rho)=\sum _{x \in \Cr} \rho(x),
\label{eq98aa}\\
&= \sum_{\rho:\Cr \to \{ \pm 1\}} q^{\sigma(\rho)} (-q^2-q^{-2})^{l(\rho)} (D _\rho)' \label{eq98ab},
\end{align}
where 
\begin{itemize}
\item $D _\rho$ is the result of $\rho(x)$-smoothing at all crossings $x$ of $D$,
\item $(D _\rho)'$ is the result of removing all trivial loops from $D_\rho$, and
\item $l(\rho)$ is the number of trivial loops in $D_\rho$.
\end{itemize}
Note that each $(D_\rho)'$ is in $\B(\fS)$. Since $\sigma(\rho) = |\Cr| \mod 2$, we have the lemma.
\epr

\no{

 Since the relations \eqref{eq.skein} and \eqref{eq.loop} preserve the $\BZ/2$-homology class, we have a $H_1(M;\BZ/2)$-grading of $\SM$:
\be \label{homgrading}
\SM = \bigoplus_{u\in H_1(M;\BZ/2)} \KzuM,
\ee
where $\KzuM$ is the subspace spanned by  framed links $\al\in M$   representing $u\in H_1(M;\BZ_2)$. 

\begin{definition}\label{Kzero} There is an action of $H^1(M;\mathbb{Z}/2)$ on $\SM$ described as follows. Given $c:H_1(M;\mathbb{Z}/2)\rightarrow \{\pm 1\}$, it acts as multiplication by $c(u)$ on  $\KzuM$. 
\end{definition}

Suppose $M = M_1 \sqcup M_2$. We identify $\cS(M_1) \ot_\ZQ \cS(M_2)$ with $\cS(M)$ via the isomorphism given by $\al_1 \ot \al_2 \to \al_1 \cup \al_2$, 
 for framed links $\al_i\in M_i$.
 }

\subsection{Barrett's isomorphism and its even analog}\label{ssEven}

 We now recall  Barrett's isomorphism and give an analog for the even part of the skein algebra.

 Given a spin structure $s$ of $M$, Barrett \cite{Barrett} showed that there is a $\Cq$-linear isomorphism $\sigma_s: \SqM \to \cS_{-q}(M)$ such that if $L$ is a framed oriented link, then
\be 
\sigma_s([L]) = s(L) [L],   \label{eqspin}
\ee
where $s(L) \in \{\pm 1\}$ is defined as follows. First consider 
$s$ as an element of 
$$H^1(\UM;\BZ/2)= \Hom(\pi_1(\UM , \{ \pm 1\}),$$ 
where $\UM$ is the $SO(3)$ principal bundle of $M$ associated to the tangent bundle. Equip $L$ with an orientation, then lift $L$ to $\tilde L \subset \UM$ using the framing, the orientation of $ L$, and the orientation of $M$. Then  define $s(L) = s(\tilde L')\in \{\pm 1\}$.

\def\See{{\cS_\ve(\fS)^\ev}}
\def\Se{{\cS_\ve}}
\def\SeM{\cS_\ve(M)}

\def\tord{{\widetilde{\ddd}}} 
 \def\Ctau{{\cC_\tau^{(2)}   }}
 \def\tbn{{\tilde\bn}}
 \def\tbt{{\tilde\bt}}
 \def\bj{{\mathbf j}}
 \def\bi{{\mathbf i}}
 \def\kq{\kappa_q}
 \def\SevS{\cS_q(\fS)^\ev}
 \def\SqS{\cS_q(\fS)}
 \def\bM{{\overline M}}
 \def\U{\mathsf U}
 
Let us turn to surface case and
define the even part of the skein algebra.

Assume $\fS= \Sigma_{g,m}= \Sigma_g \setminus\{v_1, \dots, v_m\}$. Let $\bfS= \Sigma_g$. For simple diagrams $\al, \beta$ on $\fS$ let $I(\al, \beta)$ be their geometric intersection. 
A simple diagram $\al$ is {\bf even} if $I(\al, \beta)\in 2\BZ$ for any  simple diagram $\beta$.  Equivalently, $\al$ represents the trivial homology class of $H_1(\bfS;\BZ/2)$. The {\bf even part}  $\SevS\subset \SqS$ is the submodule spanned by even diagrams. Since  defining relations \eqref{eq.skein} and \eqref{eq.loop} preserve the even part,  $\SevS$ is a subalgebra of $\SqS$. The subset $\B(\fS)^\ev \subset \B(\fS)$ of all even $\al\in \B(\fS)$ is a basis of $\SevS$.

\begin{theorem} \label{riq}
  
   (a) There is a unique $\Cq$-algebra isomorphism $\cS_{-q} (\fS)^\ev \to \cS_{q}(\fS)^\ev$ given by $\al \to \al$ for  
   $\al\in \B(\fS)^\ev$, i.e. $\al$ is an even simple diagram.

 (b) There is a unique $\Cq$-algebra isomorphism $ \cS_{\bi q} (\fS)^\ev \to \cS_{q}(\fS)^\ev$ given by $\al \to (-1)^{ \#\al} \al$ for  
   $\al\in \B(\fS)^\ev$ with $\#\al$ components. Here $\bi$ is the complex unit.

\end{theorem}

\bpr (a) Choose an arbitrary spin $s$ of $\bfS$ and restrict it to $\fS$. It is enough to show that $s(\al)=1$ for any even simple diagram $\al$, since then  \eqref{eqspin} implies $\sigma_s(\al) = \al$.

 Extend canonically $s$ to $\bM = \bfS \times [-1,1]$. Since $\al$ is even, it bounds a surface $S\subset \bfS$. Lift $S$ to $\tilde S\subset \U \bM$  using the orientation of $\bfS$ followed by the orientation of $[-1,1]$.  The boundary of  $\tilde S$ is the lift $\tilde \al$ of $\al$. Hence $\tilde \al$ is trivial in $H_1(\U \bM;\BZ/2)$; and $s(\al)=1$.
 
 (b) Fix a basis $y_1, \dots, y_{2g}$ of $H_1(\bfS;\BZ)$. Let $Q_{ij} =\la y_i, y_j\ra$, which is  the algebraic intersection number of $y_i$ and $y_j$. Let $\bT_{-1}(Q)$ be the quantum torus of matrix $Q$ and $q=-1$:
 $$ \bT_{-1}(Q)= \BC \la x_1^{\pm1}, \dots, x_{2g}^{\pm 1} \ra /( x_i x_j = (-1)^{ Q_{ij}} x_j x_i  ).$$
 We use $x_i$ to denote also its projection in the quotient $A := \bT_{-1}(Q)/( x_i^2=1)$. By \cite[Theorem 5.4]{FKL3}, there is an injective $\Cq$-algebra homomorphism
 $$ f:\cS_{\bi q}(\fS) \embed \cS_{q}(\fS) \ot_\BC A, \ f(\al) = (-1)^{ \# \al} \al \ot x^{ \bn(\al) }$$
 where $\bn(\al)\in \BZ^{2r}$ is the vector with entries $\bn(\al)_i = \la \al, y_i\ra$. 
 When $\al$ is even, $\bn(\al)_i$ is even, hence $x^{\bn(\al)} =1\in A$. Hence $f$ restricts to an algebra homomorphism $g: \cS^\ev _{\bi q}(\fS) \embed \cS_{q}^\ev(\fS)$. Since $g$ maps the basis $\B(\fS)^\ev$ of the domain bijectively to the basis $\B(\fS)^\ev$ of the codomain, it is bijective.
 \epr
\brem In \cite{FKL3}, which  extends the work of March\'e \cite{Marche}, $A$ is defined using a basis of $H_1(\fS;\BZ/2)$, not $H_1(\bfS;\BZ/2)$. But the proof there used only the algebraic intersection numbers, which are the same if we replace $H_1(\fS;\BZ/2)$ by $H_1(\bfS;\BZ/2)$. 
\erem

\subsection{Relation to character variety} We explain now a connection between  skein modules and  character varieties.

 For a non-zero complex number $\xi\in \BC$ we use $\SxS$ to denote $\cS_\xi(\fS;\BC)$.
 If $\ve\in \{\pm 1\}$ then $\SeM$ has a natural algebra structure where the product is given by disjoint union, see~\cite{Bullock}. Recall that  $\ChiM$ is the $SL_2(\BC)$ character variety of $M$.

\def\xar{\xrightarrow}

\begin{theorem}[See \cite{Bullock,PS2}]  \label{charcor} Let $M$ be a connected oriented 3-manifold.

(a) There is a surjective $\BC$-algebra homomorphism
\be 
T: \cS_{-1}(M)   \onto \BC[\ChiM], \ T(\al) = - T_\al,
\ee
where $\al\in \cS_{-1}(M)$ is a framed knot,  and $T_\al$ is defined by \eqref{eqCXM}. The kernel of $T$ is the nilradical of $\cS_{-1}(M)$.

(b) {\rm  (See  \cite{CM,PS2})} If $\fS=\Sigma_{g,m}$  then  $T$ is bijective.

\end{theorem}
 The map $T$  restricts to a $\BC$-algebra isomorphism 
 \be 
 T^\ev: \cS_{-1}(\fS)^\ev \to \BC[\X(\fS)^\ev].
 \label{eqTev}
 \ee

Using Theorem \ref{charcor} together with Barrett's isomorphism and its even analog we  define maps $\kappa^\bullet$ from  $\cS_{\ve}(\fS)^\bullet$ to $\BC[\X(\fS)^\bullet]$, where $\bullet=\emptyset$ or $\bullet=\ev$, as follows.

 For $\ve=\pm 1$ define $\kappa: \Se(M) \to \cS_{-1}(M)$ by
 \be 
 \kappa= \begin{cases} \Id  \quad & \text{if} \ \ve=-1\\
 \sigma_s & \text{if} \ \ve=1, \ \text{using a spin structure $s$.   }
 \end{cases}
 \label{eqkappa}
 \ee
 
For surfaces we will also consider even parts.
For $\ve\in \{\pm 1, \pm \bi  \}$, define 
$$\kappa^\ev: \See \to \cS_{-1}(\fS)^\ev$$ such that for an even simple diagram $\al\subset \fS$,
 \be 
 \kappa^\ev (\al) = \begin{cases} \kappa(\al)  & \text{if} \ \ve=\pm 1\\
 (-1)^{\# \al}\al  & \text{if} \ \ve=\pm \bi,  \ \text{where} \  \# \al =\text{number of components of } \ \al.
 \end{cases} \label{eqkappaev}
 \ee
By Theorem \ref{riq}  the map $\kappa^\ev$ is an algebra isomorphism.

 \def\bi{{\mathbf i}}
 \def\fr{{\mathrm{fr}}}
 \subsection{The Chebyshev-Frobenius map} We recall the Chebyshev-Frobenius homomorphism and its transparency property.
 
The $k$-th Chebyshev polynomial $T_k\in \BZ[z]$ of type 1 is defined so that 
$$
T_0(z)= 2, \ T_1(z)=z, \ T_k(z) = z T_{k-1}- T_{k-2} (z), \ k \ge 2.
$$
Alternatively, $T_k(z)$ is characterized by $T_k(x+ x^{-1}) = x^k + x^{-k}$.

Fix $\xi\in \BCx$.
For a framed knot $\al\subset M$ and an integer $k\ge 0$, write $\al^{(k)} \in \cS_\xi(M)$ for the {\bf $k$th framed power of $\al$} obtained by stacking $k$ copies of $\al$ in a small neighborhood of $\al$ along the direction of the framing of $\al$. Given a polynomial $P(z) = \sum c_i z^i \in \BZ[z]$, the {\bf threading of $\al$ by $P$} is  $P^{\fr}(\al):= \sum c_i \al^{(i)} \in \cS_\xi(M)$.

 \def\bi{{\mathbf i}}
 \def\Sx{{\cS_\xi}}
 \def\order{\ord}
 \def\MN{M}

\begin{theorem}[See \cite{BW2,Le:Frobenius} ]  \label{thmChFM}
Suppose $M$ is an oriented 3-manifold and $\xi$ is a complex root of unity. Let $N=\ord(\xi^4)$ and $\ve =\xi^{N^2}$. Also let $N' = \ord(\xi^2)$.
 
(a)  There exists a unique $\BC$-linear map $\Phi_\xi: \cS_\ve(M) \to \cS_\xi(M)$, called the {\bf Chebyshev-Frobenius homomorphism}, such that for any framed link $\al$ with components $\al_1, \dots \al_l$,
\begin{align*}
\Phi_\xi(\al) & = (T_N)^{\fr}(\al_1) \cup \cdots \cup (T_N)^{\fr}(\al_l) \\
&:= \sum_{0\le j_1, \dots, j_l\le N} c_{j_1} \dots c_{j_l}   \al_1^{(j_1)} \cup \cdots \cup \al_l^{(j_l)},\  \text{where} \ T_N(z) = \sum_ {j=0}^N c_j z^j.
\end{align*}

(b)  The image of $\Phi_\xi$ is $(-1)^{N'+1}$-transparent in $\cS_\xi\MN$ in the sense that
\be  \raisebox{-14pt}{\incl{1.2 cm}{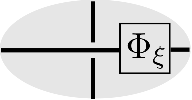}}  =  (-1)^{N'+1}\, \raisebox{-14pt}{\incl{1.2 cm}{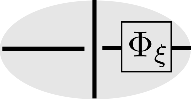}}. 
\label{eqTrans}
\ee
Here the box enclosing $\Phi_\xi$ means one applies $\Phi_\xi$ to the component of the link containing the horizontal line, while the vertical line belongs to another component of the link.


\end{theorem}
\bpr 
  Part (a), and also part (b) with  $N'$ odd, were proved by \cite{BW2} for  surfaces (and implicitly for 3-manifolds). For 3-manifolds, and all roots in part (b), the theorem is proved in \cite{Le:Frobenius}. See also \cite{LP,BL}.
\epr

\def\Se{{\cS_\ve}}
\def\ZZ{{ \mathfrak U}_\xi(\fS)}
\def\tZZ{{ \tilde{\mathfrak U}_\xi(\fS)}}
\def\tPhi{\tilde \Phi}
\def\VV{{\mathsf V}}
\def\SMz{{ \cS_\xi(M)^0}}

\subsection{Action of $\Se(M)$ on $\Sx(M)$ }  
\label{ssAction} Let $M$ be a connected oriented 3-manifold and $\xi$ be a root of 1 with 
 $N'= \ord(\xi^2)$ odd. Then $(-1)^{N'+1}=1$ and $\ve:= \xi^{N^2} \in \{ \pm 1\}$. Hence $\Se(M)$ has the structure of a commutative algebra.

Define an action of $\Se(\fS)$ on $\Sx(\fS)$ as follows. Suppose $\al$ is framed link, considered as an element of $\Se(M)$, and $\beta$ is a framed link, considered as an element of $\Sx(\fS)$. Define
\be  \al * \beta = \Phi_\xi(\al) \cup \beta. \label{eqActiono}
\ee
 The transparency \eqref{eqTrans} implies that the right hand side depends only on the  isotopy classes  of $\al$ and $\beta$. Linearly extending this action,  we get an action of $\Se(M)$ on $\Sx(M)$.
 
 \no{
 Assume now $N'= \ord(\xi^2)$ is even, and further assume that $M$ is closed. Then  $\ve:= \xi^{N^2} \in \{ \pm 1, \pm \bi\}$.
 Let $\SMz\subset \SxM$ be the $\BC$-subspace spanned by framed links $\al$ which are 0 in $H_1(M; \BZ/2)$. Let $\al, \beta\in M$ be links representing 0 in $H_1(M; \BZ/2)$. Define $I_2(\al, \beta) \in \BZ/2$ by 
 $$ I_2(\al, \beta)  \equiv |\al \cap S| \mod 2,$$
 where $S$ is any surface, not necessarily oriented, bounded by $\beta$ and transverse to $\al$. Poincare duality shows $I_2(\al, \beta)$ does not depend on the choice of $S$.  Now consider $\al$ as an element of $\cS_\ve(M)^0$ and $\beta$ as an element of $\cS_\xi(M) ^0$, define
 \be \al * \beta = (-1)^{I_2(\al,\beta)}\Phi_\xi(\al) \cup \beta.
 \ee
 Identity \eqref{eqTrans} and the behavior of $I_2$ show that this gives a well-defined action of $\cS_\ve(M)^0$ on $\cS_\xi(M)^0$.
 
Note that for a surface $\fS=\Sigma_{g,m}$ with $m >0$  we have $\cS_\xi(\fS)^0 \subsetneq \cS_\xi\fS)^\ev$.
}

\def\SxS{{\Sx(\fS)}}
\def\See{\cS^\ev_\ve(\fS)}
\def\SeS{\cS_\ve(\fS)}
\def\Seb{{\cS^\bullet_\ve}}
 \subsection{The center of the skein algebra} \label{ssCenter} 
 We now describe the center $Z(\SxS)$  of the skein algebra $\SxS$ when $\xi$ is a root of 1. The following was proved in \cite[Section 4]{FKL}.
 
 \bthm \label{thmCent1}
 Let $\fS= \Sigma_{g,m}$ with peripheral loops $\fd_1, \dots, \fd_m$, and let $\xi\in \BC$ be a root of 1, with $N'= \ord(\xi^2)$ and $N= \ord(\xi^4)$. Let $\ve:= \xi^{N^2}$.  For any object $A$ let $A^\bullet=A$ if  $N'$ is odd, and $A^\bullet=A^\ev$ if $N'$ is even.
 
 There is a $\BC$-algebra isomorphism
\be 
\tPhi_\xi^\bullet: \Seb(\fS) [d_1,\dots, d_m]/(T_N(d_i) = \fd_i, i=1, \dots, m)  \xar{\cong} Z(\SxS) 
 \ee
such that
$
\tPhi_\xi^\bullet (x) = \Phi_\xi(x)$ for $ x\in \Seb(\fS)$, and 
$
\tPhi_\xi^\bullet (d_i)= \fd_i, i=1, \dots, m.
$

\ethm

\no{ 
 \bthm \label{thmCent0}
 Let $\fS= \Sigma_{g,m}$ with peripheral loops $\fd_1, \dots, \fd_m$, and $\xi\in \BC$ is a root of 1, with $N'= \ord(\xi^2)$ and $N= \ord(\xi^4)$. Let $\ve:= \xi^{N^2}$.  For any object $A$ let $A^\bullet=A$ if  $N'$ is odd, and $A^\bullet=A^\ev$ if $N'$ is even.
 
 There is a $\BC$-algebra isomorphism
\be 
\tPhi_\xi^\bullet: \Seb(\fS) [d_1,\dots, d_m]/(T_N^\bullet(d_i)) = \fd_i^\bullet, i=1, \dots, m)  \xar{\cong} Z(\SxS) 
 \ee
such that
$
\tPhi_\xi^\bullet (x) = \Phi_\xi(x)$ for $ x\in \Seb(\fS)$, and 
$
\tPhi_\xi^\bullet (d_i)= \fd_i, i=1, \dots, m.
$

\ethm

}

We can now use the isomorphism $\kappa^\bu$ defined by Equations (\ref{eqkappa}) and (\ref{eqkappaev}) to relate the center $Z(\SxS)$ to the character variety. Consider the composition of  isomorphisms
 \be 
\begin{split}
&\Psi: Z(\SxS) \xra{ (\tPhi_\xi^\bullet)^{-1} } \Seb[ d_1, \dots, d_m]/ ( T_N(d_i)=\fd_i) ) 
\xra{\kappa^\bullet}\\
 &\cS_{-1}^\bullet(\fS) [ d_1,\dots, d_m]/(T_N(d_i) = \fd_i) \xra{T^\bullet} \BC[ \X^\bullet(\fS)] [ d_1,\dots, d_m]/(T_N(d_i) = - T_{\fd_i}).
\end{split}
\label{eqPsi}
\ee
For one case, namely when $\ve=1$, we need a spin structure $s$ of $\fS$ to define $\kappa$. To have a uniform treatment for all cases, we choose $s$ such that $s(\fd_i)=-1$ for all $i$. For example, $s$ can be any spin induced from a spin of $\bfS$.

Let  $\VV_\xi(\fS):= \MaxSpec(\SxS)$  be the classical shadow variety  of $\SxS$. The dual of the isomorphism $\Psi$ gives an identification
\be  \VV_\xi(\fS) = \{([\rho], \vw) \in \X^\bullet(\fS) \times \BC^m \mid  T_N(\w_i) =- \tr(\rho(\fd_i)).  \}
\label{eq26}
\ee
Taking the projection onto the first component, we get a finite morphism of degree $N^{m}$
\be  p_\xi: \VV_\xi(\fS)  \to \X^\bullet(\fS).
\ee
Note that $p_\xi$ is Poisson, because $\BC[\X(\fS)^\bullet] \embed \BC[ \X^\bullet(\fS)] [ d_1,\dots, d_m]/(T_N(d_i) = - T_{\fd_i})$ is Poisson, since $d_i$'s are Casimir.
\no{ The projection  $\X^\bullet(\fS)\times \BC^m \to \BC^m$ restricts to a regular morphism
\be 
\pr_2: \VV_\xi(\fS) \to \BC^m.
\ee
}

\subsection{Comparison of Poisson structures} We  show that the isomorphism $\Psi$ of \eqref{eqPsi} is Poisson.

The domain  $Z(\SxS)$ of $\Psi$  has the quantization Poisson structure given by Lemma \ref{rPoisson}. The codomain has 
the Poisson structure, which  
extends the ABG-Poisson structure from $\BC[\X^\bullet(\fS)]$ to $\BC[ \X^\bullet (\fS)] [ d_1,\dots, d_m]/(T_N(d_i) = - T_{\fd_i})$ so that $d_i$ is Casimir for $i=1,\dots, m$.

\def\ZxS{{Z(\SxS)}}
\def\tq{{\theta_q}}
\def\tx{{\theta_\xi}}
\def\tz{{\theta_\zeta}}
\def\TqO{\bT_q(Q;\Lambda)}
\def\TxO{\bT_\xi(Q;\Lambda)}
\def\TeO{\bT_\ve(Q;\Lambda)}
\def\hy{{\hat y}}
\def\by{{\bar y}}
\def\bal{{\mathrm{bl}}}
\bpro \label{rMaxSpec}
 The isomorphism $\Psi$ defined in \eqref{eqPsi} is Poisson.  
\no{
(b) For $\vw \in \BC^m$ the level surface $(\pr_2)^{-1}(\vw)$ is Poisson isomorphic to the sliced character $\X^\bullet (\fS; - T_N(\vw))$, and consequently is an irreducible symplectic variety with singularity.
}
\epro

\def\SqS{{\cS_q(\fS)}}

\bpr The idea is simple. We use the quantum trace map to reduce the general case to the case $\xi\in \{\pm 1, \pm \bi\}$, which is further reduced to the case $\xi=-1$ via the map $\kappa$. Finally, the case $\xi=-1$ was done in \cite{BFK}. Here are the details.

(a) As $\Psi= T^\bullet\circ \kappa^\bullet\circ (\tPhi_\xi^\bullet)^{-1} $, we will show that each of $T^\bullet, \kappa^\bullet, \tPhi_\xi^\bullet$ is Poisson.


{\bf Step 1.} \  $T^\bullet: \cS_{-1}^\bullet(\fS) \to \BC[\cS^\ev]$ is Poisson. For $N'$ odd this was proved in \cite{BFK}. When $N'$ is even, the map $T^\ev$ is a restriction of $T$, and hence it is Poisson.

  {\bf Step 2.} Let us prove $\kappa^\bullet$ is Poisson. \begin{enumerate}
   \item If $\ve=-1$  then $\kappa^\bullet=\Id$ and the statement is trivial.
  \item If $\ve=1$ then $\kappa^\bullet=\sigma_s$ or its restriction onto the even part. The statement follows from Proposition \ref{rPoiMor}, as $\sigma_s$ can be defined as an isomorphism for generic $q$,  from $\SqS$ to $\cS_{-q}(\fS)$.  
  \item If $\ve=\pm \bi$ then $\bullet= \ev$. The map $\kappa^\ev$ can be defined for generic $q$ by Proposition \ref{riq}. Hence we get the statement from Proposition \ref{rPoiMor}.
  \end{enumerate}

{\bf Step 3.} The map $\tPhi_\xi^\bullet: \Seb[ d_1, \dots, d_m ] /( T_N(d_i) = \fd_i  ) \to \ZxS$ is Poisson.

 The proof is much more involved. We will use the quantum trace, see \cite{BW}, and also  \cite{Le:triangular,CL}, to embed $\SxS$ into the quantum torus, where the calculation is easier.

Recall that $\tPhi_\xi^\bullet$ extends $\Phi_\xi^\bullet$ where the extension values are $\tPhi(d_i)= \fd_i$. Since $d_i$ and $\fd_i$ are Casimir elements of the Poisson brackets, we need only to prove that $\Phi_\xi^\bullet$ is Poisson. 
Denote $Z_0= \Phi_\xi^\bullet(\Seb)$. We consider two cases: $m \ge 1$ and $m=0$.

 {\bf Case 1: $m\ge 1$}.
Let $r:= 3g-3+3m$. 
 \blem \label{rQXtrace}
 There exist  an antisymmetric $r\times r$-matrix $Q$,  subgroups $\Lambda^\ev \subset \Lambda$ of $\BZ^r$ of rank $r$, and a $\Cq$-algebra embedding $\tq: \SqS \embed \TqO$, where $\TqO$ is the monomial algebra defined in Subsection \ref{ssQtorus} with ground ring $\Cq$, with the following properties.
 
(a) For  $\zeta\in \BCx$, the map $\tz$ is injective. Besides, the following diagram is commutative
 \be
\begin{tikzcd}
\SeS \arrow[r,hook, "\theta_\ve"]
\arrow[d,hook,"\Phi_\xi"]  
&  \mathbb{T}_\ve(Q;\Lambda)\arrow[d,hook,"\Phi_N"] \\
\SxS\arrow[r,hook,"\theta_\xi "] & \mathbb{T}_\xi(Q;\Lambda)
\end{tikzcd}\qquad \text{where $\Phi_N(x_i) = x_i^N$.}
\label{eqDiag1}
\ee

(b)  Diagram \eqref{eqDiag1} restricts to 
 \be
\begin{tikzcd}
\Seb \arrow[r,hook, "\theta_\ve"]
\arrow[d,hook,"\Phi_\xi"]  
&  \mathbb{T}_\ve(Q ;\Lambda^\bullet)\arrow[d,hook,"\Phi_N"] \\
Z_0\arrow[r,hook,"\tx"] & \bT_\xi(Q ;N \Lambda^\bullet)
\end{tikzcd}
\label{eqDiag12}
\ee
Moreover 
\begin{align}
 \mathbb{T}_\ve(Q ;\Lambda^\bullet) &\subset Z(  \mathbb{T}_\ve(Q ;\Lambda)  ),
 \label{eqincl1} \\
  \mathbb{T}_\xi(Q ;N \Lambda^\bullet) &\subset Z(  \mathbb{T}_\xi(Q ;\Lambda)  ).  \label{eqincl2}
\end{align}


\elem
 The proof, based on the quantum traces, is given later. 
 

 By Lemma \ref{rQXtrace}(a) and Proposition \ref{rPoiMor}, the two horizontal maps of Diagram \eqref{eqDiag12} are Poisson. Let us prove the right vertical one is Poisson. The quantum torus $\bT_q(Q ;\Lambda^\bullet)$ has a presentation
 $$ \bT_q(Q ;\Lambda^\bullet) = \Cq\la \hy_1^{\pm 1}, \dots, \hy_r^{\pm1}\ra /(\hy_i \hy_j = q^{P_{ij}}\hy_j \hy_i)$$
 for a certain antisymmetric $r\times r$ matrix $P$. Let $y_i=\pi_\ve(\hy_i)\in\mathbb{T}_\ve(Q ;\Lambda^\bullet) $, where $\pi_\ve: \bT_q(Q ;\Lambda^\bullet)\onto
 \bT_\ve(Q;\Lambda^\bullet ) $ is the natural projection. 
 Denote $\{ \cdot, \cdot \}_\ve$ the quantization Poisson bracket in $Z(\mathbb{T}_\ve(Q;\Lambda^\bullet ))$, which contains $\bT_\ve(Q ;\Lambda^\bullet)$. 
 Since $y_i y_j = y_j y_i$, we have that $q^{P_{ij} }-1$ is divisible by $q-\ve$ in $\Cq$.
  L'Hopital's rule shows that
 $$ \pi _\ve \left ( \frac{q^{P_{ij}}-1}{q-\ve}\right)= \ve^{-1}P_{ij}.   $$
 From the definition,
 \be 
 \{ y_i, y_j\}_\ve 
  =\pi_\ve  \left ( \frac{[\hy_i, \hy_j]}{q-\ve}\right) 
  =\pi _\ve \left ( \frac{q^{P_{ij}}-1}{q-\ve}  \hy_j\hy_i\right) 
  = \ve^{-1}P_{ij} y_j y_i .\label{eqPoi5}
 \ee
 Let us calculate the quantization Poisson bracket  $\{ \cdot, \cdot \}$ of $\bT_\xi(Q ;N \Lambda^\bullet)$. Let $\by_i= \pi_\xi(\hy_i)\in \bT_\xi(Q ;\Lambda^\bullet)$.  Then $\bT_\xi(Q ;N\Lambda^\bullet)$ is generated by $(\by_i)^{\pm N}$.
 The same calculation shows
 \be 
 \{ \by_i^N, \by_j^N\}_\xi 
  =\pi_\xi  \left ( \frac{[\hy_i^N, \hy_j^N]}{q-\xi}\right) 
  =\pi _\xi \left ( \frac{q^{N^2P_{ij}}-1}{q-\xi}  \hy_j^N\hy_i^N\right) 
  = \xi^{-1}N^2 P_{ij} \by_j \by_i .\label{eqPoi6}
 \ee
 As $\Phi_N(y_i) = (\by_i)^N$, from \eqref{eqPoi5} and \eqref{eqPoi6} we have
 $$ \Phi_N (\{ \cdot, \cdot\}_\ve   )  = (\ve \xi^{-1} N^2) \{ \cdot, \cdot\}_\xi ,$$
 proving that $\Phi_N$ is Poisson. In Diagram \eqref{eqDiag12}, the two horizontal maps and the right vertical one are Poisson. It follows that the left vertical one, which is $\Phi_\xi$, is Poisson. 
 
 {\bf Case 2. $m=0$.}
 By removing one puncture from $\fS$ we get $\fS'$ which is triangulable. The embedding $\iota: \fS' \embed \fS$ induces a surjective $\Cq$-algebra homomorphism $\iota_q: \cS_q(\fS') \onto \SqS$. Hence $\iota_\xi$ is Poisson. The lemma for $\fS$ follows from  the case of $\fS'$. 
 \epr

\bpr[Proof of Lemma \ref{rQXtrace}]Choose an ideal triangulation of $\fS$, with the set of edges $\cE= \{ e_1, \dots, e_r\}$. We can assume that for $i=1, \dots, 2g$ each $e_i$ is a loop in $\Sigma_g$, and that $\{e_1, \dots, e_{2g}\}$ is a basis of $H_1(\Sigma_g, \BZ)$.

Let $Q$ be the $r\times r$ matrix defined as in \cite{LY2}, which is equal to the matrix $\sigma_{ij}$ of \cite{BW}. The explicit formula of $Q$ is not important for us. For a simple diagram $\al\in \B(\fS)$ let $\bn(\al) \in \BZ^r$ be the vector defined by $\bn(\al)_i = I(\al, e_i)$. 
Let $\Lambda\subset \BZ^r$ be the subgroup generated by all $\bn(\al), \al \in B(\fS)$, and $\Lambda^\ev\subset \BZ^r$ be the subgroup generated by all  $\bn(\al)$ with even $\al \in \B(\fS)$. Then
 $\Lambda\subset \BZ^r$ is the subset of $\bk= (k_1, \dots, k_r)$ such that $k_i+k_j+k_l\in 2\BZ$ whenever $e_i, e_j, e_l$ are edges of a triangle. The quantum trace \cite{BW} is an algebra embedding  $\theta_q: \SqS  \embed \bT_q(Q;\Lambda)$. 

(a) By  \cite[ Theorems 21 \& 29]{BW}, $\theta_\zeta$ is injective and Diagram \eqref{eqDiag1} is commutative.

(b)  We begin with a lemma, which is of independent interest.
\blem (i) If $\al\in \B(\fS)$ is even, then $\theta_q(\al) \in \bT_q(Q; \Lambda^\ev)$.

(ii) We have $\la \Lambda, \Lambda \ra_Q \subset 2 \BZ$ and $\la \Lambda^\ev, \Lambda \ra_Q \subset 4 \BZ$. 

(iii) $\bT_\xi(Q;N \Lambda^\bullet)$ is in the center of $\bT_\xi(Q;\Lambda)$.
\elem
\bpr (i) Note that $\al\in \B(\fS)$ is even if and only if it represents 0 in $H_1(\Sigma_g;\BZ/2)$, or equivalently $\bn(\al)_i\in 2 \BZ$ for $i=1, \dots, 2g$.

By \cite[Lemma 6.11 ]{LY2},  any monomial $x^\bk$ appearing in $\theta_q(\al)$ satisfies $\bk - \bn(\al) \in (2\BZ)^r$. Hence $\bk_i\in 2 \BZ$  for $i=1, \dots, 2g$. Thus $\bk\in \Lambda^\ev$, and $\theta_q(\al)\in \bT_q(Q;\Lambda^\ev)$.

(ii) For $\al, \beta \in \B(\fS)$ we have $\la \bn(\al), \bn(\beta)\ra \in 2\BZ$, see \cite[Remark 3.11]{FKL}. As $\Lambda$ is spanned by $\bn(\al)$, we have $\la \Lambda, \Lambda \ra_Q \subset 2 \BZ$. By \cite[Proposition 3.14]{FKL},
$$
I(\al,\beta) \equiv \frac 12 \la \bn(\al), \bn(\beta) \ra_Q \mod 2.
$$
Hence if $I(\al, \beta)\in 2\BZ$ we have $\la \bn(\al), \bn(\beta) \ra_Q \in 4\BZ$. It follows that $\la \Lambda^\ev, \Lambda \ra_Q \subset 4 \BZ$.

(iii) By  \eqref{eq.prod} we have
$ x^\bk x^\bl = \xi^{\la \bk, \bl\ra_Q} x^\bl x^\bk$. Let $\ord(\xi)= N''$. It is enough to show that $\la N \bk ,\bl\ra_Q \subset N'' \BZ$ for $\bk\in \Lambda^\bullet$ and $\bl\in \Lambda$.

Note that  $N= \frac{N'}{gcd(N', 2)}= \frac{N''}{gcd(N'', 4)}$, and $N'= \frac{N''}{gcd(N'', 2)}$.

Assume $N'$ is odd.  Then $N''= 2N'= 2N$, and $\Lambda^\bullet=\Lambda$. From (ii) we have
$$\la N \Lambda ,\Lambda\ra_Q \subset  N (2  \BZ) = N'' \BZ.$$
Assume $N'$ is even. Then $N''= 4 N$. From (ii) we have
$$\la N \Lambda^\ev ,\Lambda\ra_Q \subset  N (4 \BZ) = N'' \BZ. $$
This completes the proof of the lemma.
\epr
Let us continue with the proof of Lemma \ref{rQXtrace}(b). By (i) we have
$$ \theta_\ve(\Seb )\subset  \mathbb{T}_\ve(Q ;\Lambda^\bullet).$$
The target space is mapped to $\mathbb{T}_\xi(Q ;N\Lambda^\bullet)$ by $\Phi_N$. Hence the commutativity of Diagram \eqref{eqDiag1} shows that we have the restriction as in Diagram \eqref{eqDiag12}.

By (iii) we have \eqref{eqincl1}. If we replace $\xi$ by $\ve$, then $N$ would be 1, and we get \eqref{eqincl2}.
\epr

\brem Proposition \ref{rMaxSpec}, for the case $m=0$ and the order of $\xi$ is odd, was proved in \cite{GJS}, using a similar (but still different) method. The even order case makes the proof much more technical.
\erem

\def\uSS{{\underline{\cS}(\fS)}}
\def\Bp{{B^+}}
\def\uBp{\uB^\uparrow}

\def\utr{{\underline{\tr}}}
\def\bbl{{\mathbf l}}
\def\mP{{\mathsf P}}
\def\uD{{\underline{\Delta}}}
\def\uA{{\underline{\cA}}}
\def\uP{{\underline{\mP}}}
\def\uv{\overleftarrow }
\def\cb{{\check b}}
\def\pPPb{\partial \PPb}

\def\uS{{\underline{\cS}}}

\newcommand{\reldown}[2]{
\begin{tikzpicture}[scale=0.8,baseline=0.3cm]
\fill[gray!20!white] (0,0)rectangle(1,1);
\draw[-stealth] (1,1)--(1,0);
\draw[very thick] (0,0.67)--(1,0.67) (0,0.33)--(1,0.33);
\draw[inner sep=1pt,right] (1,0.67)node{\tiny #1} (1,0.33)node{\tiny #2};
\end{tikzpicture}
}

\def\pPPb{\partial\PPb}

\def\Bp{{B^+}}

\newcommand{\relcross}[2]{
\begin{tikzpicture}[scale=0.8,baseline=0.3cm]
\fill[gray!20!white] (0,0)rectangle(1,1);
\draw[-stealth] (1,0)--(1,1);
\begin{knot}[clip width=4,background color=gray!20!white]
\strand[very thick] (0,0.3)--(1,0.67);
\strand[very thick] (0,0.7)--(1,0.33);
\end{knot}
\draw[inner sep=1pt,right] (1,0.67)node{\tiny #1} (1,0.33)node{\tiny #2};
\end{tikzpicture}
}

\section{Sliced  skein algebra}
\label{secMRY}

We define the sliced skein algebra in a more general setting. Namely we will involve the boundary in the definition. This extension is needed later in the proof of Theorem  \ref{thm12}.


\subsection{A quick definition for $\Sigma_{g,m}$} 
Assume $\fS= \Sigma_{g,m}$, with peripheral loops 
$\fd_1, \dots, \fd_m$. Let the ground ring $\cR$ be a commutative domain, with invertible $\xi \in \cR$ and $\vw=(\w_1, \dots, \w_m) \in \cR^m$. The sliced skein algebra is
$$ \SslR := \SSR/(\fd_i = \w_i). $$

It turns out that technically it is more convenient to work with surfaces having boundary, and to consider interior punctures as formal variables. We will extend the notion of sliced skein algebra in this direction.
\subsection{Punctured bordered surfaces} In this subsection we define  punctured bordered (pb) surfaces and tangle diagrams on them.
\def\pal{{\partial \al}}

\bdf 
\label{defpb}
A {\bf punctured bordered (pb) surface} is a surface of the form $\mathfrak{S}=\bfS\backslash \cM$ where $\overline{\mathfrak{S}}$ is a compact oriented surface with possibly empty boundary $\pbfS$ and $\cM$ is a  finite set whose elements are called {\bf punctures}, such that each connected component of $\pbfS$ has at least one puncture. Thus
 each component of the boundary $\pfS$ is an open interval and is called a {\em boundary edge}.

 A {\bf $\pfS$-arc} is a smooth  proper embedding $a: [0,1] \embed \fS$. A $\pfS$ arc $\al$ is {\bf near boundary} if, as a map into $\bfS$, it can be homotoped relative its endpoints into $\pbfS$.


A {\bf loop on $\fS$} is  a simple closed curve on $\fS$. A loop is {\bf trivial} if  it bounds a disk in $\fS$. A loop is {\bf peripheral} if  it bounds a punctured disk, i.e. a disk in $\bfS$ containing exactly one point of $\cM$. 
\edf
Note that $\bfS$ can be  uniquely recovered from $\fS$.

From now until the end of the section fix a pb surface $\fS$. 
\bdf (a) A {\bf $\pfS$-tangle diagram} is a proper immersion $\al: C \to \fS$, from a compact 1-dimensional manifold into $\fS$ having  only double point singularities, called crossings, and equipped with  an over/undercrossing information of the two strands at every crossing like in an ordinary knot diagram. Denote by $\pal=\al \cap \pfS$ the set of boundary points of $\al$.

(b) A  $\pfS$-tangle diagram $\al$ is  {\bf boundary ordered}, if it is equipped with a linear order on each set $\pal \cap b$, for each boundary edge $b$. The boundary order is also called a height order. Together the height orders define  a partial order on the set $\pal=\al \cap \pfS$ of endpoints of $\al$. Two points $x,y\in \pal$ are consecutive if there is no $z\in \pal$ such that $ x < z <y$.

(c) The empty set is considered a $\pfS$-tangle diagram.

(d) Isotopies of $\pfS$-tangle diagrams are ambient isotopies of $\fS$.

\edf


\def\tfS{ { \widetilde{\fS}  }}

\def\ptfS{ {\partial { \tfS }}}

\def\ori{{\mathfrak o}}

A boundary edge $b$ inherits an orientation from the orientation of $\fS$. 
A boundary order of $\al$ is {\bf positive} if   on a boundary edge $b$ the height order  is  increasing when one goes along $b$ in  the direction of the orientation of  $b$.

\brem  One should consider a $\pfS$-tangle diagram as the diagram of a framed tangle in the thickened surface $\fS 
\times (-1,1)$. 

\erem

\subsection{Universal sliced skein algebra } \label{ssMRY} We introduce the  {\bf universal sliced skein algebra}. 

Let $\Mo$ be the set of all interior punctures of $\fS$. The {\bf universal ground ring} is
\be
\RS = \Zq[ v^{\pm 1}, v \in \Mo],
\ee
which is the ring of Laurent polynomial in variables  $v\in \Mo$ with coefficients in $\Zq$.  For $v\in \Mo$ let $\fd_v$ be the peripheral loop surrounding $v$. 

\bdf 
\label{defslSkein}
The universal sliced skein algebra $\SslS$ is the  $\RS$-module freely spanned by isotopy classes of boundary ordered $\pfS$-tangle diagrams  subject  relations \eqref{eq.skein} and \eqref{eq.loop}, and 
\begin{align}
\al&=0 \text{ for any $\pfS$-tangle diagram $\al$  having a near boundary arc}
\label{eqNearB} \\
\pv &= v + v^{-1}  \label{eqPeri}\\
\relcross{}{}&= \reldown{}{} \label{eqisotopy}
\end{align} 
The product is defined by stacking as in the case of the ordinary skein algebra.
\edf 
Here is the convention used in relation \eqref{eqisotopy}, as well as in many other figures of the paper: In the shaded square we have a part of the diagram, having two endpoints there. The directed line is a part of a boundary edge. The order of the two endpoints is given by the arrow: i.e.,  going along 
the line increases the height. The two depicted endpoints are consecutive in the height ordered. Thus relation \eqref{eqisotopy} is simply the isotopy in the thickened surface.


\bpro Considered as a $\Zq$-algebra, $\SslS$ has a reflection $\omega$ defined so that $ \omega(q^{1/2})= q^{-1/2}, \omega(v)= v$ for all $v \in \Mo$, and for a boundary ordered tangle diagram $D$ the image $\omega(D)\in \SslS$ is obtained from $D$ be switching all crossing and reversing the height order on each boundary edge.
\epro
\bpr 
It is easy to check that $\omega$ respects all the defining relations of $\SslS$.
\epr

\brem \label{remuSS}
 If we drop the relation \eqref{eqPeri} then we get the {\bf boundary simplified skein algebra} $\uSS$ of \cite{Le:Qtrace,BKL}. When there is no interior puncture, $\SslS$ the reduced skein algebra studied in \cite{Le:Qtrace,PS2}, and it is a quotient of  the Muller skein algebra \cite{Muller}
\erem

\def\BCx{{\BC^\times}}

\subsection{Height exchange, basis}
Changing the height order results in a $q$-factor.
\blem[Height exchange rule ] \label{rHeight}
   In $\SslS$ one has
\be 
\begin{array}{c}\includegraphics[scale=0.7]{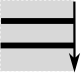}\end{array} = q^{-1} \begin{array}{c}\includegraphics[scale=0.7]{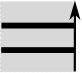}\end{array}
\label{eqHeight}
\ee
\bpr The algebra $\SslS$ is a quotient of the positive stated skein algebra (or the Muller skein algebra), for which \eqref{eqHeight} holds, see  \cite[Proposition 2.4]{Le:triangular}.
\epr
\elem

Suppose a boundary ordered $\pfS$-tangle diagram $\al$ is the disjoint union of its connected components. As the reflection $\omega(\al)$ differs from $\al$ by height order, Lemma \ref{rHeight} shows  that there is a unique $k\in \BZ$ such that $\omega(\al) = q^k \al$ in $\SslS$. Define
 the reflection normalization
\be 
[ \al]_\omega = q^{k/2} \al. \label{eqRefl}
\ee
The normalized $[\al]_\omega$ is reflection invariant.

\no{ 

It is easy to check, though we don't need it in the sequence, that
\be 
[ \al]_\omega = q^{ - \sum \frac{\bn_\al(e) (\bn_\al(e) -1)}{2} } (\al, h^+),
\ee
where $(\al, h^+)$ is $\al$ with the positive boundary order, and the sum is over all boundary edges. 
}

\def\Bsl{{\mathsf B^\sli}}
\def\vBsl{{\vec{\mathsf B}^\sli}}
 
 \bdf\label{defSsimple}
  A tangle diagram $\al$ on $\fS$ is {\bf strongly simple} if it does not have any trivial loop, any peripheral loop,  nor any near boundary arc. Let $\Bsl(\fS)$ denote the set of all isotopy classes of strongly simple $\pfS$-tangle diagrams on $\fS$. Let $\vBsl(\fS)$ be the same set $\Bsl(\fS)$, where each diagram is given the positive boundary order.
 \edf

\bpro
 \label{rslBasis}The set $\vBsl(\fS)$ is a free $\RS$-basis of $\SslS$. In particular when $\pfS= \emptyset$ the set $\Bsl(\fS)$ is a $\Zq$-basis of $\SslS$.
\epro
\bpr Let $\uSS$ be the algebra defined as $\SslS$, but without the peripheral relation \eqref{eqPeri}, see Remark \ref{remuSS}. This is known as the boundary simplified skein algebra in \cite{BKL}. By  \cite[Proposition 6.10]{BKL} the set 
 $\uB(\fS)$  of all isotopy classes of simple diagrams without near boundary arcs is a free $\RS$-basis of $\uSS$. 
 Note that if $\al \in B^\sli$ then $(\prod_{i=1}^m \fd_{v_i} ^{k_i}) \al$ can be realized as the disjoint union of $\al$ and $k_i$ copies of $\fd_{v_i}$ which are in a small neighborhood of $v_i$,
  and hence it is an element of $\uB(\fS)$. 
  It follows that $\uSS$ has $\B^\sli(\fS)$ as a basis over   $\RS[\fd_{v_1}, \dots, \partial _{v_m}]$. Changing the ground ring,  we see that
 $$ \SslS = \uSS / ( \fd_{v_i} = v_i + v_{i}^{-1})=\uSS \ot_ { \RS[\fd_{v_1}, \dots, \partial _{v_m}]  } \RS$$
 has $\B^\sli(\fS)$ as a basis over $\RS$. Here $\RS$ is considered as an 
 $\RS[\fd_{v_1}, \dots, \partial _{v_m}]$-algebra, with the module structure given by $ \fd_{v_i} \to v_i+ v_i^{-1}$.
\epr

\def\vuB{{\overrightarrow{\uB}}}

\def\pPPb{\partial \PPb}

\section{Dehn-Thurston coordinates and modification}

\label{secDT} The set $\Bsl(\fS)$ of isotopy classes of strongly simple diagrams on the  surface 
$\fS=\Sigma_{g,m}$ is a basis of $\SslS$, and in this
section we will parameterize this basis by a submonoid of $\BZ^{2r}$, where $r= 3g-3+m$. This parametrization is a modified version of
 the standard Dehn-Thurston (DT) coordinates. The standard DT-coordinates do not behave well under the skein product. That is the reason why we need the modified version, which will be shown later to capture  the highest degree term in the product of skein algebra. 

We will first construct our DT-coordinates for elementary pieces, called the basic DT pairs of pants. Then we combine them together to get the global DT-coordinates for $\Bsl(\fS)$.

Throughout the section we fix $\fS=\Sigma_{g,m}$ with  exceptions
\be (g,m) \not \in \{ (1,0), (0,k), k\le 3 \} 
\ee
This is equivalent to $r:= 3g-3 + m >0$.

\subsection{Three basic DT pairs of pants, DT-datum}

\bdf

(a) The {\bf DT pair of pants $\PP_3$} is the compact oriented surface of genus 0 with 3 boundary components $b_1, b_2, b_3$, equipped with  a $Y$-graph as in Figure  \ref{fig:P123}. By removing $b_3$ from $\PP_3$ we get $\PP_2$, which is a surface having one puncture (in place of $b_3$).  By removing $b_2$ from $\PP_2$ we get $\PP_1$,  which is a surface having  two punctures.

\FIGc{P123}{From left to right: $\PP_3, \PP_2, \PP_1$. The $Y$ graph is in red. }{1.8cm}

(b) The group $\BZ/3$ generated by the counterclockwise $2\pi/3$ rotation  acts on $\PP_3$, cyclically permuting the indices $1\to 2 \to 3\to 1$ in  $b_1, b_2, b_3$. Let $\ell_i$ and $a_{ij}$, with $1\le i, j \le 3$ be the curves given in Figure \ref{fig:StandC} and their images under the action of $\BZ/3$, with the convention that $a_{ij} = a_{ji}$. The curves/arcs $\ell_i$  and $a_{ij}$  are called {\bf standard.}

\FIGc{StandC}{Standard curves $\ell_1, a_{23}, a_{11}$. }{1.8cm}

(b)  For $\PP_k$ with $k=1,2$  the standard curves are $a_{ij}, \ell_i$ with $1\le i,j \le k$. 
\edf

 
 \def\SC{{\fS_\cC}}
\def\cCt{{\cC^{(2)}}}
\def\tiF{{\tau \in \cF}}


A {\bf pants decomposition} of $\fS=\Sigma_{g,m}$ is a maximal collection $\cC$ of disjoint non-trivial, non-peripheral loops which are pairwise non-isotopic.
By removing a small open tubular neighborhood of each $c\in \cC$ from $\fS$ we get a surface $\SC$ whose connected components are called {\bf faces} of the pants decomposition. Let $\cF=\cF(\cC)$ be the set of faces. We can assume that there is a projection $\pr: \SC \onto \fS$ which identifies pairs of circle boundary components of $\SC$. For $c\in \cC$ let $c'$ and $c''$ be the components of $\pr^{-1}(c)$, and let $\cCt$ be set of all components of $\partial \SC$.

A triple $(a,b,c)\in \cC$ is {\bf triangular} if the curves are the images of the three boundary components of a face under the projection $\SC \onto \fS$. Note that  two of the three $a,b,c$ might be equal. 


A {\bf dual graph $\Gamma$ of $\cC$} is a trivalent graph embedded into $\fS$, transverse to each $c\in \cC$  such that its preimage in each face $\tau$ is a $Y$-graph, making $\tau$ isomorphic to one of the three standard DT pair of pants. Note that $\Gamma$ does not have a vertex on any half-edge ending at a puncture.

\bdf  \label{defLength}

 A DT-datum of $\fS$ consists of a pants decomposition $\cC$ and a dual graph $\Gamma$.

 The  {\bf bold vertices} of $(\Gamma, \cC)$ are
elements of  $\Gamma \cap (\bigcup_{c\in \cC} c)$, as well as 
their lifts in $\SC$.

For a simple diagram  $\al$  on $\fS$ its 
{\bf length coordinate} $\al$ at $c\in \cC$ is $n_\al(c) :=I(\al, c)$.

A simple diagram  $\al$ is {\bf good with respect to $(\cC, \Gamma)$} if $\al$ does not contain any bold vertices and it is taut with respect to $\cC$, meaning  $|\al \cap c| = n_\al(c)$ for all $c\in \cC$.
\edf
An easy Euler characteristic count shows that
$
|\cC|= r:= 3g-3+m .
$

The following is standard, and is the basis of all the definitions of DT-coordinates.

\blem \label{rSlides}
 Two good simple diagrams are isotopic in $\fS$ if and only if they are related by a sequence of t-slides and loop-slides as seen in Figure \ref{fig:slide}.

\elem

\begin{figure}[h]
    \centering
    \includegraphics[width=330pt]{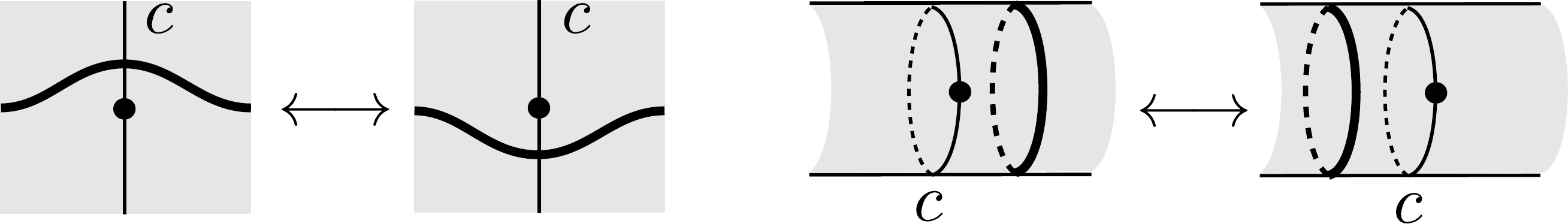}
    \caption{Left: $t$-slide, Right: loop-slide. Here $c\in \cC$.}
    \label{fig:slide}
\end{figure}

\def\hC{{ \hat \cC}}
\def\tiC{{ \tilde \cC}}

\def\tiC{{ \tilde \cC}}

\subsection{Punctured bordered surfaces associated to faces} To each face of the pants decomposition of $\fS$ we associate a pb-surface and analyze  how to patch curves on faces to get a global diagram on $\fS$.

\def\Pb{\check \tau}
\def\pPb{\partial{\Pb  }}
\def\PPb{{\check{\PP  }}}

For a face $\tau\in \cF$ let  $\Pb$ be the result of removing all the bold vertices from $\tau$. Such $\Pb$ is a punctured bordered surface.
If $\al\subset \fS$ is a simple diagram in good position and $\tau$ is a face, then $\pr^{-1}(\al) \cap \tau$, called the 
{\bf lift of $\al$ in a  
 $\tau$}, is a {\em strongly simple diagram on } $\Pb$ in the sense of Definition \ref{defSsimple}. Let $\Bsl(\Pb)$ be the  set of isotopy classes of strongly simple diagrams of $\Pb$.

Let $c$ be 
a  boundary component  of $\tau$ and let $\ell_c$ be a loop in the interior of $\tau$ parallel  to $c$. Define the permutation $\theta_c: \Bsl(\Pb) \to \Bsl(\Pb)$, called 
the {\bf twist on $c$}, as follow. If $\al \cap c=\emptyset$ then $\theta_c(\al)= \al$; otherwise $\theta_c(\al)$ is given by
$$ \incl{1.3cm} {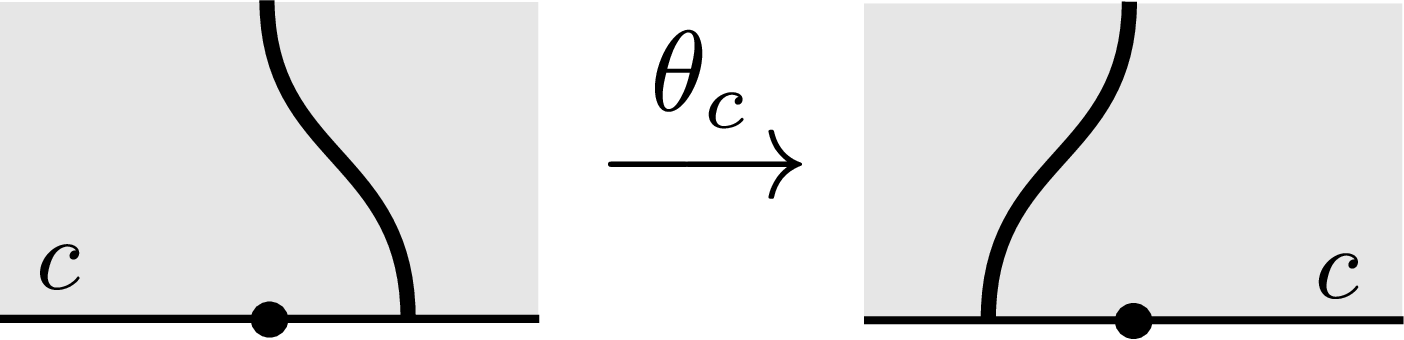}  $$

\def\patch{{\mathsf{patch}}}
Suppose $\al= (\al_\tau)_{\tau\in \cF}\in \prod_{\tau \in \cF} \Bsl(\Pb)$, where $\cF$ is the set of all faces. Let $\al^\sqcup=\sqcup _\tau \al_\tau$.
We say that $\al$ is  {\bf matched} if $|\al^\sqcup \cap c'| = |\al^\sqcup \cap c''|$ for each $c\in \cC$.
For matched $\al$  we can patch the $\al_\tau$ together to get  $\patch(\al) \in \Bsl(\fS)$. Let $(\prod_{\tau \in \cF} \Bsl(\Pb))^*$ be the set of all matched elements. From Lemma \ref{rSlides} we get 

\blem \label{rSlides2}
Two elements $\al, \al'\in (\prod_{\tau \in \cF} \Bsl(\Pb))^*$ are patched to the same element of $\Bsl(\fS)$ if and only if they are related by a sequence of 
\begin{itemize}\item  t-moves:\ 
$
\al \leftrightarrow \theta_{c'}^{-1} \circ \theta_{c''}(\al) ,\ c\in \cC.
$
\item loop-moves: Assume $\al$ contains $\ell_{c'}$, then the loop-move at $c$ is
$ \al \leftrightarrow (\al \setminus \ell_{c'} ) \cup \ell_{c''}.$

\end{itemize}
\elem
\def\Add{{\mathsf{Add}}}
 \subsection{Difference between our twist coordinates and the usual ones}
Before giving the detailed definition, let us point out the difference between our twists coordinates and the ones used  in \cite{Luo,FKL}.
 Consider the case of $\PP_3$, with boundary components $b_1, b_2, b_3$.  Denote by $\theta_i$ the twist map on $b_i$.

 For $i=1,2,3$, the arc $a_{ii}$ in $\PP_3$ will be called a {\bf standard return arc}; it has two end points on $b_i$, and we say that it 
 {\bf approaches $b_{i+1}$} if as in picture  it goes toward $b_{i+1}$ before coming back to $b_i$. As usual indices are taken modulo 3. Any arc $\theta_i^k(a_{ii})$, for $k\in \BZ$, is called a {\bf return arc approaching $b_{i+1}$}.

In the usual definition of \cite{Luo,FKL}, the twist coordinate $t'_i(\al)$, where $\al$ is a simple diagram,  is defined so that the contribution of each standard return arc  is 0. Our modified version of the twist coordinate is the following
\be 
t_i(\al)=    t'_i(\al) + \text{number of  return arcs approaching}\ b_{i}.
\ee
This modification is very important as it will allow to identify the highest degree term  in certain filtrations of the skein algebra. 

If $n_i= |\al \cap b_i|$, then the number of return arcs approaching $b_i$ is
\be 
\Add_i:= \max(0,\frac 12 ( n_{i-1}- n_i - n_{i+1})), \quad \text{indices taken mod 3}.
\ee

The modifications for $\PP_2, \PP_1$ are different.

\subsection{DT-coordinates for $\PP_3$}
\def\rot{{\mathrm{rot}   }}
\def\loo{\ell}
\def\arcc{\omega}

We now present DT-coordinates for  $\PPb_3$.

\bpro 
\label{rP3}
(a) There is a unique injective  map
$$ \nu: \Bsl( \PPb_3) \to \BN^{3} \times \BZ^{3}, \quad \nu(\al)= (n_1(\al),n_2(\al), n_3(\al), t_1(\al), t_2(\al), t_3(\al)), $$
satisfying all of the following (1)-(6).
\begin{enumerate}

\item length coordinates: 
$ 
n_i(\al) = |\al \cap b_i|.
$
\item disjoint additivity: If $\al_1,\dots, \al_k$ are components of $\al$ then $\nu(\al) = \sum_{j=1}^k \nu(\al_j)$.
\item The twist increases the twist coordinate by 1: The twist $\theta_i$ changes only the coordinate $t_i$, and the change is given by
\be 
t_{ i }(\theta_i \al) =  \begin{cases} t_i( \al) +1,  & n_i(\al) >0 \\
t_i( \al),  & n_i(\al) =0.\end{cases}
\ee
\item Twist coordinates of loops: \
$
t_i(\ell_j)=  \delta_{ij}. 
$

\item Standard straight arcs: For all applicable $i,j,k$ with $j\neq  k$,
\ $
t_i(a_{jk}) = 0 . 
$

\item Standard return arcs: With indices taken modulo 3,
\begin{align}
t_i(a_{jj}) &=   \delta_{i-1,j}. \label{eqDT3}
\end{align}
\end{enumerate}

(b) The image  $\Lambda_3 :=\nu(\Bsl(\PPb_3))$ inside $\mathbb{N}^3\times \mathbb{Z}^3$ is  the subset  constrained by 

\begin{itemize}
\item[(i) ] $n_1+n_2+n_3$  is even, and

  \item[(ii)] if $n_i=0$ then $t_i\geq \Add_i(n_1, n_2, n_3):=\frac 12 \max(0,n_{i-1} -n_i- n_{i+1} )$ for  $i\in \{1,2,3\}$.

\end{itemize}

Moreover, the set $\Lambda_3$ is a submonoid of $\BZ^6$.
\epro

\def\Lo{{\mathring \Lambda}}

This is \cite[Proposition 7.4]{BKL}, where our $\PPb_3$ and $\Bsl(\PPb_3)$ are respectively $\PPb_3$ and $\uB(\PPb_3)$ of \cite{BKL}. Moreover, we rescale so that our $t_i$ is one half of the $t_i$ in \cite{BKL}. Otherwise our Proposition \ref{rP3} is identical to \cite[Proposition 7.4]{BKL}.

\def\cb{{\check b}}

\subsection{DT-coordinates for $\PP_2$} \label{secP2}

We now introduce   DT-coordinates for  $\PPb_2$. Note that our $\PP_2$ is not the same as the $\PP_2$ of \cite{BKL} and our $\Bsl(\PPb_2)$ is not the same as $\uB(\PPb_2)$ of \cite{BKL}. Instead we modify the restriction of $\nu$ from $\PPb_3$ to $\PPb_2$. 

The boundary components of $\PPb_2$  are $\cb_1, \cb_2$, where $\cb_i$ is $b_i$ without the bold vertex.


\bpro\label{rP2}
 (a) There is a unique injective  map
$$ \nu: \Bsl( \PPb_2) \to \BN^{2} \times \BZ^{2}, \quad \nu(\al)= (n_1(\al),n_2(\al),  t_1(\al), t_2(\al)), $$
satisfying the conditions (1)-(5) of Proposition \ref{rP3}, and in addition, the twist values  of the standard return arcs are
\begin{align}
{ (t_1(a_{11}), t_2(a_{11})))=(0,1), \quad (t_1(a_{22})), t_2(a_{22})))=(-1,1). }
\label{eqP2}
\end{align}

(b) The image  $\Lambda_2:=\nu(\Bsl(\PP_2))$ is the submonoid of $\mathbb{N}^2\times \mathbb{Z}^2$ constrained by:  
\begin{itemize}
\item [(i)]  $n_1 + n_2 \in 2 \BN$ and,

\item [(ii)] if $n_i=0$ then 
$t_i\ge \Add_i(n_1, n_2)$, where 
$\Add_1(n_1, n_2) = -\frac{n_2}2$ and  
$\Add_2(n_1, n_2) = \frac{n_1}2$.
\end{itemize}

\epro
\bpr  The identity embedding $f: \PPb_2 \embed \PPb_3$  injects $\Bsl(\PPb_2)$  into $\Bsl(\PPb_3)$, with image 
\be 
f(\Bsl(\PPb_2))=\{ \al \in \Bsl(\PPb_3 )\mid n_3(\al) = 0, \al \ \text{has no $\ell_3$ components}\}.
\label{eq16}
\ee

For $\al \in \Bsl( \PPb_2)$ define $ n_1(\al) = n_ 1(f(\al)),\ n_2(\al) = n_ 2(f(\al)))$, and 
\begin{align}
  t_1(\al) & = t_ 1(f(\al)) - \frac 12 \max(0, n_2-n_1)  \label{eqp22}\\ 
   t_2(\al) &= t_ 2(f(\al)) +  \frac 12 \max(0, n_2-n_1). \label{eqp2222}
\end{align}
This means that we have \eqref{eqP2} and properties (1) and (3)-(5) of Proposition \ref{rP3}. Let us prove the remaining property (2). Note that $r_2(\al)= (n_2(\al)-n_1(\al))/2$ is the number of return arcs based at $b_2$. As such, the number $r_2$ is additive, meaning $r_2(\al) = \sum_{j=1}^k r_2(\al_i)$. From here and the definitions \eqref{eqp22} and \eqref{eqp2222} we get the additivity of $t_1$ and $t_2$. The additivity of $n_1$ and $n_2$ is clear. Thus we have property~(2).

As $t_1(f(\al)), t_2(f(\al))$, and thus all the DT-coordinates of $f(\al)$, are determined by $\nu(\al)$, we see that $\nu$ is injective.

(b) Let us prove the image of $\nu$ is $\Lambda_2$. This can be done easily by examining the image of $f(\Bsl(\PP_2))$ in \eqref{eq16} and the description of $\Lambda_3$ given in Proposition \ref{rP3}. Here are the details.

Let $\al\in \Bsl(\PPb_2)$. Since $n_3(f(\al))=0$ we have $n_1(\al) + n_2(\al)\in 2 \BZ$, which is condition (i). Let us prove (ii). First assume $n_1(\al)=0$.
By \eqref{eqp22} and Proposition \ref{rP3}(b), 
\begin{align*}
t_1(\al) \ge \frac 12 \max(0, n_3- n_1-n_2)  - \frac 12 \max(0, n_2-n_1)= -\frac{n_2}2,
\end{align*}
where the last identity follows since $n_1=n_3=0$. Similarly one can easily show that if $n_2=0$ then $t_2 \ge n_1/2$. This proves $\nu(\Bsl(\PPb_2)) \subset \Lambda_2$.

By reversing the above argument one can also easily show that $\Lambda_2 \subset 
\nu(\Bsl(\PPb_2))$. This proves $\nu(\Bsl(\PPb_2)) = \Lambda_2$.

\no{
 Then $\al$ consists of $n_2/2$ parallel copies of a twist $ \theta^k(a_{22})$ for some $k\in \BZ$, and some, say $l$ copies of $\ell_1$. From the construction $t_1(\al) = l-n_2/2\ge n_2/2 = \Add_1(n_1, n_2)$. Assume now $n_2(\al)=0$. Then $\al$ consists of $n_1/2$ parallel copies of a twist $ \theta^k(a_{11})$, and some, say $l$ copies of $\ell_2$. From the construction $t_1(\al) = l+n_1/2\ge n_2/2 = \Add_2(n_1, n_2)$.

Let us show $\nu: \Bsl(\PPb_2\to \Lambda_2$ is surjective. 
Reversing the above argument, we see that if $n_i=0$ and $t_i\ge \Add_i(n_1, n_2)$, then there is $\al \in \Bsl(\PPb_2)$ such that $\nu(\al) = (n_1, n_2, t_1, t_2)$.
Assume $n_1, n_2 >0$, with $n+1+ n_2\in 2 \BZ$, and $t_1, t_2\in \BZ$ arbitrary.
By letting $n_3=t_3=0$.

The constraints of $\Lambda_3$ given in Proposition \ref{rP3} gives us the image $\Lambda_2$, constrained by (i) and (ii).}
 Let us prove that 
$\Lambda_2$ is a monoid.  For any $(\bn', \bt'), (\bn'', \bt'') \in \Lambda_2$ we need to show $(\bn, \bt)= (\bn'+ \bn'', \bt'+ \bt'')\in \Lambda_2$.  Observe that
 $ \{ (\bn, \bt)  \in (\BN_{>0})^2\times \BZ^2 \mid n_1 + n_2  \in 2 \BN\} $ is  subset of $\Lambda_2$.
 
Hence if $\bn \in (\BN_{>0})^2$ then $(\bn, \bt) \in \Lambda_2$.

Assume $\bn = (n_1, n_2) \not \in (\BN_{>0})^2$. Then $n_1=0$ or $n_2=0$. Say $n_1=0$. Then using the linearity of $\Add_1$, we have $t_1= t'_1+ t''_1 \ge -n'_2/2+ -n''_2/2= - n_2/2$. 
Similarly if $n_2=0$ then the linearity of $\Add_2$ implies $t_2 \ge \Add_2(n_1, n_2)$. This completes the proof. \epr

\subsection{DT-coordinates for $\PP_1$}

We now introduce  the  DT-coordinates for $\Bsl(\PPb_1)$.  

\bpro\label{rP1}
 (a) There is a unique injective  map
$$ \nu: \Bsl( \PPb_1) \to \BN \times \BZ, \quad \nu(\al)= (n_1(\al),  t_1(\al)), $$
satisfying the conditions (1)-(4) of Proposition \ref{rP3}, and besides,
$
 t_1(a_{11}) = 1.
$

(b) The image  $\Lambda_1:=\nu(\Bsl(\PP_1))\subset 2\mathbb{N}\times \mathbb{Z}$ is a submonoid defined by: If $n=0$ then $t\ge 0$.

\epro
\bpr Let $\al \in \Bsl(\PPb_1)$. If $|\al \cap b_1|= 0$, then $\al$ is the disjoint union of $l$ copies of $\ell_1$, and define $(n_1(\al), t_1(\al))=(0,l)$.

If $|\al \cap b_1|=n\neq 0$ then $n$ must be even. There is a unique $k\in \BZ$ such that  $\theta_1^k (\al) $ is the disjoint union of $n/2$ parallel copies of $a_{11}$. Define $(n_1(\al), t_1(\al))=(n,-k+n/2)$. It it easy to check all the claims of the lemma. 
\epr

\subsection{Dehn-Thurston Coordinates} We  define the DT-coordinates for strongly simple diagrams on $\fS= \Sigma_{g,m}$, equipped with a DT-datum $(\cC, \Gamma)$, by combining the faces.

For each face $\tau$ we choose an identification of $\tau$ with one of the three DT pairs of pants $\PP_1, \PP_2, \PP_3$, and call it the {\bf characteristic map} of $\tau$. For a type $\PP_3$  face there are three possibilities of the characteristic maps (up to isotopies), but the definition given below does not depend on the choice, as all relevant functions are $\BZ/3$-equivariant.

Let $\al$ be a strongly simple diagram in good position with respect to $(\cC, \Gamma)$. The length function $n_\al: \cC \to \BN$ was given in Definition \ref{defLength}. Define
$ t_\al: \cC \to \BZ$
as follows. For a face $\tau$ let $\al_\tau= \pr^{-1}(\al) \cap \tau$. For $c\in \cC$ recall that $c', c''\in \cCt$ are the lifts of $c$. Define
\begin{align}
 t_\al(c)&= t(\al';c') + t(\al'';c'').
\end{align}
Here 
\begin{itemize}
\item $\tau'$ (respectively $\tau'')$ is the unique face containing $c'$ (respectively $c''$),
\item $t(\al;c')= t_i(\al_{\tau'})$, where the characteristic map sends $c'$ to the $i$-th boundary component $b_i$.
Define  $t(\al''; c'')$ similarly.
\end{itemize}

\def\fA{{\mathfrak A}}
\def\fB{{\mathfrak B}}

\bpro \label{rDTcoord}
Let $(\cC, \Gamma)$ be a DT-datum of $\fS= \Sigma_{g,m}$ where $(g,m) \neq (1,0),(0,k)$ for $k \le 3$. The following map $\nu$ is well-defined and injective
$$ \nu: B(\fS) \to \BN^\cC \times \BZ^\cC,\ \nu(\al) = (n_\al, t_\al).$$

\epro

\begin{proof}   To show the well-definedness 
we need show that if $\al$ and $\al'$ are strongly simple diagrams in good position, and they are isotopic then $\nu(\al) = \nu(\al')$. Since $n_\al(c) = I(\al, c)$, clearly $n_\al(c)= n_{\al'}(c)$.
By Lemma \ref{rSlides} the two diagrams $\al$ and $\al'$ are related by a sequence of t-slides and loop-slides, which preserve  $t_\al$ by Properties (3) and (4) of Proposition \ref{rP3}. This shows $\nu$ is well-defined.

Now we prove the injectivity. Assume $\nu(\al)= \nu(\beta)$. Consider the collections $\fA=(\al_\tau)_{\tau\in \cF}$ and $\fB=(\beta_\tau)_{\tau\in \cF}$. Let $c\in \cC$. Since
$$t(\al;c') + t(\al; c'') = t_\al(c)= t_\beta(c)= t(\beta;c') + t(\beta; c'')$$ we have 
$$ t(\al; c') - t(\beta;c')= t(\beta;c'')-  t(\al; c''):=k.$$
 Then after $k$ twists on $c$, we can bring $\beta$ to $\beta'$ with $t(\beta';c')= t(\al;c')$ and $t(\beta';c'')= t(\al;c'')$. Repeating this procedure to all $c\in \cC$,  we can isotope $\beta$ to $\gamma$ with $t(\gamma;c')= t(\al;c')$ and $t(\gamma;c'')= t(\al;c'')$ for all $c\in \cC$. This shows $\gamma_\tau$ is isotopic to $\al_\tau$ for all faces $\tau$. Hence $\al$ is isotopic to $\gamma$, and whence to $\beta$.
\end{proof}

To describe the image of $\nu$ let us introduce the following notion.
For $\bn: \cC\to \BN$ and  $b\in \cCt$ we define $\Add(b; \bn)$ as follows. If the face containing $b$ is  of type $\PP_1$ then let $\Add(b; \bn)=0$. If $\tau$ is of type $\PP_2$ or $\PP_3$, then choose a characteristic map of $\tau$ such that $b= b_i$, then define $\Add(b; \bn)=\Add_i(\bn)$. From the images of $\nu$ given in Propositions \ref{rP3}, \ref{rP2}, and \ref{rP1} we get

\bpro \label{rDTcoord2}
The image $\Lambda_{\cC,\Gamma}= \nu( B(\fS)  )$  is the submonoid of $ \BN^\cC \times \BZ^\cC$ consisting of $(\bn, \bt)$ satisfying
\begin{itemize}
\item if $(a,b,c)\in \cC$ is a triangular triple then $\bn(a)+ \bn(b) + \bn(c) \in 2 \BZ$,

\item if $\bn(c) =0$ then $\bt(c) \ge \Add(c'; \bn) + \Add(c''; \bn)$.

\end{itemize}
\epro

\def\uSS{{\underline{\cS}(\fS)}}
\def\cI{{\mathcal I}}
\def\uS{{\underline{\cS}}}
\def\Left{{\mathsf{Left}}}

\def\ctau{{\check \tau}}
\def\Ytau{{\cY(\ctau)}}

\section{A quantum trace for the DT pairs of pants}
\label{secQT}

For each DT pair of pants $\tau=\PPb_j$, we show that there is an embedding of the sliced skein algebra $\cS^\sli(\ctau)$ into a quantum torus $\Ytau$. This map is a quotient of the so-called $A$-version of the quantum trace developed in \cite{LY2, BKL}, composed with a linear map.
 These maps will be later glued together to give the global map on $\Sigma_{g,m}$.
 
 In this section $j\in \{1,2,3\}$ and it is the index of $\PP_j$. Recall that the boundary components of $\PP_3$ are $b_1, b_2, b_3$, and  $\PPb_j= \PP_j \setminus \{\text{ bold vertices}\}$. Each $\PPb_j$ is a punctured bordered surface.
 
\subsection{Quantum tori associated to $\PP_j$}\label{ssQtPj} We now  define a quantum torus $\cY(\ctau)$ for  $\ctau=\PPb_j$.

As $\PP_2= \PP_3\setminus b_3$, we also use 
$b_3$ to denote the only interior puncture of 
$\PPb_2$. Similarly the interior punctures of $\PPb_1$ 
are $b_3$ and $b_2$.

Recall that in Section \ref{secMRY}  the universal  sliced skein algebra  of $\fS$ is defined over the ring $\RS$.  For brevity, for $\tau= \PP_j$, we denote the ground ring $\cR_{\ctau}$ by $\RR_j$. Thus, by definition,
$$\RR_3 = \Zq, \ \RR_2= \Zq[ b_3^{\pm 1}], \ \RR_1= \Zq[b_3^{\pm1}, b_2^{\pm1}].$$

Define  $\cY(\PPb_3)$ by the following, where  indices are taken mod 3:
\be 
\cY(\PPb_3)= \RR_3\la x^{\pm 1}_i, u_i^{\pm 1}, i=1,2,3 \ra/ (x_{i+1} x_i = q x_i x_{i+1}, u_i u_j = u_j u_i, u_i x_j = q^{2\delta_{ij}} x_ju_i   ).
\ee

Define $\cY(\PPb_2)$  by dropping $x_3$ and $ u_3$, but changing $\RR_3$ to $\RR_2$:
\begin{align}
\cY(\PPb_2)& :=  \RR_2\la x_1^{\pm 1},x_2^{\pm 1},  u_1^{\pm1},  u_2^{\pm1}\ra /( u_i x_j = q^{2\delta_{ij}} x_j u_i, x_2 x_1 =q x_1 x_2, u_i u_j = u_j u_i)
\end{align}
Finally, further dropping $x_2$ and $u_2$ and changing $\RR_2$ to $\RR_1$,
\begin{align}
\cY(\PPb_1) &:=  \RR_1\la x_1^{\pm 1}, u_1^{\pm1} \ra /( u_1 x_1= q^2x_1 u_1).
\end{align}
Each $\cY(\PP_j)$ is a $\Zq$-algebra with reflection $\omega$ which fixes each of $x_i, u_i, b_i$.

Each quantum torus $\cY(\PP_j)$  has a free $\RR_j$-basis consisting of monomials
\be B(\cY(\PPb_j)):=\{Y^\bnu \mid \bnu \in \BZ^{2j}\}, \ Y^\bnu=[x^\bn u^\bt]_\Weyl \ \text{for } \ \bnu=(\bn, \bt)\in \BZ^j \times \BZ^j.
\label{eqBaseY}
\ee
For $i\le j$ let
$$ \cY(\PPb_j)_{\deg_i =k} = \RR_j\text{-span of monomials having  degree of $x_i =k$}.$$
 
\def\Ssl{{\cS^\sli}}
\subsection{Filtrations}  

\label{ssFiltr6} We now define  function $\ddd_\tau$, for $\tau= \PP_j$, which helps to define a preorder (a linear order without the antisymmetric property) on $\Bsl(\PPb_j)$ and related filtrations.

Our DT-coordinates parameterize the $\RR_j$-basis $\Bsl(\PPb_j)$ of the sliced skein algebra $\Ssl(\PPb_j)$:
$$ \nu: \Bsl(\PPb_j) \xra{\cong} \Lambda_j \subset \BN^j \times \BZ^j.$$

 To compare elements of $\Lambda_j\subset \BN^j \times \BZ^j$, define
 $\ddd_{\PP_j}: \BN^j \times \BZ^j\to  \BZ^3$ by
\begin{align}
\ddd_{\PP_3} ((n_1, n_2, n_3, t_1, t_2, t_3))&= (n_1+n_2+n_3, t_1+t_2+t_3,0) \\
\ddd_{\PP_2} ((n_1, n_2, t_1, t_2))&= (n_1+n_2, t_1+t_2, t_2) \\
\ddd_{\PP_1}(n_1, t_1) &= (n_1, t_1,0) 
\end{align}
The reason  why $\ddd$ is defined as above will be clear in the proof of Theorem \ref{thmbtr}. Note that only $\ddd_{\PP_2}$ has non-zero third coordinates in its image. This is dictated by the complexity of the quantum trace for $\PP_2$  that will be constructed.

Using the lexicographic order on  $\BZ^3$, for $\bk\in \BZ^3$ define
\be 
F_\bk (\cY(\PPb_j)) = \cR\text{-span of } \{ Y^\bbl \mid \ddd_{\PP_j}(\bbl) \le \bk\}.
\ee

For two monomials $Y^\bnu,Y^{\bnu'}\in B(\cY(\PPb_j))$, we say $Y^\bnu$ has {\bf higher $\ddd_{\PP_j}$-degree} if $\ddd_{\PP_j}(\bnu) > \ddd_{\PP_j}(\bnu')$. Note that this is not an linear order, but a preorder on $B(\cY(\PPb_j))$.

\subsection{Quantum traces for $\cS^{\mathsf{sl}} (\PP_j)$ } \label{ssUtr}
We now formulate main results of this section.

Recall that for $i\le j$, the loop $\ell_i$ is in the interior of $\PPb_j$ parallel to $b_i$, and  $\theta_i$  is the slide (or twist) on $\cb_i$. Also if $\al\in \Bsl(\fS)$ then $[\al]_\omega\in \SslS$ is the element $\al$ normalized by \eqref{eqRefl} so that it is reflection invariant.

\def\YPj{{ \cY(\PPb_j )  }}
\def\SPj{\Ssl(\PPb_j)}
\def\Stau{\Ssl(\ctau) }
\def\Ytau{\cY(\ctau) }

\bthm \label{thmbtr} For $\tau= \PP_j$, where $j=1,2,3$,
 there is a reflection invariant $\RR_j$-algebra homomorphism
$$ \utr_\tau: \Stau\to \Ytau$$
having the following properties: For $i\le j$ and $\al \in \Bsl(\ctau)$ with $|\al \cap b_i| = k$,

\begin{enumerate}

\item Boundary grading:
\be \utr_\tau  
\left( 
[\al]_\omega 
\right) 
\ \in \ \Ytau_{\deg_i=k} .  \label{eqGrade_utr}
\ee

\item Near boundary loop: The value of the loop $\ell_i$ is given by
\be \utr_\tau (\ell_i) = u_i + u_i^{-1}. \label{eqLoop_utr}
\ee

\item Twist:  if $|\al \cap b_i|=k\neq 0$ then
\be 
\utr_\tau ([\theta_i(\al)]_\omega)=  q^{-k} \, u_i\, \utr(\al) = [  u_i\,\, \utr(\al) ]_\Weyl.
\label{eqTwist_utr}
\ee

\item Highest order term: 
\be 
\utr_\tau([\al]_\omega) \qeq Y^{\nu(\al)} + F_{<\, \ddd_{\tau}(\nu(\al)) }(\Ytau). \label{eqHighdeg_utr}
\ee
\end{enumerate}

\ethm
Parts (1)--(3) will be used to glue the maps $\utr_\tau$ along boundary components to get a global map. Then part (4) will be used to determine the top degree part of the global map.

For $\tau= \PP_3$, Theorem \ref{thmbtr} is the $\PP_3$ case of Theorem \cite[Theorem 8.1]{BKL}, where  our $\PP_3, \Bsl(\PPb_3), u_i, \cY(\PPb_3)$ are respectively $\PP_3, \uB(\PP_3), u_i^2, \cY(\PPb_3)$ of \cite{BKL}. 
We will prove the theorem for the remaining cases $\tau= \PP_2, \PP_1$ in  subsequent subsections.

\def\Rtau{{\cR_{\ctau}}}

\brem All the conclusions of Theorem \ref{thmbtr} are transferable from the ground ring $\Rtau=\RR_j$ to any commutative $\Rtau$-domain.
\erem

\def\vcup{{ \overrightarrow{\cup}}}

\subsection{Reduction to 1-component case}

\blem \label{rReduct}
Assume there is a reflection invariant $\RR_j$-algebra homomorphism 
$$ \utr_\tau: \Stau \to \Ytau$$
satisfying all the conditions (1)-(4) of Theorem \ref{thmbtr} for all 1-component $\al \in \Bsl(\ctau)$. Then we also have (1)-(4) for all $\al \in \Bsl(\ctau)$, i.e. we have the Theorem \ref{thmbtr}.
\elem

\bpr The proof is identical to that of \cite[Lemma 8.3]{BKL}. {Here are the details.

Since \eqref{eqLoop_utr} concerns only 1-component element, it holds true by assumption.

Let $\al\in \Bsl(\ctau)$ have  connected components $\al_1, \dots, \al_k$ with $k \ge 2$. By height exchange,
\be \al \qeq \al_1 \dots \al_k \label{eq1comp}
\ee 
The multiplicative nature of $\utr_\tau$ shows that \eqref{eqGrade_utr}  and \eqref{eqHighdeg_utr}
also hold for $\al$.

It remains to prove \eqref{eqTwist_utr}. Since $k=|\al \cap b_i| >0$ one of components of $\al$, say $\al_1$ after  re-indexing, intersects $b_i$.
Recall that $u_i$ commutes with all  variables except for $x_i$, for which $u_i x_i= q^2 x_i  u_i$. Hence from \eqref{eqGrade_utr} we get $u_i   \utr_\tau([\al]_\omega) = q^{2k} \utr_\tau([\al]_\omega)\, u_i$. By definition 
$$[u_i^2\,\utr_\tau ([\al]_\omega) ]_\Weyl = q^{-k} u_i^2\, \utr_\tau([\al]_\omega).$$
Since $\utr_\tau([\al]_\omega)$ is reflection invariant, $[u_i^2\,\utr([\al]_\omega) ]_\Weyl$ is reflection invariant by Lemma \ref{rReflection}.

From the definition of $\theta_i$ and the height exchange rule, we get 
$$ \theta_i(\al_1\al_2 \dots \al_k)\qeq \theta_i(\al_1) \al_2 \dots \al_k.$$
Using \eqref{eq1comp} and  the above identity, we get
$$ \utr([\al]_\omega) \qeq \utr(\theta_i(\al_1)) \utr( \al_2 \dots \al_k) \qeq u_i^2\, \utr([\al]_\omega).$$
Using Identity \eqref{eqTwist_utr} for $\al_1$, we get
$$ \utr([\al]_\omega)  \qeq u_i^2\, \utr([\al]_\omega).  $$
By Lemma \ref{rReflection}(b), in the above  we can replace $\qeq$ by $=$,  and get \eqref{eqTwist_utr}.
}
\epr

\def\aaa{{\mathbf a}}


As said, when  $\tau=\PP_3$ Theorem  \ref{thmbtr}  was  part of \cite[Theorem 8.1]{BKL}. We record here the exact values of 1-component simple diagrams.

\blem[ Lemma 8.5 of \cite{BKL}]
 \label{rIndent1}
 For $i\neq k$ in $\{1,2,3\}$ and $m\in \BZ$  we have
\begin{align} 
\utr(\ell_i) & = u_i + u_i^{-1}\\
\utr(\theta_i^m( a_{ik}))&= [u_i^{m}x_i x_k]_\Weyl \label{eqa23} \\
 \utr( [\theta_i^m(a_{ii})]_\omega)) &=  [u_i^{m} u_{i+1} x_i^2]_\Weyl+[u_{i}^{m+1} u_{i+2}^{-1} x_i^2]_\Weyl. \label{eqgamma1}
\end{align}
\elem

\def\vBsl{\vec{\mathsf B}^\sli}
\subsection{Proof of Theorem \ref{thmbtr} for $\PP_2$} We deduce the case $\tau=\PP_2$  from $\tau=\PP_3$.

\begin{proof} By Proposition \ref{rslBasis},  $\vBsl(\ctau)$ is a free $\RR_j$-basis of $\Ssl(\ctau)$. 
The embedding $\PPb_2 \embed \PPb_3$ gives
 $$\vBsl(\PPb_2)= \{ \al \in \vBsl(\PPb_3)\mid n_3(\al)=0, \ell_3 \not \subset \al \} \subset \vBsl(\PPb_3).$$

Let $S''\subset \Ssl(\PPb_3)$ be the $\RR_3$-span of $\vBsl(\PPb_2)$ and $\cY''= \cY(\PPb_3)_{\deg(x_3)=0} $.  By comparing the basis over $\RR_2$ we have $ S'' \ot_{\RR_3} \RR_2= \Ssl(\PPb_2)$. 
From \eqref{eqGrade_utr} for $\PP_3$ we have
\be \utr_{\PP_3}( S'') \subset \cY''.
\ee 
Hence $\utr_{\PP_3}$ gives rise to the first  of the following algebra homomorphisms
\be  \Ssl(\PP_2) \to \cY''\ot_{\RR_3}  \RR_2 \onto  (\cY''\ot_{\RR_3}\RR_2)/(u_3= b_3) \equiv \cY(\PPb_2).
\label{equtr2}
\ee
Let $\utr_{\PP_2}$ be the composition of the sequence in \eqref{equtr2}. From Lemma \ref{rIndent1} we have

\blem \label{rIndent2}   Let $\tau = \PP_2$. For $i=1,2$ and $m\in \BZ$  we have
\begin{align} 
\utr_\tau(\ell_i) & = u_i + u_i^{-1} \label{eqell12} \\
\utr_\tau\theta_i^m( a_{12}))&= [u_i^{m}x_1 x_2]_\Weyl  \label{eqa12}
 \\ 
\utr_\tau [\theta_1^m(a_{11})]_\omega)) &=  [u_1^{m} u_{2} x_1^2]_\Weyl+ b_{3}^{-1}[u_{1}^{m+1}  x_1^2]_\Weyl. \label{eqa11}\\
\utr_\tau( [\theta_2^m(a_{22})]_\omega)) &=  [u_2^{m+1} u_1^{-1} x_2^2]_\Weyl+ b_3 [ u_2^m x_2^2]_\Weyl . \label{eqa22}
\end{align}
\elem

Let us now prove the Theorem. As each map in  \eqref{equtr2} is reflection invariant, so is $\utr_\tau$.
Suppose  $\al\in \Bsl(\PPb_2)$ has one component. Then $\al$ is one of  the curves on the left hand sides of \eqref{eqell12}-\eqref{eqa22}. On the right hand side of each of \eqref{eqell12}-\eqref{eqa22}, the first monomial has $\ddd_{\PP_2}$-degree higher than any other monomial. The definition of $\ddd_{\PP_2}$ was {\em designed} for this property. The exponent of the first monomial is  equal to the DT coordinate of the curve of the left hand side. Our twist coordinates are {\em defined} for this to  be true.

Thus the identities of
Lemma \ref{rIndent2} show that all conditions (1)-(4) of Theorem~\ref{thmbtr} are satisfied for $\al$. By Lemma \ref{rReduct} we have the theorem for $\PP_2$. .
\end{proof}

\def\IM{{\mathrm{Image}}}

\subsection{Proof of Theorem \ref{thmbtr} for $\tau=\PP_1$} The proof is similar to the case   of $\PP_2$.

\begin{proof} Using $\PPb_1\subset \PPb_2$ we consider $\vBsl(\PPb_1)$ as a subset of $\vBsl(\PPb_2)$. 
Let $S'\subset \Ssl(\PPb_2)$ be the $\RR_2$-span $\vBsl(\PP_1)$, and $\cY'= \cY(\PPb_2)_{\deg(x_2)=0} $. From \eqref{eqGrade_utr} for $\PP_2$ we have
$ \utr_{\PP_2}( S') \subset \cY'.
$
This shows that $\utr_{\PP_2}$ gives rise to the first 
of the following algebra homomorphisms
\be  \Ssl(\PP_1) \to (\cY'\ot_{\RR_2}\RR_1)   \onto (\cY'\ot_{\RR_2}\RR_1)/(u_2= b_2)  
\equiv \cY(\PPb_1).
\label{equtr1}
\ee
Let $f$ be the composition of the maps of \eqref{equtr1}. From Lemma \ref{rIndent2} we have
\begin{align} 
f(\ell_1) &= u_1+ u_1^{-1}   \label{eqa2323a}\\
f(  [\theta_1^m(a_{11})]_\omega)) &=  b_3^{-1} [u_1^{m+1} x_1^2]_\Weyl + b_2 [ u_1^m x_1^2]. \label{eqa1111a}
\end{align}
The curves on the left hand sides of \eqref{eqa2323a} and \eqref{eqa1111a} generate the $\RR_1$-algebra $\Ssl(\PPb_1)$, while the right hand sides have even order in $x_1$. It follows that $\IM(f) \subset \cY(\PPb_1)^\ev$, which is the subalgebra having even degree in $x_1$. Let $\utr_\tau$, with $\tau=\PP_1$, be the composition
$$ \Ssl(\PP_1) \xar{f}   \cY(\PPb_1)^\ev \xar{g} \cY(\PPb_1)^\ev \subset \cY(\PPb_1).$$
where $g$ is the algebra homomorphism given on generators $u_1$ and $x_1^2$ by $g(u_1)=u_1$ and $g(x_1^2) = b_3 x_1^2$. Identities \eqref{eqa2323a} and \eqref{eqa1111a} translate to
\begin{align} 
\utr_\tau(\ell_1) &= u_1+ u_1^{-1}   \label{eqa2323}\\
\utr_\tau(  [\theta_1^m(a_{11})]_\omega)) &=  [u_1^{m+1} x_1^2]_\Weyl + b_2 b_3 [ u_1^m x_1^2], \ m \in \BZ. \label{eqa1111}
\end{align}

We can now finish  proving the Theorem. As each map in \eqref{equtr1} and $g$ is reflection invariant, so is $\utr_\tau$.
Suppose  $\al\in \Bsl(\PPb_1)$ has one component. Then $\al$ is  one of the curves on the left hand sides of \eqref{eqa2323}-\eqref{eqa1111}.  On the right hand side of each of \eqref{eqa2323}-\eqref{eqa1111}, the first monomial has $\ddd_{\PP_1}$-degree higher than any other monomial.
The explicit values of $\utr_\tau(\al)$ given in 
\eqref{eqa2323}-\eqref{eqa1111} and the definitions of $\ddd_{\PP_1}$ and the DT-coordinates show that all conditions (1)-(4) of Theorem \ref{thmbtr} are satisfied for $\al$. Again the definitions of $\ddd_{\PP_1}$ and the twist DT-coordinates are defined for this purpose. By Lemma \ref{rReduct} we have the theorem in the remaining case.
\end{proof}

\def\xar{\xrightarrow}
\def\dd{\mathbb D  }
\def\tEEE{\tilde{\mathsf E}}
\def\EEE{{\mathsf E}}
\def\tPhi{{\tilde \Phi}}
\def\barm{{\bar m}}

\def\tmQ{{\tilde{\mQ}}}

\def\xxx{\red{xxx}}
\section{Degenerations of sliced skein algebras and quantum tori}
\label{secDegen}

In this section we prove that the sliced skein algebra $\SslR$, where $\fS= \Sigma_{g,m}$ and $\cR$ is a commutative domain, has a degeneration  which is a monomial subalgebra of a quantum torus. Consequently $\SslR$ is
 a domain. 

 We will also construct an algebra $\BN$-filtration of $\SslR$ such that the associated graded algebra embeds into the same quantum torus.


 In this section fix $\fS= \Sigma_g \setminus \{ v_1, \dots, v_m\}$, with  $(g,m) \neq (0,k), (1,0), (1,1)$ for $k \le 3$; or equivalently $r:= 3g-3+m >0$.
 We also fix a commutative ring $\cR$ having an invertible element $\xi$, and an $m$-tuple $\vw=(\w_1, \dots, \w_m)\in \cR^m$.

\def\YS{\cY(\fS)  }
 \def\LT{{\mathsf{LT}}}
 \def\led{\le_\ddd }
 \def\YxS{\cY_\xi(\fS;\cR)}
 \def\TxL{\bT_\xi(\tmQ,\cR;\Lambda)}
\subsection{A special DT-datum} \label{ssN1}

To make it technically easier we choose the following DT-datum $(\cC, \Gamma)$ for $\fS$ as described in Figure \ref{fig:cC}.

\FIGc{cC}{The DT-datum $(\cC,\Gamma)$, with $\Gamma$ in red. Left: the case $g\ge 1$. Right: $g=0$. The small white disks are the punctures. }{3cm}

When  $g=0$, both $\cC$ and $\Gamma$ are fully described. When $g\ge 1$,  first choose an arbitrary collection of loops forming a pants decomposition of $\Sigma_g$, then in  a face $\tau$ remove $m$ punctures  near a boundary component $b$, and add $m$ parallels of $b$ (on $\tau$, which separate the punctures) to the collection to get $\cC$. 

Let $\cF$ be the set of faces of the pants decomposition.
There is no  $\PP_1$ face if $g\ge 1$, and there are two $\PP_1$ faces if $g=0$. There are $\barm$  faces of type $\PP_2$ where $\barm = m$ if $g \ge 1$, and $\barm=m-4$ if $g=0$. We will number elements of $\cC=\{ c_1, \dots, c_r\}$ so that the first $\barm +1$ elements are as in  Figure \ref{fig:cC}.

Recall that $\SC$ is the disjoint union of the faces, and $\pr: \SC \onto \fS$ is the natural projection. Let $\cCt$ be the
components of $\partial \SC$. Each element $c\in \cC$ has two lifts $c'$ and $c''$ in $\cCt$. 

 What is important with the choice of $\cC$  is the following property:\\
{\bf (N1)} In a $\PP_2$ face with characteristic map $f$, the curve $\pr(f(b_2))$ comes before $\pr(f(b_1))$.

Condition {\bf (N1)} is used to transfer certain inequalities from faces to  $\fS$, see  Lemma~\ref{rOrder2}.

Recall that  $\Bsl(\fS)$ is  a basis of the sliced skein algebra $\SslR$.  We have  the DT-coordinates map $\nu: \Bsl(\fS) \xra{\cong} \Lambda= \Lambda_{\cC,\Gamma} \subset \BZ^{\cC} \times \BZ^\cC$. Let 
$$\mu= \nu^{-1}: \Lambda \to \Bsl(\fS).$$
 To define a linear order on $\Bsl(\fS)$, consider the following additive embedding $\ddd$:
 \begin{align*}
\ddd&: \BN^\cC \times \BZ^{ \cC} \embed \BZ^{2r}, \ \ddd(\bn, \bt) = ( \sum_{c\in \cC} \bn(c), \sum_{c\in \cC}  \bt(c), \bt(c_1), \dots, \bt(c_{r-1}),  \bn(c_1), \dots, \bn(c_{r-1})).
\end{align*}
We give  $\Lambda$ a structure of an ordered monoid, by equipping it with the {\bf $\ddd$-order}, defined by
$$ \bk \le_\ddd \bl \ \text{if} \ \ddd(\bk) \le \ddd(\bl) \ \text{in lexicographic order of } \ \BZ^{2r}.$$
The bijection  $\nu$  carries the $\ddd$-order to $\Bsl(\fS)$, meaning
 $\al \le_\ddd \beta$ if $ \nu(\al)   \le_\ddd \nu(\beta)$.

The linear order on $\Lambda$ allows to define a good  $\Lambda$-filtration on $\SslR$:
\begin{align}
E_\bk(\SslR) &:= \cR\text{-span of } \ \{ \mu(\bl) \mid \bl  \led \bk\}, \ \bk\in \Lambda. \label{eqFil2}
\end{align}
We will show that the filtration  $(E_\bk(\SslR)_{\bk\in \Lambda}$  respects the product. 

\no{
For $0\neq x \in \SslR$  let the {\bf $\ddd$-lead term} $\LT^\ddd(x)= C_\al \al$, where $0 \neq C_\al \in \cR$ and $\al\in \Bsl(\fS)$ is the largest (in the $\ddd$-order) basis element in the linear presentation of $x$.
  Note that $\LT^\ddd(x)$ is in $\SslR$, not in any associated graded algebra.
  }

\subsection{The associated quantum torus}\label{ssQTs}

Define the antisymmetric matrix (or map)  
\begin{align*} \mQ: \cC \times \CC\to \BZ, \ 
\mQ(a,c)=
\# \left(\begin{array}{c}\includegraphics[scale=0.17]{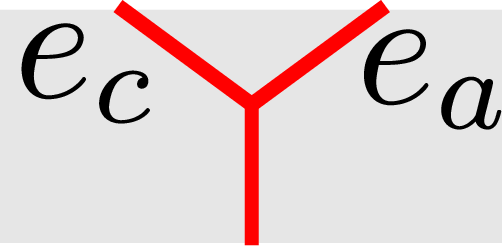}\end{array}\right)-\#\left(\begin{array}{c}\includegraphics[scale=0.17]{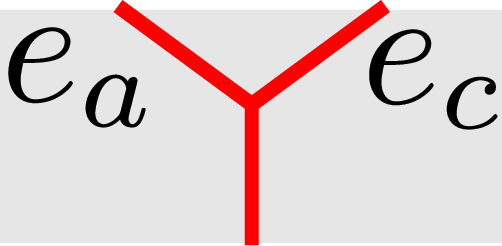}\end{array}\right).
\end{align*}
Here $e_a$ is the edge of $\Gamma$ dual to $a\in \cC$, and the right hand side is the signed number of times a half-edge of $e_a$ meets a half-edge of $e_c$ at a vertex of $\Gamma$, where the sign is $+1$ if $e_a$ is clockwise to $e_c$, and $-1$ otherwise. 

Define the {\bf symplectic double } $\tmQ$ of $\mQ$ by
$$ \tmQ =  \begin{bmatrix} \mQ &  2\Id_\cC \\ -2\Id_ \cC & 0_\cC
\end{bmatrix}$$
where $\Id_\cC$ and $0_\cC$ are respectively the $\cC\times \cC$ identity matrix and the $\cC\times \cC$ 0-matrix.

 Let $\cY_\xi(\fS;\cR)$ be the quantum torus  over $\cR$,  associated with $\tmQ$:
$$ \cY_\xi(\fS;\cR):= \cR\la y_c^{\pm 1}, u_{c}^{\pm 1}, c \in \cC \ra / ( y_a y_c = \xi^{\mQ(a,c)} y_c y_a, u_a u_c= u_c u_a, u_a y_c= \xi^{2\delta_{ac}} y_c u_a  ). $$
Here for brevity we suppress the dependence on $\cC, \Gamma$ in the notation $\cY_\xi(\fS;\cR)$.

The set of normalized monomials $\{ Y^\bm \mid \bm \in \BZ^{2r}\}$ is a free $\cR$-basis of $\YxS$. 
Here if $\bm = (\bn, \bt)\in \BZ^\cC \times \BZ^\cC$ then $Y^\bm= [y^\bn u^\bt]_\Weyl$.
The monomial subalgebra $\TxL\subset \YxS$ is the $\cR$-submodule spanned by $\{Y^\bm \mid \bm\in \Lambda\}$.

\subsection{Main result}  We formulate the main result of this section which yields that  the sliced skein algebra is always a domain and hence the sliced character variety is an irreducible affine variety
for any $\vw\in \BC^m$.

\bthm \label{thmMain2} Let $\fS=\Sigma_{g,m}$ with $r:= 3g-3+m > 0$.
 Let $(\cC, \Gamma)$ be a DT-datum described in Subsection \ref{ssN1}, and  $\Lambda=\Lambda_{\cC,\Gamma}$ be the submonoid of possible DT-coordinates.
 Let $\cR$ be a commutative domain, $\xi\in \cR$ be invertible, and $\vw\in \cR^m$.

(a) For $\bk, \bl  \in \Lambda$ we have $\la \bk, \bl\ra _{\tmQ} \in 2\BZ$, and in $\SslR$,
\be 
\mu(\bk) \mu(\bl) =  \xi^{\frac 12 \la \bk, \bl\ra _{\tmQ}  } \mu (\bk + \bl) 
\mod E_{ <_\ddd  \, (\bk + \bl) }(\SslR).  \label{eqProdTop}
\ee

(b) The $\Lambda$-filtration $(E_\bk(\SslR))_{\bk \in \Lambda}$ is compatible with the product, and its associated graded algebra is isomorphic to the monomial algebra $\bT_\xi(\tmQ,\cR;\Lambda)$:
\be 
\Gr^E(\SslR) \xar{ \cong} \bT_\xi(\tmQ,\cR;\Lambda).
\ee

(c)  The sliced skein algebra $\SslR$ is a domain.
\ethm

\def\DDD{\mathsf{D}}

\def\uTR{{{\mathsf{TR}}}}
\def\UniF{\bigsqcup _{\tau\in \cF}}
\def\tY{{\tilde \cY}}
\def\cSC{{ \check\fS_\cC}}
\def\EEE{{\tY}^\diamond}

For  $\xi=-1\in\cR=\BC$, we have $\SslR = \BC[\X^\vw(\fS)]$. Hence we have
\bcor[See \cite{Whang}] \label{rInt} The sliced character variety $\X^\vw(\fS)$ is an irreducible affine variety
for any $\vw\in \BC^m$.  
\ecor


\def\tcS{\tilde {\mathcal S}}
\def\cCtt{{\cC^{(2)}_\tau}}

\subsection{An $\BN$-filtration} \label{ssFilN}

To prove Theorem \ref{thmMain2} we use an $\BN$-filtration, which is defined now.

Let $\ddd_1: \BN^\cC \times \BZ^{ \cC} \to \BN$ be the first component of $\ddd$, meaning $\ddd_1(\bn, \bt) = \sum_{c\in \cC} \bn(c)$.
Define an $\BN$-filtration of $\SslR$ by
\begin{align}
F_k(\SslR) &:= \cR\text{-span of } \ \{ \mu(\bl) \mid \ddd_1(\bl) \le k \} , \  k\in \BN.  \label{eqFil1}
\end{align}
This is an algebra $\BN$-filtration  whose associated graded algebra is denoted by $\Gr^F(\SslR)$. For $k\in \BN$ the set 
$$\Bsl_k(\fS)= \{ \mu(\bl) \in \Bsl(\fS) \mid \ddd_1(\bl) \le k\}$$
 is a free $\cR$-basis of $F_k(\SslR)$. Hence the lead term function $\lt^F$ defined in Subsection~\ref{ssFiltr}  maps $\Bsl(\fS)$ bijectively onto a free $\cR$-basis of $\Gr^F(\SslR)$. Identifying $\al$ with $\lt^F(\al)$, we will consider $\Bsl(\fS)$ as a free $\cR$-basis of $\Gr^F(\SslS)$. For $\al, \beta \in \Bsl(\fS)$ denote by $\al * \beta$ their product in $\Gr^F(\SslR)$, while the product in $\SslR$ is denote by $\al \beta$. 

Parallel to the $\Lambda$-filtration of $\SslR$, define the following $\Lambda$-filtration of $\YxS$:
\begin{align}
E_\bk(\YxS) &:= \cR\text{-span of } \ \{ Y^\bl \mid \bl   \led \bk\}, \ \bk\in \Lambda.
\label{eqFil3}
\end{align}
This is an algebra filtration whose associated graded algebra is canonically isomorphic to, and will be identified with, the quantum torus $\YxS$.

\bthm \label{thmMain3} (a)  With the assumption of Theorem \ref{thmMain2},
there is an $\cR$-algebra homomorphism
\be 
\phi_\cR: \Gr^F(\SslR) \to \YxS
\ee
such that for $\bk\in \Lambda$ we have  
\be
 \phi_\cR(\mu(\bk)) = Y^{\bk} + E_{<_\ddd \bk   }(\YxS). \label{eqtopdeg}
\ee

(b) The map $\phi_\cR$ is injective. 
Consequently $\SslR$ is a domain.
\ethm 

Theorem \ref{thmMain3} will be proved in Subsection \ref{ssMain3}.

\subsection{Combining the faces} In this subsection, 
using $\SC = \UniF \tau$, we show how to patch the quantum traces $\utr_\tau$ together to get a global map.

We assume now the ground ring is $\cR=\RS$ defined in Subsection \ref{ssMRY}, so that we work with the universal sliced skein algebras. The reason is we want to use the reflection to simplify some calculations.
For brevity we shorten $\cY_q(\fS; \RS)$ to $\cY(\fS)$.

Let $\cSC = \UniF \ctau$. 
Identifying punctures of $\SC$ with those of $\fS$ via $\pr$, we assume
 $\cR_\SC= \RS$. Let
 \begin{align*}
 \tcS := \Ssl(\cSC)=\bigotimes_{\tau \in \cF} \Ssl(\Pb), \quad 
 \tY:=  \bigotimes_{\tau \in \cF} \cY(\tau)
 \end{align*}
 where both tensor products are over $\Zq$.  Both
 $\tcS$ and $\tY$ are $\RS$-algebras.

 Taking the tensor product of $\utr_\tau$, we have an $\RS$-algebra homomorphism
$$ \uTR= \bigotimes_{\tau \in \cF} \utr_\tau \, :  \tcS \to \tY.$$
Recall that $\cY(\ctau)$, for $\tau=\PP_j$, is the $\RR_j$-algebra of  polynomials in variables $x_a^{\pm 1}, u_a^{\pm 1}$, where $a$ runs through the set $\cCtt$ of boundary components of $\tau$, and these variables $q$-commute as described in Subsection \ref{ssQtPj}.
We identify  $\tY$ with
the $\RS$-algebra of polynomials in  variables $x^{\pm1}_a, u_a^{\pm1}, a \in \cCt$, where the variables corresponding to boundary components of a face $\tau$ are $q$-commuting by the rules of $\cY(\ctau)$, otherwise they  commute.

Let $\EEE $ be the $\RS$-subalgebra of $\tY$ spanned by monomials in which the degrees of $x_{c'}$ and $x_{c''}$ are equal for $c\in \cC$. Equivalently, $\EEE$ is the $\RS$-subalgebra
generated by  $\{ u^{\pm 1}_a, a\in \cCt\} \cup \{(x_{c'} x_{c''})^{\pm1},c\in \cC\}$. The following explains where the matrix $\mQ$ comes from.

\def\DDD{\tcS^\diamond   }
\blem \label{rIsoY}
There is a $\RS$-algebra isomorphism
\be 
\cY(\fS) \xar{\equiv} \EEE/(u_{c'} = u_{c''}, c\in \cC), \ \text{with } y_c \to [x_{c'} x_{x''}], \ u_c\to u_{c'} = u_{c''}.
\ee
\elem
\bpr This follows right away from the definition of $\cY(\fS)$. 
\epr
We will identify $\cY(\fS)$ with  $\EEE/(u_{c'} = u_{c''}, c\in \cC)$ via  Lemma \ref{rIsoY}.

 Recall that an element $\al=(\al_\tau)_{\tiF} \in \prod_{\tau \in \cF} \Bsl(\Pb)$ is matched if $|\al \cap c'| =|\al \cap c''|$ for all $c\in \cC$. 
Let $\tcS^\diamond $ be the $\RS$-submodule of $\tcS$ spanned by all boundary ordered {\em matched} elements. Theorem \ref{thmMain2}(a) implies that $\uTR(\DDD)\subset \EEE$.
Consider the composition
\be 
\Phi: \DDD  \xra{\uTR} \EEE \onto \cY(\fS).
\ee

\subsection{Definition of $\phi$} We will first define $\phi: \Gr^F(\SslS) \to \cY(\fS)$ as an $\RS$-linear map, then we show that it respects the product.

\def\Lift{{\mathsf{Lift}}}
\def\cSC{{\check\fS_\cC}}
\def\pcSC{\partial \cSC}
\def\pSC{\partial \cSC}
\def\tiF{{\tau \in \cF}}
\def\ttD{{\tilde D}}


Let $\al\in \Bsl(\fS)$. Choose a diagram $D$ in good position representing $\al$. Let $h$ be a choice of a linear order on each set $D \cap c, c\in \cC$. Let $\Lift(D, h)\in \prod_{\tau \in \cF} \Bsl(\Pb)$ be the collection of lifts of $D$ in faces, with the boundary order lifting $h$. As $\Lift(\al,h)$ is matched, we can define $\phi(D,h):=\Phi(\Lift(D,h))$.

\blem $\phi(D,h)$ depends only on $\al$ as an element of $\Bsl(\fS)$.
\elem
\bpr
Let us show $\phi(D,h)$ does not depend on $h$. We only need to show that $\phi(D,h)$ does not change if
 we exchange the $h$-order of two consecutive points on $\al \cap c$.
 The height exchange formula \eqref{eqHeight} tells us the change in height on $c'$ results in a factor $q^\epsilon$ for an $\epsilon\in \{\pm 1\}$, while the change in $c''$ is in the opposite direction and results in the factor  $q^{-\epsilon}$. Thus $\phi(D,h)$ does not change, and hence does not depend on $h$. We will drop $h$ in the notation  $\phi(D,h)$.
 
Let us show  $\phi(D)= \phi(D')$ when $D,D'$ are  good diagram representatives $\al$. Lemma~\ref{rSlides2} reduces the proof to the case when $D'$ is obtained from $D$ by a t-move or a loop-move.
 
$\bullet$ the loop-move at $c\in \cC$. For the lift this move is $\ell_{c'} \leftrightarrow \ell_{c''}$.
 By \eqref{eqLoop_utr}, we have 
 $$\uTR(\ell_{c'})= u_{c'} + u_{c'}^{-1}, \ \uTR(\ell_{c''})= u_{c''} + u_{c''}^{-1}. $$ 
 As $u_{c'}= u_{c''}$ in $\cY(\fS)$, we get that $\uTR(\ell_{c'})= \uTR(\ell_{c''})$, implying $\phi(D) = \phi(D')$.
 
$\bullet$ t-move at $c\in \cC$. We have 
 $$ \Lift(D',h) = \theta_{c'}^{\epsilon}  \theta_{c''}^{-\epsilon} ( \Lift(D,h)),$$
 where $\epsilon\in \{ \pm 1\}$, depending on the direction of the t-move. Identity \eqref{eqTwist_utr} shows that $\uTR(\Lift(D',h)) = \uTR( \Lift(D,h))$, implying $\phi(D)= \phi(D')$.
\epr
\def\YfS{\cY(\fS)}
The lemma shows we have a map $\phi: \Bsl(\fS) \to \cY(\fS)$. As $\Bsl(\fS)$ is a $\RS$-basis of $\Gr^F(\SslS)$, we extend   $\phi$  to a $\RS$-linear map, also denoted by $\phi: \Gr^F(\SslS) \to \YfS$.

In the proof of the following lemma, we use the filtration $F$ in an essential way.
\blem The  map $\phi: \Gr^F(\SslS) \to \YfS$ is an $\RS$-algebra homomorphism.
\elem

\def\Cr{\mathsf{Cr}}
\def\Res{{\mathsf {Res}}}
\def\dcup{{\vec \cup   }}
\bpr
Let $\al, \al'\in \Bsl(\fS)$. We need to show
\be
\phi(\al *\al')=\phi(\al) \phi(\al').
\label{eqAlgPhi}
\ee 
Let $D$ and $D'$ be respectively good position representatives of $\al$ and $\al'$. We can assume that $\Cr:=D\cap D'$ does not intersect any $c\in \cC$. Choose an order $h$ (respectively $h'$) on  $D \cap c$ (respectively $D'\cap  c$), for each $c\in \cC$. Using the projection $\pr$,  identify $\Lift(D) \cap \Lift(D')$ with~$\Cr$.

The product $\al\al'$ in $\SslS$ is presented by the diagram $D \dcup D'$, the union of $D$ and $D'$ with $D$ above $D'$. Similarly, the product $\Lift(D,h)\,  \Lift(D', h')$ in $\Ssl(\cSC)$ is given by the diagram $\Lift(D,h)\, \dcup\,  \Lift(D', h')$.

Assume $\sum_{c\in \cC} |\al \cap c| + |\al' \cap c| =k$. Then $\al * \al' \in \Gr^F_k(\SslS)$. 

Let $\Res= \{ \rho:\Cr \to \{ \pm 1\}\}$. 
From Identity \eqref{eq98aa},  we have, in $\SslS$, 
\begin{align}
\al\, \al' &= \sum_{\rho\in \Res} q^{\sigma(\rho)}  (D \dcup D')_\rho, \quad \sigma(\rho)=\sum _{x \in \Cr} \rho(x),
\label{eq98a}
\end{align}
where $(D \dcup D')_\rho$ is the result of $\rho(x)$-smoothing at all crossings $x$ of $(D \dcup D')$. 
If 

(*) there is a face $\tau$ such that $\tau \cap (D \dcup D')_\rho$ has a near boundary arc,\\
 then $(D \dcup D')_\rho$ is in $F_{k-1}(\SslS)$, as the total geometric intersection number with $c\in \cC$ drops. Hence, if $\Res^*$ is the set of all $\rho: \Cr\to \{\pm 1\}$ not satisfying (*), then
 \be 
\al * \al' = \sum_{\rho \in \Res^*} q^{\sigma(\rho)}  (D \dcup D')_\rho. \label{eq98}
\ee

The same resolution $\rho$ applies to $\Lift(D,h)\,  \dcup\,  \Lift(D', h')$, and the result is 0 in $\Ssl(\cSC)$ if we have (*), because of the near boundary arc relation. Hence in $\Ssl(\cSC)$ we have
\begin{align}
\Lift(D,h)\,  \Lift(D', h') &=  \sum_{\rho\in \Res^*} q^{\sigma(\rho)}   (\Lift(D,h) \dcup \Lift(D', h'))_\rho\, .\label{eq97a}
\end{align}
When $\rho\in \Res^*$, with the height order $h$ on boundary of $D$ and $h'$ on the boundary of $D'$, the lift of $(D \dcup D')_\rho $ is exactly $(\Lift(D,h) \cup \Lift(D', h'))_\rho$. Hence  \eqref{eq98} and \eqref{eq97a} imply~\eqref{eqAlgPhi}.
\epr

\subsection{Proof of Theorem \ref{thmMain3}} \label{ssMain3} 

(a) We consider two cases.

\noindent {\bf Case 1:}
Assume  $\cR=\RS$.
As $\phi$ has been defined, let us prove  \eqref{eqtopdeg}. 
 Identity \eqref{eqtopdeg} follows fairly easily from \eqref{eqHighdeg_utr} and Property (N1) of the numeration of $\cC$ given in Subsection \ref{ssN1}. Here are the details.

 Let $\Ctau\subset \cCt$ be the set of boundary components of face $\tau$. For a function $\bn: \Ctau\to \BZ$ let $\tbn: \cCt \to \BZ$ be its 0 extension (the extension taking value 0 on $\cCt\setminus \Ctau$).

Consider the following  additive functions
 \begin{align*}
& \ddd_\tau: \BN^{\Ctau} \times \BZ^{\Ctau}  \to \BZ^3, \quad \text{defined in Subsection \ref{ssFiltr6}} \\
& \tord: \BN^{\cCt} \times \BZ^{\cCt}  \to (\frac 12 \BZ)^{2r}, \quad \text{given by}\\ 
&(\bn, \bt) \to (\sum_{a\in \cCt}\! \frac {\bn(a)}{2}, \sum_{a\in \cCt}\! \bt(a), \bt(c_1') + \bt(c_1''), \dots, \bt(c_{r-1}') + \bt(c_{r-1}''), \bn(c_1') + \bn(c_1''), \dots,  \bn(c_{r-1}') + \bn(c_{r-1}'')  ).
 \end{align*}
Property (N1) is designed specifically so that  the extension from $\Ctau$ to $\cCt$ has property:
\blem \label{rOrder2}
If
$\ddd_\tau((\bn, \bt)) \le \ddd_\tau((\bn', \bt'))$ then $\tord((\tbn, \tbt)) \le  \tord((\tbn', \tbt'))$ in lexicographic order.
\elem
\bpr For $\tau= \PP_3$ or $\PP_1$ the statement is obvious. For $\PP_2$ the statement follows from (N1).
\epr
 By \eqref{eqHighdeg_utr}, we have the following top degree in face $\tau$:
$$ \utr_\tau(\al_\tau) \qeq (Y_\tau)^{\nu(\al_\tau)} + E_{ < \ddd_\tau(\nu(\al_\tau)  ) }(\cY(\Pb)).$$
Considering $\Ytau$ as a subspace of $\tY$ and using Lemma \ref{rOrder2}, we have
$$ \utr_\tau(\al_\tau) \qeq  Y^{\widetilde{\nu(\al_\tau)}} + E_{ < \tord(\widetilde{\nu(\al_\tau)}  ) }(\tY).$$
Taking the product over all faces and denoting $N(\al) = \sum _{\tau \in \cF} \widetilde{\nu(\al_\tau)}$, we get
$$ \uTR(\al) \qeq Y^{N(\al)  } +  E_{ < \tord(N(\al) ) }(\tY).$$
The projection $\EEE \onto \cY(\fS)$ sends $Y^{N(\al)  }$ to $Y^{\nu(\al)}$. Besides $\tord(N(\al))= \ddd(\nu(\al))$. Hence
$$ \phi(\al) \qeq Y^{\nu(\al)  } +  E_{ < \ddd(\nu(\al) ) }(\cY(\fS)). $$
Using the reflection invariance of Lemma \ref{rReflection}, we can replace $\qeq$ by $=$, and obtain \eqref{eqtopdeg}.

\noindent {\bf Case 2a:}
  $\cR$  contains $\xi^{1/2}$ and elements $s_i$ such that $s_i+ s_i^{-1} = \w_i$. (For example $\cR=\BC$.)
  Then $\cR$ is a $\RS$-algebra by $q^{1/2} \to \xi^{1/2}, v_i \to s_i$. Changing the ground ring, define $\phi_\cR= \phi \ot_\RS \cR$. As the basis $\Bsl(\fS)$, which is used in the definition of the $E$ and $F$ filtrations, remains a basis under changing of ground rings, we still have \eqref{eqtopdeg}.
  
  When $\cR$ does not have $\xi^{1/2}$, or $s_i$ such that $s_i+ s_i =\w_i$ we can reduce to the above case by a sequence of degree two ring extensions. The point is that in the definition of $\SslR$ we use only integral power of $\xi$ but we never use odd powers of $\xi^{1/2}$, similarly we use only the sum $ \w_i= s_i + s_i^{-1}$. Below is a sketch.

\noindent {\bf Case 2b:}
$\cR$ does not have $\xi^{1/2}$, but (a) is true for the ground ring $\cR'$ which is the extension of $\cR$ by adding $\xi^{1/2}$. 
The Galois algebra involution $\sigma: \cR' \to \cR'$  which is identity on $\cR$ and 
$\sigma(\xi^{1/2}) = - \xi^{-1/2}$ can be extended to $\cS_\xi^{\sli, \vw} (\fS; \cR')$ by setting $\sigma(\al) = \al$ for all $\al \in \Bsl(\fS)$. That the extended $\sigma$ is an algebra map follows from the fact that in the defining relation of sliced skein algebra we use integral powers of $\xi$, but not odd powers of $\xi^{1/2}$. 
Similarly we extend $\sigma$ to an algebra involution of $\cY_\xi(\fS;\cR')$. By taking the $\sigma$-invariant part of $\phi_{\cR'}$ we get the result for the ground ring $\cR$.

\noindent {\bf Case 2c:} For some $i$,
$\cR$ does not have $s_i $ such that $s_i+s_i^{-1}= \w_i$, but (a) is true for the ground ring $\cR'$ which is the extension of $\cR$ by adding $s_i$. The Galois involution is given by $\sigma(s_i) = s_i^{-1}$.
The proof is similar to the case 2b. 

The general case of $\cR$ is reduced to Case 2a by  induction using 2b and 2c.

(b) 
Let us now prove  $\phi_\cR$ is injective.
 Define the $\Lambda$-filtration of  $\Gr^F(\SslR)$ by
$$ E_\bk(\Gr^F(\SslR)) = 
\cR\text{-span of } \ \{ \al\in \Bsl(\fS) \mid \nu(\al)  \le_\ddd \bk\}, \bk \in \Lambda.  $$

Equation \eqref{eqtopdeg} implies that $\phi_\cR$ respects the $\Lambda$-filtrations of its domain and target space, and $\Gr^E(\phi_\cR)(\al) = Y^{\nu(\al)}$. Thus
 $\Gr^E(\phi_\cR)$ maps the $\cR$-basis $\Bsl(\fS)$ of
 $\Gr^F(\SslR)$ bijectively onto the basis 
$\{Y^\bk \mid \bk \in \Lambda\}$ of $\bT_\xi(\tmQ,\Lambda;\cR)$. Hence $\Gr^E(\phi_\cR)$ is a linear $\cR$-isomorphism. By Lemma~\ref{rLift1}, the map $\phi_\cR$ is injective. 

 Since $\Gr^F(\SslR)$ embeds into the domain $\bT_\xi(\tmQ,\Lambda;\cR)$, it is a domain. By Proposition~\ref{liftfacts}, $\SslR$ is a domain. This completes the proof of Theorem \ref{thmMain3}. \qed

\subsection{Proof of Theorem \ref{thmMain2}} \label{ssMain2} We are ready to complete the proof of the theorem.
\bpr[Proof of Theorem \ref{thmMain2}] 

(a) For notation brevity, we let $\cR=\RS$, as the proof for general $\cR$ is identical. By \eqref{eqtopdeg},
\begin{align*}
\phi(\mu(\bk)) & = Y^{\bk} + E_{<_ \ddd \, \bk  }(\cY(\fS)) \\
\phi(\mu(\bl)) & = Y^{\bl} + E_{<_ \ddd \,\bl  }(\cY(\fS)) \\
\phi(\mu(\bk+ \bl) & = Y^{\bk+ \bl} + E_{<_ \ddd \, (\bk+ \bl )  }(\cY(\fS))
\end{align*}
Using  $Y^\bk Y^\bl = q^{\frac 12 \la \bk, \bl\ra   } Y^{\bk+\bl}$, where $\la \bk, \bl\ra = \la \bk, \bl\ra_\tmQ$, we get
\be 
\phi\left(\mu(\bk) * \mu(\bl) - q^{\frac 12 \la \bk, \bl\ra  } \mu (\bk + \bl)\right) \ \in E_{<_ \ddd \, (\bk+ \bl )  }(\cY(\fS))
\label{eq73s}
\ee
By Lemma \ref{rLift1}, 
$$\phi^{-1} (E_{ < _\ddd\, (\bk+ \bl ) } (\cY(\fS))) = E_{ < _\ddd\, (\bk+ \bl ) } (\Gr^F(\SslS)).$$ 
Hence by taking $\phi^{-1}$ of \eqref{eq73s}, we have, in $\Gr^F(\SslS)$,
\be 
\mu(\bk) * \mu(\bl) =  q^{\frac 12 \la \bk, \bl \ra  } \mu(\bk+\bl) +   \sum c_\bj\, \mu(\bj), \quad c_\bj\in \RS, \bj) <_ \ddd\, (\bk+ \bl).
\label{eq101}
\ee 
By definition of $\Gr^F(\SslS)$, we have $\ddd_1(\bj)= \ddd_1(\bk+\bl)$. Lifting \eqref{eq101} to $\SslS$, we have
\be 
\mu(\bk) \mu(\bl) =  q^{\frac 12 \la \bk, \bl \ra  } \mu(\bk+\bl) +   \sum c_\bj\, \mu(\bj) + \sum c_\bi\, \mu(\bi), \quad c_\bi\in \RS, \ddd_1(\bi) < \ddd_1(\bk+ \bl). \label{eq101h}
\ee
As $\ddd_1$ is the first component of $\ddd$, we have that $ \ddd_1(\bi) < \ddd_1(\bk+ \bl)$ implies $\bi <_ \ddd\, (\bk+ \bl)$. Hence last two sums in \eqref{eq101h} are in $E_{ <_\ddd\, (\bk+\bl)}(\SslS)$, and we get 
 \eqref{eqProdTop}.

(b) follows right away from Identity \eqref{eqProdTop}, and (c) is Theorem \ref{thmMain3}(b). 
\epr

\def\cZ{{\mathcal Z}}
\def\SSe{{\cS^{\sli, \vw}_\xi(\fS)}}
\def\SSee{{\cS^{\sli, \vw}_\xi(\fS)^\ev}}
\def\bal{{\mathrm{bl}}}
\def\bode{{\boldsymbol{\delta}}}
\def\Bev{{\Bsl(\fS)^\ev}}
\def\BslS{{\Bsl(\fS)}}
\def\hLast{{\hat \Lambda^\bullet}}

\def\PIdeg{{\mathrm{PI\text{-}deg}}}
\def\Bast{\Bsl(\fS)^\bullet}
\def\SSx{{\SxS}}

\def\hL{{\hat \Lambda}}
\def\hLe{{\hat \Lambda^\ev}}
\def\Lev{{\Lambda^\ev}}
\def\Last{\Lambda^\bullet }
\def\LE{\mathsf{LE}}
\section{Center and dimension over center} 
\label{secCenter}

In this section we calculate the center of the sliced skein algebra at roots  of 1, show that it is an almost Azumaya $\BC$-domain, and calculate its PI-degree.

\subsection{Setting} \label{ssSetting} In this section $\fS= \Sigma_g \setminus \{ v_1, \dots, v_m\}$, with  $r:= 3g-3+m>0$. 
Also  $\cR= \BC$, $\xi\in \BC$ is a root of unity, and $\vw= (\w_1, \dots, \w_m) \in \BC^m$. Let 
 \begin{align*}
 N''=\ord(\xi), N'= \ord(\xi^2), 
 N= \ord(\zeta^4), \ve:= \zeta^{N^2}.
 \end{align*}

Recall that  the set $\Bsl(\fS)$  of isotopy classes of strongly simple diagrams is a $\BC$-basis of the sliced skein algebra $\SSe$. Let $\Bsl(\fS)^\ev\subset \BslS$ be the subset of even elements, and $\SSee$ be the subalgebra generated by $\Bev$. Recall that $A^\bullet$ is $A$ if $N'$ is odd, and $A^\bullet=A^\ev$ if $N'$ is even.

 The Frobenius map $\Phi_{\xi}: \SeS\to \SxS$ descends to a $\BC$-algebra homomorphism
$$ \Phi_{\xi}^\vw: \cS_\ve^{\sli, T_N(\vw)} (\fS) \to  \SslS, \ \text{where} \ T_N(\vw)= ( T_N(\w_1, \dots, T_N(\w_m)).
$$



\subsection{Center and PI-degree} Following is the main result of the section.

\bthm \label{thmCenter} Assume the assumptions of Subsection \ref{ssSetting}.

(a) The center $\cZ=\cZ(\SSe)$ of $\SSe$ is 
\be \cZ =  \BC\text{-span of } \ \{ \Phi_{\xi}^\vw(\al) \mid \al \in  \Bast \}= \Phi_{\xi}^\vw\left( \cS_\ve^{\sli, T_N(\vw)} (\fS)^\bullet \right ).
\label{eqZ}
\ee

Consequently for the natural projection $\pr_\xi^\vw: \SxS \onto \SSe $,
\be 
\pr_\xi^\vw (Z(\SxS)  ) = Z ( \SSe).
\label{eqCenterPr}
\ee
\no{where 
$$ \Bast = \begin{cases} \Bsl(\fS) & \text{if $N'$ is odd}\\
(\Bsl)^\ev(\fS) & \text{if $N'$ is even}.
\end{cases} 
$$
}

(b) The algebra $\SSx$ and its quotient  $\SSe$ are almost Azumaya $\BC$-domains with the same PI-degree given by
\be
\PIdeg (\SSe) = \PIdeg (\SSx) = \begin{cases} N^{r} & \text{if} \  N' \not \in 2 \BZ  \\
2^{g} N^{r} & \text{if} \  N' \in 2 \BZ.
\end{cases}
\label{eqPIdeg}
\ee

\ethm

\bpr

Using the $\BC$-basis, any  non-zero $x\in \SSe$ has a unique presentation
\be  x = \sum _{\al \in S} c_\al \al, \ 0 \neq c_\al \in \BC, S \subset \Bsl(\fS), 1\le  |S | < \infty.
\label{eqPres}
\ee

Fix a DT-datum $(\cC, \Gamma)$ of subsection \ref{ssN1}. As $\cC=\{c_1, \dots, c_r\}$, we identify $\BZ^\CC$ with $\BZ^r$. 
We have the DT-coordinate map
$ \nu: \Bsl(\fS) \xra{ \cong } \Lambda\subset \BZ^{2r}$, and let $\mu= \nu^{-1}: \Lambda\to \Bsl(\fS)$.
Let $\Lambda^\ev = \nu( \Lambda^\ev)$. Let  $\hL$ be  the $\BZ$-spans of $\Lambda$ in $\BZ^{2r}$. Recall that $\ddd$ and $\ddd_1$ were defined in Subsections \ref{ssN1} and \ref{ssFilN}. 

(a) Let $\cZ'$ be the right hand side of \eqref{eqZ}. 
Then  $\cZ'\subset \cZ$ by Theorem \ref{thmCent1}. We will prove the converse, that if $x\in \cZ$ then $x\in \cZ'$ by induction on $(d(x), k(x))$ defined as follows. With the presentation \eqref{eqPres},
let $d(x)= \max \{ \ddd_1(\al) , \al \in S\}$ and $S'= \{ \al \in S \mid \ddd_1(\al)=d(x)\}$. Let $k(x)= |S'| \ge 1$.

Let $\al\in S'$ be the largest element in  $\ddd$-order. Since $x$ is central, for any $\bl\in \Lambda$ we have  $x \mu(\bl)= \mu(\bl) x$. Hence, from \eqref{eqProdTop} we have
$$ \xi^{\frac  12 \la \nu (\al), \bl  \ra _\tmQ}  =  \xi^{\frac  12 \la \bl , \nu (\al)  \ra _\tmQ}, \quad \text{or} \ \xi^{ \la \nu (\al), \bl  \ra _\tmQ} =1 .$$
 As $\ord(\xi)=N''$, this means $\la \nu (\al), \bl  \ra _\tmQ\in N'' \BZ$ for all $\bl\in \Lambda$.  
 By Lemma \ref{rKer}(a) below,   $ \nu(\al) \in  N \Lambda^\bullet $. Then $x':= x - \Phi_{\xi, \vz}(\beta)$ is central and $(d(x'), k(x')) < (d(x), k(x))$ in the lexicographic order of $\BN^2$. By induction we have $x\in \cZ'$. Thus $Z$ is given by \eqref{eqZ}.
 
 Using the description of $Z(\SxS)$ of Theorem \ref{thmCent1}, we get \eqref{eqCenterPr}.

(b)
Let $D$ be the number on the right hand side of \eqref{eqPIdeg}. Then
  $\SSx$ is an almost Azumaya $\BC$-domain of  PI-degree  $D$ by  \cite[Theorem 6.1]{FKL2}.
Since $\SSe$ is a $\BC$-domain by Theorem \ref{thmMain2}, and is a quotient of $\SSx$, by \cite[Lemma 17.7.2(i)]{MR} it is an almost Azumaya $\BC$-domain with PI-degree $\le D$. We will prove the converse inequality 
by proving that  $\dim_\cZ(\SSe) \ge D^2$.

 By Lemma \ref{rKer}(b) we have
  $[\hL: N\hLast] = D^2$. Choose $D^2$ elements 
$\bk_1, \dots, \bk_{D^2} \in \Lambda$ representing $[\hL: N\hLast]$. We will show that $\mu(\bk_1), \dots, \mu(\bk_{D^2})$ are linearly independent over~$\cZ$. Let   $y_1, \dots, y_{D^2} \in \cZ$ such that $I:=\{ i \mid y_i \neq 0\}$ is non-empty. We want to show $\sum y_i \mu(\bk_i) \neq 0$.

For $0 \neq x\in \SSe$ with presentation \eqref{eqPres}, define its {\bf lead exponent} $\LE(x):= \nu(\al)$, where $\al$ is the maximum of $S$ in the $\ddd$-order. By \eqref{eqtopdeg} we have $\LE(xx')= \LE(x) + \LE(x')$. 
Note that polynomial $T_N$ has degree $N$. From part (a) we see that $\LE(y_i)\in N\Last$. Hence  the elements $\LE(y_i\,  \mu(\bk_i)), i \in I,$ are pairwise distinct modulo $N\hLast$. It follows that $\sum y_i \mu(\bk_i) \neq 0$. This shows $\dim_\cZ(\SSe) \ge D^2$ and completes the proof.
\epr

\blem
\label{rKer} (a) 
Suppose   $\bk\in \Lambda$ then we  have $\bk \in N \Last$ if for all $\bl\in \Lambda$,
\be  \la \bk, \bl \ra_\tmQ \in  N'' \BZ. 
 \label{eqKer}
 \ee

(b) One has $[\hL : N\hLast]= D^2$, where $D$ is the number on the right hand side of \eqref{eqPIdeg}.
\elem
\bpr The lemma, in another form, was actually proved in \cite{FKL2}. Since the formulation is different here, we give a proof. 

(a) 
As \eqref{eqKer} holds for all $\bl\in \Lambda$, it holds for $\bl\in \hL$, which, from Proposition \ref{rDTcoord2}, is
$$ \hL= (\BZ^r)^\bal \times \BZ^r, \ 
(\BZ^r)^\bal =\{\bn \in \BZ^r \mid n_i+ n_j+ n_k \in 2\BZ \ \text{for triangular }  \ c_i, c_j, c_k \}.$$

Let $\bk=(\bn, \bt)\in \Lambda$. Let $\bn \cdot \bt$ denote the dot product.
By definition, for  $(\bn', \bt') \in \hL$, 
 \be 
\la (\bn,\bt), (\bn',\bt') \ra_\tmQ=   \la \bn, \bn' \ra_\mQ +  2\bn \cdot \bt' -  2\bn' \cdot \bt. \label{eq78}
\ee
Choose $\bn'=0$ and $\bt'= \bode_i$, which has all entries 0 except that the $i$-th one is 1. Since $(\bn', \bt')\in \hL$, 
from \eqref{eq78} we get
 $2 n_i \in N'' \BZ$, which is equivalent to
  \be 
n_i \in N' \BZ .\label{eqs25}
\ee
In particular $n_i \in n \BZ$. Choose $\bn'= 2 \bode_i$ and $\bt'=0$, we get that $t_i \in N\BZ$. Thus $(\frac{\bn } {N}, \frac \bt {N})\in\BZ^{2r}$. We will show that $(\frac{\bn } {N}, \frac \bt {N})\in\Lambda$. We need only to check that if $c_i, c_j, c_k$ are triangular, then 
\be \frac 1N(n_i+ n_j + n_k )  \in 2\BZ. 
\label{eqs26}
\ee 

Assume $N'$ is odd. Then $N= N'$. Then \eqref{eqs26} follows because
$n_i+ n_j + n_k\in 2 \BZ$  and $n_i+ n_j + n_k\in N \BZ$.
The lemma is proved in this case, as $\Lambda^\bullet = \Lambda$.

Assume now  $N'$ is even. Then $N''= 2N'= 4N$. Hence \eqref{eqs26} follows from \eqref{eqs25}. Let $\beta= \mu( \frac1{N}(\bn, \bt) ) \in \Bsl(\fS)$. Let us prove that $I(\beta,\gamma)$ is even for all $\gamma\in \Bsl(\fS)$. Note that modulo 2 we have $I(\beta, \gamma) = |\beta \cap \gamma|$. Hence
from Lemma \ref{rCross} and \eqref{eqProdTop}, for any $\gamma \in \Bsl(\fS)$, modulo 2 we have
\be 
I(\beta, \gamma) \equiv \frac{1}{2} \la \nu(\beta), \nu(\gamma) \ra_\tmQ  = \frac 1{2N} \la \bk  , \nu(\gamma) \ra_\tmQ. \label{eq83z}
\ee
By \eqref{eqKer} we have $ \la \bk  , \nu(\gamma) \ra\in N'' \BZ = 4N \BZ$. 
 It follows that  $I(\beta, \gamma) $ is even, or $\beta \in \Lambda^\bullet $.
 
(b) Case 1: $N'$ is odd. Then $\Last=\Lambda$. We have
$$ [\hL: N \hLast]=  [\hL: N \hL]= N^{2r},$$
where the last identity follows since $\hL$ is a free abelian group of rank $2r$.

Case 2: $N'$ is even. Then $\Last=\Lev$. We have
$$   [\hL: N \hLe]= [\hL: N \hL]\,  [N\hL: N \hLe]= N^{2r} [N\hL: N \hLe]= N^{2r} [\hL: \hLe].$$
It remains to prove $[\hL: \hLe] = 2^{2g}$. First let us prove that the composition
$$ f: \Lambda \xra{\mu  } \BslS \xra{h} H_1(\Sigma_g,\BZ/2)$$
is a monoid homomorphism, where $h(\al)$ is the homology class of $\al$. Let $H$ be a collection of $2r$ simple closed curves representing a basis of $H_1(\Sigma_g; \BZ/2)$. It is enough to show that for each $\gamma\in H$, the function $f_\gamma:\Lambda \to \BZ/2$ defined by  $f_\gamma(\bk)= I(\mu(\bk),\gamma) \mod 2$ is additive. By \eqref{eq83z},
$$  f_\gamma(\bk)= \frac 12  \la \bk,  \nu(\gamma)\ra_\tmQ  \mod 2,  $$
which shows that $f_\gamma$ is additive. Thus $f$ is additive.

Let $\hat f: \hL \to H_1(\Sigma_g,\BZ/2)$ be the $\BZ$-linear extension of $f$. Clearly $\hat f$ is  surjective. Note that $\hLe = \ker \hat f$. It follows that $[\hL: \hLe] = | H_1(\Sigma_g;\BZ/2)|= 2^{2r}.$
\epr

\section{Smooth locus and Azumaya locus} 
\label{secthmAzu}

 In this section we show that, for $\fS=\Sigma_{g,m}$ and  a root of unity $\xi$, a sliced-smooth  character belongs to the Azumaya locus of the sliced skein algebra $\SSe$, and also the Azumaya locus of the full skein algebra $\SSx$. When  $m=0$   and  $\ord(\xi)$ is   odd, this was proved by Ganev-Jordan-Safronov \cite{GJS}. Our result includes not only all roots of 1, but also all $m \ge 0$.

\subsection{The result}

\def\XSb{\X(\fS)^\bullet}
\def\XwSb{\X^\vw(\fS)^\bullet}
\def\Sslw{\cS_\xi^{\sli, \vw}(\fS)}
\def\VwS{\VV_\xi^\vw(\fS)}
\def\VS{\VV_\xi(\fS)}
Recall that  $\VV_\xi(\fS):= \MaxSpec(\SxS)$  is the classical shadow variety  of $\SxS$. Let $N= \ord(\xi^4$), and $\XSb= \XS$ if $\ord(\xi^2)$ is odd, and $\XSb= \XS^\ev$ if $\ord(\xi^2)$ is even.  In Subsection \ref{ssCenter}, 
the dual of $\Psi$ of \eqref{eqPsi} gives an identification
\be  \VV_\xi(\fS) = \{([\rho], \vw) \in \XSb \times \BC^m \mid  T_N(\w_i) =- \tr(\rho(\fd_i))  \}
\label{eq26a}
\ee
which gives rise to a Poisson, finite morphism of degree $N^{m}$
\be  p_\xi: \VV_\xi(\fS)  \to \XS^\bu, \ ([\rho], \vw) \xra{p_\xi} [\rho].
\ee
Recall  the statrification of $\XSb$ by the sliced character varieties
\be 
\XSb = \bigsqcup_{ \vw\in \BC^m } \XwSb.
\ee
Let $\VwS$ be the 
classical shadow variety  of sliced skein algebra $\Sslw$.
\bthm 
\label{thmSliAzu} Let $\fS=\Sigma_{g,m}$ and $\xi\in \BC$ be a root of 1, and $\vw\in \BC^m$. 

(a) Any smooth point of $\VwS$ is in the Azumaya locus $\cA(\Sslw)$ of $\Sslw$.

(b) For a smooth point $x$ of $\XwSb$, the fiber $p_\xi^{-1}(x)$ is in the Azumaya locus of $\SxS$.
 In other words,
\be 
p_\xi^{-1} \left( \bigsqcup_{\vw\in \BC^m} (\XwSb)^\smooth \right) \subset  \cA(\SxS).
\ee
\ethm

\bpr (a) Since $\Psi$ is Poisson by Proposition \ref{rMaxSpec}, the identification \eqref{eq26a} is Poisson.

By \eqref{eqCenterPr},
the natural projection $\pr^\vw_\xi: \SxS \onto \Sslw$ maps the center of the domain onto the center of the codomain, and the domain and codomain have the same PI-degree. Hence by Corollary \ref{rQuotient} the embedding $\VwS \embed \VS$, considered as the
dual of the restriction of $\pr^\vw_\xi$ onto the centers induces an, sends $\cA(\Sslx  )$ into $\cA(\SxS)$. We will identify $\VwS$ as a subset of $\VS$ using this embedding. Since $\pr^\vw_\xi$ can be defined for generic $q$, Proposition \ref{rPoiMor} shows that $\VwS$ is a Poisson subvariety of $\VS$.

 Using the identification \eqref{eq26a} followed by the projection $p_\xi$, which is Poisson, we get a Poisson isomorphism
\be p_\xi^\vw :  \VwS  \xra{\cong  }  \{[\rho] \in \XSb\mid  -T_N(\w_i) = \tr(\rho(\fd_i))  \}= \X^{-T_N(\vw)  }(\fS)^\bullet. \label{eqIso5}
\ee
Since $\X^{-T_N(\vw)  }(\fS)$ is a symplectic variety with singularity
by Proposition \ref{rSliceSym}, we conclude by Proposition \ref{rPoiAzu}
\be 
 (\VwS)^\smooth \subset \cA(\SSe).  
 \label{eqSmAzu}
 \ee

 (b)  By Lemmas \ref{rPoiAzu} 
we have 
$ \cA(\SSe) \subset  \cA(\SSx). $  Combining with \eqref{eqSmAzu}, we have
\be (\VwS)^\smooth \subset \cA(\SSe).
\label{eqincl4}
\ee
Using Isomorphism \eqref{eqIso5}, we have
$$ p_\xi^{-1} \left( (\XwSb)^\smooth \right) = \bigsqcup_{\vw'\in \BC^m, T_N(\vw')= - \vw} \VV_\xi^{\vw'}(\fS)^\smooth \subset \cA(\SSe),$$
where the inclusion follows from \eqref{eqincl4}. This completes the proof of the theorem. \epr 

\no{

\blem The classical shadow variety $\VV_{\xi, \vw}(\fS)$ of $\SSe$ is Poisson isomorphic to the sliced character variety $\X^\bullet(\fS;- T_N(\vw)$, and hence is a symplectic variety with singularity.
\elem

Let us now prove the lemma.

(1) The quotient map $f_q: \Sq(\fS) \onto \cS_{q; \vw}(\fS) = \Sq(\fS)/(\fd_i= \eta_i)$ induces a surjective $\BC$-algebra homomorphism $f_\xi: \SSx \onto \SSe$. By Theorem \ref{thmCenter}, both $\SSx$ and $\SSe$ are  almost Azumaya $\BC$-domains of the same PI-degree.  By Proposition \ref{rPoiMor} $f_\xi$ restricts to a surjective, Poisson homomorphism
\be \bar f_\xi: Z(\SSx) \onto Z(\SSe)
\ee
 Its dual gives an embedding of Poisson varieties
$$ \bar f^*:\VV_{\xi, \vw} \embed \VV_\xi.$$

(2) Consider the isomorphism $\Psi$ of \eqref{eqPsi}, and the following  diagram

 \be
\begin{tikzcd}
Z(\SSx) \arrow[r, "\Psi"]
\arrow[d,two heads,"\bar f_\xi"]  
&  \BC[\X^\bullet][d_1, \dots, d_m]/ (-T_N(d_i)=  T_{\fd_i})\arrow[d,two heads,"\bar f"] \\
Z(\SSe) \arrow[r, "\Psi_\vw"] & \BC[\X^\bullet]/ (-T_N(\eta_i)=  T_{\fd_i})
\end{tikzcd}
\label{eqDiag6}
\ee
Here $\bar f$ is the quotient map sending $d_i$ to $\eta_i$, and $\Psi_\vw$ is the quotient of $\Psi$. The kernel of $\bar f_\xi$ is generated by $\fd_i-\eta_i$, while $\Psi(\fd_i)= d_i$. Hence $\Psi(\ker \bar f_\xi) = \ker \bar f$. It follows that $\Psi$ descends to a $\BC$-algebra isomorphism $\Psi_\vw$.

As $d_i$ is Casimir in the upper right algebra, the quotient map $\bar f$ is Poisson. Since $\Psi$ is a Poisson isomorphism by Proposition \ref{rMaxSpec}, and the two vertical quotient maps are Poisson, the map $\Psi_\vw$ is a Poisson isomorphism. 

Taking the dual of $\Psi_\xi$, we get a Poisson isomorphism 
$$ (\Psi_\xi)^*: \X^\bullet(\fS;- T_N(\vw)) \to \VV_{\xi;\vw}. $$
The sliced character variety $\X^\bullet(\fS;- T_N(\vw))$ is a symplectic variety with singularities by Proposition \ref{rSliceSym}.
}

When $g=0$, we have $\XSb= \XS$ since $H_1(\bfS; \BZ/2)$ is trivial.
From the smoothness in Propositions \ref{rSmooth},
 we have the following.
\bcor  \label{r0Azu}
Suppose $\rho: \pi_1(\Sigma_{0,m}) \to SL_2(\BC)$ is an irreducible representation
such that $T_N(\tr(\rho(\fd_i) ) \neq \pm 2$, then any lift of $[\rho]$ to $\VV_\xi(\fS)$ is in the Azumaya locus of  $\SSx$.
\ecor

Interestingly, we can use the Azumaya locus result to show that  certain characters, even irreducible, are not smooth in the sliced character variety.

\bcor \label{rNonSmooth}
For $m < m'$ the  embedding $\Sigma_{g,m} \embed \Sigma_{g,m}$ induces the embedding
 $f: \X(\Sigma_{g,m} )\embed \X(\Sigma_{g,m'})$. For each $[\rho]\in \X(\Sigma_{g,m})$ its image $f([\rho])$ is not sliced-smooth.

\ecor 
\bpr
Let $\xi$ be a root of 1 of odd order. The embedding $\Sigma_{g,m'} \embed \Sigma_{g,m} $ induces a projection $\cS_\xi(\Sigma_{g,m'}) \onto \cS_\xi(\Sigma_{g,m})$. From the PI-degree formula \eqref{eqPIdeg}
 we have 
$$D:= \PIdeg(\cS_\xi(\Sigma_{g,m} ))< D' : = \PIdeg(\cS_\xi(\Sigma_{g,m'})). $$
Since $\PIdeg (\cS_\xi(\Sigma_{g,m}))=D$,  
 there is a non-zero irreducible representation $\gamma$ of $\cS_\xi(\Sigma_{g,m})$ having classical shadow $[\rho]$, with dimension $\le D$. Pulling back $\gamma$ we get a non-zero irreducible representation of
  $\cS_\xi(\Sigma_{g,m'})$ whose classical shadow is in  $p_\xi^{-1}(f([ \rho]))$. Since $D < D'$, Theorem~\ref{thmSliAzu} implies
   that $ f([\rho])$ is not sliced-smooth.
\epr

\no{
\brem If $m=0$ and  $\ord(\xi)$ is odd, then Theorem \ref{thmSliAzu} was proved in \cite{GJS}. For some results concerning the abelian characters see \cite{KK}.
\erem
}

\def\Sxr{\cS_{\xi,\rho}}
\def\SxrM{\cS_{\xi,\rho}(M)}
\def\SxH{\cS_\xi(H_g)}
\def\trho{{\tilde \rho}}
\def\SnM{\cS_{-1}(M)}

\section{Geometric representations of the skein algebra}

\label{secGeoKau}

We show that the irreducible Azumaya representation of certain characters of a closed surface can be realized by a geometric action. We introduce 
the character reduced skein module of a closed 3-manifold, and  show that it has dimension 1 when the character is irreducible. We only consider roots of unity $\xi$ where   $\xi^2$ has odd order, deferring the  case where  $\ord(\xi^2)$ is even to another work, as it requires a strengthened version of \cite{FKL3}.

\def\kz{{\kappa^{(0)}}}
\no{
\subsection{The map $\kappa$ and $\kz$} Let $M$ be a connected  oriented  3-manifold. For $\ve \in \{\pm 1\}$  we defined by \eqref{eqkappa} a $\BC$-algebra isomorphism
$$ \kappa: \SeM \xra{\cong} \SnM,$$
which depends on a spin structure of $M$ when $\ve=1$.

For $\ve \in \{\pm 1, \pm \bi \}$ define a $\BC$-algebra isomorphism $\kz: \cS_\ve(M) ^0 \xra{\cong} \cS_{-1}(M) ^0$ as follows. For $\ve=\pm 1$ let $\kz$ be the restriction of $\kappa$. Clearly $\kz$ is still a $\BC$-algebra isomorphism.

For $\ve=\pm \bi $ and $\al\subset M$ be a framed link representing 0 in $H_1(M;\BZ/2)$ let 
\be
\kz(\al) = (-1)^{\# \al} \al, \ \#\al= \text{ number of components of $\al$ }.
\ee
In Appendix we show that $\kz$ is a well-defined $\BC$-algebra isomorphism.
}

\subsection{Character-reduced skein modules}  
Let $M$ be a connected  oriented  3-manifold, and let $\xi$ be a root of 1 with  $\ord(\xi^2)$ odd. We define the character-reduced skein module of $M$. We prove that if $\rho$ is an irreducible $\G$-representation of the fundamental group of a handlebody then the skein module of the handlebody reduced at the character of  $\rho$ is an irreducible representation of the skein algebra of its boundary. We use this to show that  the $\rho$-reduced skein module of a closed $3$-manifold $M$ is isomorphic to $\mathbb{C}$ when $\rho:\pi_1(M)\rightarrow \G$ is irreducible.

Let $N =\ord(\xi^4)$, and $\ve \xi^{N^2}$. This implies that $\ve\in \{\pm 1\}$. Recall that
 $\kappa: \SeM \to \SnM$ is the identity if $\ve=-1$, and $\kappa$ is the Barrett isomorphism if $\ve=1$. In the latter case $\kappa$ depends on the choice of  a spin structure on $M$.

Let $\rho\in \Hom(\pi_1(M), SL_2(\mathbb {C}))$, 
then $[\rho]$ is a point of $\X(M)$,  and hence defines
 a maximal ideal $\fm'_\rho\leq \BC[\X(M)]$. Let $\fm_\rho$ be the pull-back of $\fm'_{\rho}$ under the composition
 \be 
 \SeM \xra {\kappa } \SnM \xra{T} \BC[\X(M)] 
 \label{eqkappa1}
 \ee
 where $T$ is the mapping from Theorem \ref{charcor}.

 \no{
Explicitly, from the definition and \eqref{eqCXM},
 \be 
 \fm_\rho = \BC\text{-span} \{  (-1)^ {\# \al} \kappa (\al) -  \tr(\rho(\al)) \mid \al \ \mid  \text{framed links in}\  M \}.
 \ee
 }
 
 Recall that $\SxM$ is a $\SeM$-module, by \eqref{eqActiono}. The {\bf $\rho$-reduced skein module} is
 \be 
 \SxrM= \SxM/\fm_\rho \SxM. 
 \ee 

\no{
Now assume $\ord(\xi^2)$ is even. Let $\rho\in \Hom(\pi_1(M), SL_2)^0$, which is the subset of all representations fixed by the action of $H^1(M;\BZ/2)$. Then $[\rho]$  is  a point of $\X(M)^0$, and  determines
 a maximal ideal $\fm'_\rho\lhd \BC[\X(M)^0]$. Let $\fm_\rho$ be the pull-back of $\fm'_{\rho}$ under the composition
 \be  \label{eqkappa1}
 \SeM^0 \xra {\kz } \SnM^0\xra{T} \BC[\X(M)]^0. 

 \ee

 Consider $\SxM^0$ as a $\SeM^0$-module, see Subsection \ref{ssAction}. Define the {\bf $\rho$-reduced skein module}
 \be 
 \SxrM= \SxM^0/\fm_\rho \SxM^0. 
 \ee 

}

Given a handlebody and an irreducible representation of its fundamental group, the character of its restriction to the boundary surface is an Azumaya point. The following theorem asserts that a reduced skein module of a handlebody is an irreducible representation of the skein algebra of its boundary surface.
\begin{theorem} \label{thmAzuRep} 
Let $\xi\in \BC$ be a root of unity with  $\ord(\xi^2)$ odd, and $H$ a handlebody of genus $g>1$. Let $\Sigma= \partial H= \Sigma_g$. In the case that $\ord(\xi)$ is also odd fix a spin structure on $H$ and restrict it to $\Sigma$.
Assume $\rho:\pi_1(H) \rightarrow \G$ is an irreducible representation. 
Denote its pullback to $\pi_1(\Sigma)$ under the inclusion map by  $\bar \rho: \pi_1(\Sigma)\to SL_2(\BC)$.
The action of $\Sx(\Sigma)$ on $\Sx(H)$ descends to an irreducible action of $\Sx(\Sigma)$ on $\cS_{\xi, \rho}(H)$. The dimension of $\cS_{\xi, \rho}(H)$ is the PI-degree of $\Sx(\Sigma)$, which is
 $N^{3g-3}$. The classical shadow of this representation is the character of $\bar \rho$. In other words, $\cS_{\xi, \rho}(H)$ is the Azumaya representation of $[\bar \rho]$.
\end{theorem}

Since every closed $3$-manifold has a Heegaard splitting, a consequence of this theorem is that a  skein module of a closed $3$-manifold reduced at an irreducible representation of its fundamental group is $\mathbb{C}$.

\begin{theorem}
 \label{thmGeoKau} Let $\xi\in \BC$ be a root of unity with $\ord(\xi^2)$ odd, and $M$ a connected oriented closed 3-manifold. 
 If $\rho: \pi_1(M) \to SL_2(\BC)$ is an irreducible representation then the $\rho$-reduced skein module $\cS_{\xi, \rho}(M)$ is isomorphic to $\BC$.

\end{theorem}

\no{ 
$$ \hat \al = (-1)^ {\# \al} \kappa (\al).$$
$$\fm_\rho = \BC\text{-span} \{  \ \hat \al -  \tr(\rho(\al)) \mid \al \ \mid  \text{framed links in}\  M \}  $$

Fix a spin structure $s$ of $M$ so that we can define the map $\kappa$ by \eqref{eqkappa}. 
}

\subsection{Surface spine of $H$}

Any handlebody is homeomorphic to a cylinder over a surface. We refer to such a surface as a {\bf spine} of the handlebody.

\bpro \label{rspine}
Given a handlebody $H$  of genus $g$,
let $\rho: \pi_1(H) \to SL_2(\BC)$ be an irreducible representation.  Let $\fS=\Sigma_{0,g+1}$, with peripheral loops $\fd_i, i=1, \dots, g+1$.
There is a homeomorphism $f: \fS\times [-1,1] \to  H$, such that
\be 
\tr(\rho(f(\fd_i))) \neq \pm 2, \  i=1, \dots, g+1.
\ee

\epro
\bpr We have $\pi_1(H)=\la x_1, \dots, x_g\ra$, the free group on $x_1, \dots, x_g$, where $\fd_i$ is in  the free homotopy class of $x_i$. The product $x_{g+1}:= x_1 \dots x_g$ represents $\fd_{g+1}$.
It is known that \cite{Hensel}  automorphisms of $\pi_1(H)$ can be realized by  orientation-preserving auto-diffeomorphisms of $H$. 
Consider the following automorphisms of $\pi_1(H)$:
\begin{enumerate}
\item Permutation: $x_i \to x_{\sigma(i)}$, where $\sigma$ is a permutation of $\{1, \dots, g\}$.
\item Inverse: $x_i \to x_i^{-1}$ for one index $i$; all other $x_j$ are fixed.
\item Slide:  $x_i \to x_i x_k$ or $x_i \to x_k x_i$, for some $k\neq i$, and $x_j \to x_j$ for $j\neq i$. 
\end{enumerate}
An element $X\in SL_2(\BC)$ is
\begin{itemize}
\item {\bf elliptic} if $\tr(X) \in (-2,2)$,
\item {\bf parabolic} if $\tr(X) = \pm 2$, but $X\neq \pm \Id$,
\item {\bf loxodromic} if $\tr(X) \not\in [-2,2]$.
\end{itemize}

\blem \label{rTrbound}

(a) Assume $ x, y \in SL_2(\BC)$, where $x$ is loxodromic,  
then the set of $n\in \BZ$ such that $\tr(y x^n) \in \{ 2 , -2 \}$ is finite.


(b)  Assume $ \rho(x_1)$ is parabolic, then there is an automorphism $\theta$ of $\pi_1(H)$ so that for some $i$,  $\rho(\theta(x_i))$ is   loxodromic.

\elem
\bpr
(a) Since $x$ is loxodromic, after a conjugation we can assume
$$ x= \begin{pmatrix}  u &0  \\ 0 & u^{-1}
\end{pmatrix}, |u| \neq 1, \quad 
y = \begin{pmatrix}
a & b \\ c &d
\end{pmatrix}.$$
Notice  $\tr( y x^n) = au^n + d u^{-n}$. From this we conclude
$$ \{ \tr(yx^n) \mid n\in \BZ \} = \begin{cases} \{0\} & a=d=0 \\
\{ a u^n \} \ \text{or} \ \{ d u^{-n} \} & a\neq 0, d=0,\ \text{or } \ a= 0, d\neq 0
\\
\{ au^n + d u^{-n}\}  & a\neq 0, d\neq 0.
\end{cases}  $$
In each case, using $|u| \neq 1$, one can easily get the statement.

(b) Let $v_1$ be the unique (up to a scalar factor) eigenvector of $\rho(x_1)$. Irreducibility implies that $v_1$ is not an eigenvector of one of $\rho(x_i)$, say $\rho(x_2)$. Let $v_2$ be an eigenvector of $x_2$. In the basis $(v_1, v_2)$ of $\BC^2$, we have
$$  \rho(x_1) = \pm \begin{pmatrix}  1 &1  \\ 0 & 1
\end{pmatrix}, \ \rho(x_2)=  \begin{pmatrix}  a &0  \\ b & a^{-1} \end{pmatrix}.
   $$
   We have $b\neq 0$ since $v_1$ is not an eigenvector of $\rho(x_2)$. 
This implies that  $\tr(\rho(x_1^n x_2))= a+ a^{-1} + nb$ is not in the interval $[-2,2]$  when $|n|$ is large. Hence $\rho(x_1^n x_2)$ is loxodromic. The automorphism $ x_2 \to x_1^n x_2$ and $x_i \to x_i$ for all other $i$, will satisfy the requirement.
\epr
Choose an arbitrary identification $H = \fS \times [-1,1]$. To prove the proposition it is enough to show that there is an automorphism $\gamma$ of $\pi_1(H)$ such that 
\be \tr(\rho(\gamma(x_i)) \neq \pm2, i=1, \dots, g+1.
\label{eqtr5}
\ee

(a) Case 1: One of $\rho(x_i)$, say $\rho(x_1)$, is loxodromic. Then $\tr(\rho(x_1))\neq \pm 2$. 
By Lemma~\ref{rTrbound}(a) there is $N$ such that $\tr(x'_i) \neq \pm 2$, where $x'_i= x_i x_1^N$ for $i=2,\dots, g$. Again by  Lemma~\ref{rTrbound}(a), there is $n$ such that $\tr(x_1x '_2, \dots, x'_g x_1^n) \neq \pm 2$ and $\tr(x'_g x_1^n)\neq \pm 2$. The   automorphism $\gamma: F_g \to F_g$ defined by $x_1\to x_1, x_i \to x'_i$ for $i=2,\dots, g-1$, and $x_g \to x_g x_1^{n+N}$, will satisfy~\eqref{eqtr5}.


Case 2: One of $\rho(x_i)$ is parabolic. Lemma \ref{rTrbound}(b) reduces this case to Case 1.

Case 3: All $\rho(x_i)$'s are elliptic or $\pm \Id$. 

First, we can get rid of $\pm \Id$  by automorphisms of the slide type. 

We can assume that matrices of the form $\rho(x_{i_1} \dots x_{i_k})$, with distinct $i_1, \dots, i_k$,  are not parabolic nor loxodromic, since otherwise we are reduced to Case 2 or Case 1. 

Irreducibility implies two of them, say $\rho(x_1)$ and $\rho(x_2)$, do not commute. In particular $\rho(x_2 x_1) \neq \pm \Id$. Thus $\rho(x_2 x_1)$ is elliptic.

Let $y= x_1 \dots x_g$.  If $\rho(y)\neq \pm \Id$, then $\rho(y)$ is elliptic. Since all $\rho(x_i)$ and $\rho(y)$ are elliptic, they do not have trace equal to  $\pm 2$.

Assume $\rho(y)= \pm \Id$. Let $\gamma(x_1)= x_2 x_1$ and $\gamma(x_i)= x_i$ for $i\neq 1$. Then each $\rho(\gamma(x_i))$ is elliptic, and $\rho(\gamma(y))= \pm  \rho(x_2)$ is also elliptic. Hence we have the conclusion. 
\no{

At this stage we get that all $\rho(x_i)$ and $\rho(x_1\dots x_g)$ are elliptic. Hence their traces are not $\pm 2$.

Now assume $\rho$ is not DP. Assume $\rho(x_1)= \begin{pmatrix} u &0 \\ 0& u^{-1}
\end{pmatrix}$ with $u\neq \pm 1$. 

Each $\rho(x_i)$ can be either (i) in $P$, (ii) in $D$, or (iii) not in $D \cup P$.  Non-DP implies that there must be one, say $\rho(x_g)$, which is in (iii). For each index $i$ of group (i) replace $x_i$ by $x_i x_g$, which cannot be $\pm \id$ since $\rho(x_g)$ is not in the group $DP$. Then none of $x_2, \dots, x_g$ is in $P$. Now replace each $x_i$, with $i=2, \dots, g$ by $x_1 x_g$. Again each $x_i x_g$ is not $\pm 1$, and not parabolic nor loxodromic, hence each is elliptic.

By Lemma \ref{DP}, for each $i$ we have $\tr(x_i) \neq 0$. Since each is elliptic, we have $\tr(\rho(x_i))\neq 2$.}
\epr

\subsection{Proof of Theorem \ref{thmAzuRep} } In order to finish the proof we need the following two lemmas.

\blem \label{rDim}
If $H$ is a handlebody with boundary a surface $\Sigma$ of genus $g$ and $\rho:\pi_1(H)\rightarrow SL_2\mathbb{C}$ is irreducible then  $\dim_\BC (\cS_{\xi, \rho}(H)) = N^{3g-3}$, which is the PI-degree of $\Sx(\Sigma)$.
\elem

\bpr By Proposition \ref{rspine} we can assume that $H= \fS \times [-1,1]$, with $\fS=\Sigma_{0,g+1}$  where the traces of the images of the peripheral loops of $\fS$  under $\rho$   are not equal to $\pm 2$. 

Recall  the embedding $\Phi_\xi: \SeS \embed Z:=Z(\SxS)$ and denote its image by $Z_0$. By Theorem \ref{thmCent1}, 
 $$ Z= Z_0 [d_1, \dots, d_{g+1}]/(T_N(d_i) -\Phi_\xi(\fd_i)).$$
 From here it is easy to see that  $Z$ is a free $Z_0$-module of rank $d= N^{g+1}$, see  \cite[Theorem 3.5]{FKL2}. 
 The dual morphism
$$  \Phi_\xi^\#: \MaxSpec(Z) \onto \MaxSpec(Z_0) $$
is finite of degree $d$. Topologically $\Phi_\xi^\#$ is a branched covering of order $d$.

The representation $\rho$ determines  $\fm_{\rho}\in \MaxSpec(Z_0)$ via Equation (\ref{eqkappa1}).  Since the polynomial $T_N(z) - w$ in variable $z$ is separable for $w \neq \pm 2$, we see that  the set $(\Phi_\xi^\#)^{-1} (\fm_\rho)$ has $d$ elements,
$ (\Phi_\xi^\#)^{-1} (\fm_\rho) = \{ \fm_1, \dots, \fm_d\}.$
Each $\fm_j\in \MaxSpec(Z)$ is a maximal ideal lying over $\fm_\rho$, meaning $\fm_j \cap Z_0 = \fm_\rho$. In particular $\fm_\rho  Z \subset \fm _j$.

{\bf Claim.} We have
\be 
Z/\fm_\rho Z \cong \bigoplus_{j=1}^d Z / \fm_j.
\label{eqUnram}
\ee 
Proof of Claim.  Since $Z$ is a free $Z_0$-module of rank $d$, we have 
 $$\dim _\BC(Z/\fm_\rho Z)=  \dim _\BC(Z\ot _{Z_0} (Z_0/\fm_\rho))= \dim _\BC(Z\ot _{Z_0} \BC)=d.$$ Consider the composition
\be  Z/\fm_\rho Z \onto Z / (\fm_1 \cap \fm _2 \cap \dots \cap \fm_d) \xra{\cong} \bigoplus_{j=1}^d Z / \fm_j,
\label{eqDim}
\ee
Here the isomorphism follows from  the Chinese Remainder Theorem, which works since for $i\neq j$ we have $\fm _i + \fm _j = Z$ due to the maximality. Since the dimensions of the domain and the codomain of the map in \eqref{eqDim} are both $d$, we get \eqref{eqUnram}, proving the Claim.

Let $A= \SxS$. Since $\fm_\rho A= (\fm_\rho Z)A$, we have, using \eqref{eqUnram},
\be 
A/ \fm_\rho A = A/ ((\fm_\rho Z) A) = A\ot_{Z} (Z/\fm_\rho Z)= \bigoplus_{j=1}^d A/\fm_j A.
\label{eqAfm1}
\ee
By Corollary \ref{r0Azu}, each $\fm_j$ is in the Azumaya locus of $\SxS$, hence 
\be  A/\fm_j A \cong M_{D'}(\BC), \ D' = \text{PI-degree of} \ \Sigma_{0,g+1}= N^{g-2}. \label{eqAfm2}
\ee
Here the formula of PI-degree is given by Theorem \ref{thmCenter}. Combining \eqref{eqAfm1} and \eqref{eqAfm2}, we get $\dim_\BC(A) = d (D')^2= N^{3g-3}$.
\epr

\blem  \label{rRho}
Given a handlebody $H$ of genus $g$ with boundary $\Sigma$ and spine $\fS$, an irreducible representation $\rho:\pi_1(H)\rightarrow \G$ and its pullback $\bar \rho:\pi_1(\Sigma)\rightarrow \G$,
the following holds \be\fm_{\bar \rho}\, \SxS = \fm_{\rho}\, \SxS.\ee
\elem
\bpr
 By definition
\be \fm_\rho = \BC\text{-span} \{  \Phi_\xi( \al )- (-1)^{\# \al} \varsigma(\al)  \tr(\rho(\al)) \mid \al \  \text{framed links in}\  H \},  
\ee
where $\varsigma(\al)\in \{\pm 1\}$ is defined so that $\kappa(\al) = \varsigma(\al) \al$.  With the spin of $\Sigma=\partial H$ coming from that of $H$, the sign $\varsigma(\al)$ 
remains the same whether in $H$ or in $\Sigma$.
Hence we have $\fm_{\bar \rho}\, \SxS \subset \fm_{\rho}\, \SxS$. On the other hand, for any link $\al$ in $H$, 
since $\Phi_\xi( \al )$ is transparent, we can push it through any other link to a thickening of the boundary $\Sigma$. This shows the converse inclusion $\fm_{\bar \rho}\, \SxS \supset \fm_{\rho}\, \SxS$. \epr

By Lemma \ref{rDim} the vector $\cS_{\xi,\rho}(H)$ has dimension equal to the PI-degree of $\Sx(\Sigma)$.
By Lemma \ref{rRho}, the action of $\Sx(\Sigma)$ on $\SxS$ descends to an action of $\Sx(\Sigma)/ \fm _{\bar \rho} \Sx(\Sigma)$ on $\cS_{\xi,\rho}(H)$. As $\Sigma$ is closed and  $\bar\rho$ is irreducible, $[\bar\rho]$ is a smooth point of $\X(\Sigma)$, see \cite{Goldman}. By Theorem \ref{thmSliAzu}, $\fm_{\bar \rho}$ is an Azumaya point of $\Sx(\Sigma)$. By Theorem \ref{thmAzu}(d), the module $\cS_{\xi,\rho}(H)$ is the Azumaya representation with classical shadow $\fm_{\bar \rho}$. This completes the proof of the Theorem \ref{thmAzuRep}.\qed

\subsection{Proof of Theorem \ref{thmGeoKau}}
Choose a Heegaard splitting $M = H_1 \cup H_2$, with $H_0 := H_1 \cap H_2= \Sigma_g$. Here each of  $H_1$ and $H_2$ is a handlebody of genus $g$. For each $i=0,1,2$ let $\rho_i: \pi_1(H_i) \to SL_2(\BC)$ be the composition of $\rho$ and the embedding $H_i \embed M$. Let $V_i=\cS_{\xi}(H_i)$ and $\bar V_i= \cS_{\xi, \rho_i}(H_i)$ for $i=0,1,2$.

There is a left $V_0$-module structure on $V_2$. 
The reflection anti-involution $\omega: V_0\to V_0$ transforms the left $V_0$-module structure of $V_1$ to a right $V_0$-module structure. 

From Lemma \ref{rRho}, we get  a right 
$\bar V_0$-module for $\bar V_1$, and  a left
$\bar V_0$-module for $\bar V_2$.

\blem \label{rTensor}
With the notation from above
$ 
\SxrM \cong \bar V_1 \ot _{ \bar V_0} \bar V_2.
$
\elem

\bpr  By \cite[Proposition 2.2]{Prz},
\be  \SxM \cong V_1 \ot_{V_0} V_2. \label{eqTensor}
\ee
The center of $V_0$ acts naturally on the right hand side of \eqref{eqTensor}. By Lemma \ref{rRho}, we have $\fm_{\rho_0} V_i = \fm_{\rho_i} V_i$. Hence
\begin{align*}
\cS_{\xi,\rho}(M)  & = (V_1 \ot_{V_0} V_2)/ \fm_\rho(V_1 \ot_{V_0} V_2) = (V_1 \ot_{V_0} V_2)\ot_{\Se(H_0)}  (\Se(H_0)/ \fm_\rho)\\
&=(V_1 \ot_{V_0} (V_2\ot_{\Se(H_0)}  \Se(H_0)/ \fm_\rho)= V_1 \ot_{V_0} \bar V_2= \bar V_1 \ot_{\bar V_0} \bar V_2,
\end{align*}
which proves the lemma.
\epr

 Since  $\fm_{\rho_0}$ is Azumaya, we have $\bar V_0 \cong M_D(\BC)$, where $D=N^{3g-3}$. 
 By Theorem \ref{thmAzuRep}, $\bar V_1$ and $\bar V_2$ are the right and the left Azumaya representations
  of $V_0$ with classical shadow $\fm_{\rho_0}$, and each is isomorphic of $\BC^D$. 
  From Lemma  \ref{rTensor} we have 
$$ \SxrM \cong \BC^D \ot_{ M_D(\BC)  } \BC^D \cong \BC,$$
proving the theorem. \qed

\subsection{The skein module of a $3$-manifold as a module over $\mathcal{S}_\epsilon(M)$ }\label{sheaf}
The goal of this section is to prove that if $\xi\in \mathbb{C}$ is a root of unity with $\ord(\xi^2)$ odd   the skein module $\mathcal{S}_{\xi}(M)$ of any compact, oriented  $3$-manifold $M$ is a finite rank module over $\mathcal{S}_\epsilon(M)$.

 If $M$ is a compact  oriented three-manifold,  Morse theory guarantees the existence of a surface $\Sigma \subset M$ such that every link in $M$ can be isotoped into a regular neighborhood of $\Sigma$. If $\partial M$ is empty, the surface $\Sigma$ gives a Heegaard splitting of the manifold $M$ into two handlebodies. For manifolds with non-empty boundary the surface $\Sigma$ splits $M$ into two compression bodies, which generalizes the concept of a Heegaard splitting. In any case such a $\Sigma$ is called a Heegaard surface.

\begin{theorem} \label{boundarymodule}Suppose that $M$ is a compact connected oriented $3$-manifold. Let $\xi $ be a  root of unity with $\ord(\xi^2)$ odd. Let  $m=\frac{n}{gcd(n,2)}$ and $\epsilon=\xi ^{m^2}$. 
 The skein module $\mathcal{S}_{\xi}(M)$ is a finite rank module  over $\mathcal{S}_{\epsilon}(M)$. 
\end{theorem}
 
\proof  Let $\Sigma$ be a Heegaard surface for $M$.  Note that $M$ can be obtained from $\Sigma\times [0,1]$ by adding $2$-handles to $\partial(\Sigma\times [0,1])$  followed by  $3$-handles. Since every link in $M$ can be isotoped to be near $\Sigma$, we have $\mathcal{S}_{\xi}(M)=\mathcal{S}_\xi(\Sigma)/S$  where  $S\leq \mathcal{S}_{\xi}(M)$  is the submodule spanned by handleslides  and the map  
\be \iota_*: \mathcal{S}_\xi(\Sigma) \to \mathcal{S}_\xi(M)\ee is surjective.

The threading map makes $\mathcal{S}_{\xi}(M)=\mathcal{S}_{\xi}(\Sigma)/S$ into a module  over the algebra $\mathcal{S}_\epsilon(M)$. 
Since by Theorem \ref{thmChFM}   framed links threaded with the $m$th Chebyshev polynomial of the first kind  are transparent, this implies that the action of $\mathcal{S}_\epsilon(\Sigma)$ leaves $S$ invariant. By \cite{AF}  $\mathcal{S}_{\xi}(\Sigma)$ is a finite rank module over $\mathcal{S}_\epsilon(\Sigma)$. Therefore $\mathcal{S}_{\xi}(M)=\mathcal{S}_\xi(\Sigma)/S$ is a finite rank $\mathcal{S}_\epsilon(\Sigma)$-module.  From the naturality of threading we have the following commutative diagram:

\be \begin{CD}  \mathcal{S}_\epsilon(\Sigma) @>inc>> \mathcal{S}_\epsilon(M)  \\    @V\Phi_{\xi}VV @V\Phi_{\xi}VV\\
\mathcal{S}_\xi(\Sigma) @>inc>> \mathcal{S}_\xi(M) . \end{CD} \ee

Hence the action of $\mathcal{S}_\epsilon(\Sigma)$ on $\mathcal{S}_\xi(M)$ factors through $\mathcal{S}_\epsilon(M)$.  Therefore, $\mathcal{S}_\xi(M)$ is a finite rank module over $\mathcal{S}_\epsilon(M)$.
 
 \qed

\begin{remark}If the order of $\xi^2$ is odd  then  $\mathcal{S}_\xi (M)$ is the global sections of a sheaf defined over the affine variety underlying $\mathcal{S} _\epsilon(M)$. The statement that $\mathcal{S}_{\xi,\rho}(M)$ is isomorphic to $\mathbb{C}$ for  irreducible $ \rho:\pi_1(M)\rightarrow \G$  means that the restriction of this sheaf to the characters of irreducible representations is a line bundle.\end{remark}

\bibliographystyle{hamsalpha}
\bibliography{biblio}

 \end{document}